\documentclass[10pt]{amsart}
\usepackage[cp1251]{inputenc}
\usepackage[english]{babel}
\usepackage{amsmath}
\usepackage{amssymb}
\usepackage{amsfonts}

\usepackage[linktocpage=true, colorlinks=true, linkcolor=blue, citecolor=blue, urlcolor=blue]{hyperref}

\begin{document}
\renewcommand{\refname}{References}

\thispagestyle{empty}

\title[Expansions of Iterated Stratonovich Stochastic Integrals]
{Expansions of Iterated Stratonovich Stochastic Integrals 
of Multiplicities 1 to 4. Combained Approach Based on
Generalized Multiple and Iterated Fourier series}
\author[D.F. Kuznetsov]{Dmitriy F. Kuznetsov}
\address{Dmitriy Feliksovich Kuznetsov
\newline\hphantom{iii} Peter the Great Saint-Petersburg Polytechnic University,
\newline\hphantom{iii} Polytechnicheskaya ul., 29,
\newline\hphantom{iii} 195251, Saint-Petersburg, Russia}%
\email{sde\_kuznetsov@inbox.ru}

\maketitle {\small
\begin{quote}
\noindent{\sc Abstract.} 
The article is devoted to the expansions of iterated
Stratonovich stochastic integrals of multiplicities 1 to 4
on the base of the combined approach of 
generalized multiple and iterated Fourier series.
We consider two different parts of the expansion of 
iterated Stra\-to\-no\-vich stochastic integrals.
The mean-square convergence of the first part is proved on the base
of generalized multiple Fourier series that are converge
in the sense of norm in Hilbert space
$L_2([t, T]^k),$ $k=1,2,3,4.$ The mean-square convergence
of the second part is proved on the base of generalized iterated Fourier
series that are converge pointwise. At that, we do not use the iterated Ito
stochastic integrals as a tool of the proof and directly consider
the iterated Stratonovich stochastic integrals.
The cases of
multiple Fourier--Legendre series and 
multiple trigonometric Fourier series are considered in detail.
The considered expansions contain only one operation of the limit transition
in contrast to its existing analogues.
This property is very important for the mean-square approximation 
of iterated stochastic integrals.
The results of the article can be applied to the numerical integration 
of Ito stochastic differential equations.
\medskip
\end{quote}
}


\setlength{\baselineskip}{1.6em}

\tableofcontents

\setlength{\baselineskip}{1.2em}


\vspace{2mm}

\section{Introduction}

\vspace{5mm}

Let $(\Omega,$ ${\rm F},$ ${\sf P})$ be a complete probability space, let 
$\{{\rm F}_t, t\in[0,T]\}$ be a nondecreasing right-continous 
family of $\sigma$-algebras of ${\rm F},$
and let ${\bf f}_t$ be a standard $m$-dimensional Wiener stochastic 
process, which is
${\rm F}_t$-measurable for any $t\in[0, T].$ We assume that the components
${\bf f}_{t}^{(i)}$ $(i=1,\ldots,m)$ of this process are independent. 
Consider
an Ito stochastic differential equation (SDE) 
in the integral form

\vspace{-1mm}
\begin{equation}
\label{1.5.2}
{\bf x}_t={\bf x}_0+\int\limits_0^t {\bf a}({\bf x}_{\tau},\tau)d\tau+
\int\limits_0^t B({\bf x}_{\tau},\tau)d{\bf f}_{\tau},\ \ \
{\bf x}_0={\bf x}(0,\omega).
\end{equation}

\vspace{2mm}
\noindent
Here ${\bf x}_t$ is some $n$-dimensional stochastic process 
satisfying the equation (\ref{1.5.2}). 
The nonrandom functions ${\bf a}: \mathbb{R}^n\times[0, T]\to\mathbb{R}^n$,
$B: \mathbb{R}^n\times[0, T]\to\mathbb{R}^{n\times m}$
guarantee the existence and uniqueness up to stochastic equivalence 
of a solution
of (\ref{1.5.2}) \cite{1}. The second integral on 
the right-hand side of (\ref{1.5.2}) is 
interpreted as an Ito stochastic integral.
Let ${\bf x}_0$ be an $n$-dimensional random variable, which is 
${\rm F}_0$-measurable and 
${\sf M}\{\left|{\bf x}_0\right|^2\}<\infty$ 
(${\sf M}$ denotes a mathematical expectation).
We assume that
${\bf x}_0$ and ${\bf f}_t-{\bf f}_0$ are independent when $t>0.$

It is well known that one of the effective approaches 
to the numerical integration of 
Ito SDEs is an approach based on the Taylor--Ito and 
Taylor--Stratonovich expansions
\cite{KlPl2}-\cite{Mi3}. The most important feature of such 
expansions is a presence in them of the so-called iterated
Ito and Stratonovich stochastic integrals, which play the key 
role for solving the 
problem of numerical integration of Ito SDEs and have the 
following form

\vspace{-1mm}
\begin{equation}                    
\label{ito}
J[\psi^{(k)}]_{T,t}=\int\limits_t^T\psi_k(t_k) \ldots \int\limits_t^{t_{2}}
\psi_1(t_1) d{\bf w}_{t_1}^{(i_1)}\ldots
d{\bf w}_{t_k}^{(i_k)},
\end{equation}

\begin{equation}
\label{str}
J^{*}[\psi^{(k)}]_{T,t}=
\int\limits_t^{*T}\psi_k(t_k) \ldots \int\limits_t^{*t_{2}}
\psi_1(t_1) d{\bf w}_{t_1}^{(i_1)}\ldots
d{\bf w}_{t_k}^{(i_k)},
\end{equation}

\vspace{2mm}
\noindent
where every $\psi_l(\tau)$ $(l=1,\ldots,k)$ is a 
nonrandom function 
on $[t,T],$ ${\bf w}_{\tau}^{(i)}={\bf f}_{\tau}^{(i)}$
for $i=1,\ldots,m$ and
${\bf w}_{\tau}^{(0)}=\tau,$
and 
$$
\int\limits\ \hbox{and}\ \int\limits^{*}
$$ 

\vspace{2mm}
\noindent
denote Ito and 
Stratonovich stochastic integrals,
respectively; $i_1,\ldots,i_k = 0, 1,\ldots,m.$
In this paper we use the definition of the Stratonovich 
stochastic integral from \cite{KlPl2}.

Note that $\psi_l(\tau)\equiv 1$ $(l=1,\ldots,k)$ and
$i_1,\ldots,i_k = 0, 1,\ldots,m$ in  
\cite{KlPl2}-\cite{Mi3}. At the same time
$\psi_l(\tau)\equiv (t-\tau)^{q_l}$ ($l=1,\ldots,k$; 
$q_1,\ldots,q_k=0, 1, 2,\ldots $) and $i_1,\ldots,i_k = 1,\ldots,m$ in
\cite{3}-\cite{20}.

The construction of 
effective expansions (converging in the mean-square sense)
for collections 
of the iterated Stra\-to\-no\-vich stochastic integrals
(\ref{str}) of multiplicities 1 to 6
composes the subject of this article.

The problem of effective jointly numerical modeling 
(with respect to the mean-square convergence criterion) of iterated 
Ito and Stratonovich stochastic integrals 
(\ref{ito}) and (\ref{str}) is 
difficult from 
theoretical and computing point of view \cite{KlPl2}-\cite{rr}.
The only exception is connected with a narrow particular case, when 
$i_1=\ldots=i_k\ne 0$ and
$\psi_1(s),\ldots,\psi_k(s)\equiv \psi(s)$.
This case allows 
the investigation with using of the Ito formula 
\cite{KlPl2}-\cite{Mi3}.

Consider a brief review of the mean-square
approximation methods for the iterated stochastic 
integrals (\ref{ito}) and (\ref{str}).

Seems that iterated stochastic integrals can be approximated by multiple 
integral sums of different types \cite{Mi2}, \cite{Mi3}, \cite{Al}. 
However, this approach implies partitioning of the interval 
of integration $[t, T]$ of iterated stochastic integrals 
(the length $T-t$ of this interval is a small 
value, because it is a step of integration of numerical methods for 
Ito SDEs) and according to numerical 
experiments this additional partitioning leads to significant computational
costs \cite{7}.

In \cite{Mi2} (also see \cite{KlPl2}, \cite{KPS}, \cite{Mi3}, 
\cite{KPW}, \cite{Zapad-9}) 
Milstein G.N. proposed to expand (\ref{ito}) or (\ref{str})
into iterated series of products
of standard Gaussian random variables by representing the Brownian bridge
process as a trigonometric Fourier series with random coefficients 
(version of the so-called Karhunen--Loeve expansion).
To obtain the Milstein expansion of (\ref{str}), the truncated Fourier
expansions of components of the Wiener process ${\bf f}_s$ must be
iteratively substituted in the single integrals, and the integrals
must be calculated, starting from the innermost integral.
This is a complicated procedure that does not lead to a general
expansion of (\ref{str}) valid for an arbitrary multiplicity $k.$
For this reason, only expansions of single, double, and triple
stochastic integrals (\ref{ito}), (\ref{str}) were presented 
in \cite{KlPl2}, \cite{KPS}, \cite{KPW}, \cite{Zapad-9} ($k=1, 2, 3$)
and in \cite{Mi2}, \cite{Mi3} ($k=1, 2$) 
for the case $\psi_1(s), \psi_2(s), \psi_3(s)\equiv 1;$ 
$i_1, i_2, i_3=0,1,\ldots,m.$

Moreover, the authors of the works
\cite{KlPl2}
(Sect.~5.8, pp.~202--204), \cite{KPS} (pp.~82-84),
\cite{KPW} (pp.~438-439),  
\cite{Zapad-9} (pp.~263-264) use 
the Wong--Zakai approximation 
\cite{W-Z-1}-\cite{Watanabe} (without rigorous proof) within the frames
of the Milstein approach 
\cite{Mi2} based on the series expansion 
of the Brownian bridge process. See discussion in Sect.~7 of 
this paper for details.

Note that in \cite{rr} the method (similar to the Milstein
approach) 
of expansion
of iterated (double) Ito stochastic integrals (\ref{ito}) 
($k=2;$ $\psi_1(s), \psi_2(s) \equiv 1;$ $i_1, i_2 =1,\ldots,m$) 
based on expansion
of the Wiener process using Haar functions and 
trigonometric functions has been considered.

It is necessary to note that the approach based 
on the Karhunen--Loeve expansion \cite{Mi2} excelled 
in several times (or even in several orders) 
the methods of multiple integral sums \cite{Mi2}, \cite{Mi3}, \cite{Al}
considering computational costs in the sense 
of their diminishing.

An alternative strong approximation method was 
proposed for (\ref{str}) in \cite{3}, \cite{4} (also see
\cite{11}-\cite{16}, \cite{19}, \cite{20}-\cite{20xxz}),
where $J^{*}[\psi^{(k)}]_{T,t}$ was represented as the multiple stochastic 
integral
from the certain discontinuous nonrandom function of $k$ variables, and the 
function
was then expressed as the ge\-ne\-ra\-li\-zed 
iterated Fourier series by complete systems of 
continuously
differentiable functions that are or\-tho\-nor\-mal in the space 
$L_2([t, T])$. 
As a result,
the general iterated series expansion of products
of standard Gaussian random variables was obtained in 
\cite{3}, \cite{4} (also see
\cite{11}-\cite{16}, \cite{19}, \cite{20}-\cite{20xxz}) for (\ref{str}) with an 
arbitrary multiplicity $k.$
Hereinafter, this method is referred to as the method of 
generalized iterated Fourier series.
It was shown \cite{3}, \cite{4} (also see
\cite{11}-\cite{16}, \cite{19}, \cite{20}-\cite{20xxz})
that the method of 
generalized iterated Fourier series $(k=2)$ leads to the Milstein approach
based on the Karhunen--Loeve expansion \cite{Mi2}
in the case of trigonometric system 
of functions.

As we mentioned above, the 
Milstein approach based on the Karhunen--Loeve expansion \cite{Mi2}
and the method of generalized 
iterated Fourier series \cite{3}, \cite{4} (also see
\cite{11}-\cite{16}, \cite{19}, \cite{20}-\cite{20xxz}) lead
to iterated application of the operation of limit transition.
So, these methods may not converge in 
the mean-square sense 
to the appropriate iterated Stratonovich stochastic
integrals (\ref{str}) 
for some methods of series summation.

The mentioned problem 
(iterated application of the operation of limit transition) 
not appears in the method, which 
is proposed for (\ref{ito}) in Theorems 1, 2 (see below)
\cite{7}-\cite{19}, \cite{20}-\cite{33}.
The idea of this method is as follows: 
the iterated Ito stochastic 
integral (\ref{ito}) of multiplicity $k$ is represented as 
the multiple stochastic 
integral from the certain discontinuous nonrandom function of $k$ variables 
defined on the hypercube $[t, T]^k$, where $[t, T]$ is the interval of 
integration of the iterated Ito stochastic 
integral (\ref{ito}). Then, 
the
nonrandom function of $k$ variables is expanded in the hypercube $[t, T]^k$
into the generalized 
multiple Fourier series converging 
in the mean-square sense
in the space 
$L_2([t,T]^k)$. After a number of nontrivial transformations we obtain
(see Theorems 1, 2 below) the 
mean-square convergening expansion of 
the iterated Ito stochastic 
integral (\ref{ito})
into the multiple 
series of products
of standard  Gaussian random 
variables. The coefficients of this 
series are the coefficients of 
generalized multiple Fourier series for the mentioned nonrandom function 
of $k$ variables, which can be calculated using the explicit formula 
regardless of the multiplicity $k$ of 
the iterated Ito stochastic 
integral (\ref{ito}).
We will call this method as the method of generalized multiple
Fourier series.

\vspace{5mm}

\section{Method of Expansion of Iterated Ito Stochastic 
Integrals Based on Generalized Multiple Fourier Series}

\vspace{5mm}

Suppose that every $\psi_l(\tau)$ $(l=1,\ldots,k)$ is a
nonrandom function from the space $L_2([t, T])$.
Define the following function on the hypercube $[t, T]^k$

\vspace{-3mm}
\begin{equation}
\label{ppp}
K(t_1,\ldots,t_k)=
\begin{cases}
\psi_1(t_1)\ldots \psi_k(t_k)\ &\hbox{for}\ \ t_1<\ldots<t_k\\
~\\
~\\
0\ &\hbox{otherwise}
\end{cases},\ \ \ \ t_1,\ldots,t_k\in[t, T],\ \ \ \ k\ge 2,
\end{equation}

\vspace{2mm}
\noindent
and 
$K(t_1)\equiv\psi_1(t_1)$ for $t_1\in[t, T].$

Suppose that $\{\phi_j(x)\}_{j=0}^{\infty}$
is a complete orthonormal system of functions in the space
$L_2([t, T])$. 
The function $K(t_1,\ldots,t_k)$ belongs to the 
space $L_2([t, T]^k).$
At this situation it is well known that the generalized 
multiple Fourier series 
of $K(t_1,\ldots,t_k)\in L_2([t, T]^k)$ is converging 
to $K(t_1,\ldots,t_k)$ in the hypercube $[t, T]^k$ in 
the mean-square sense, i.e.

\vspace{-1mm}
$$
\hbox{\vtop{\offinterlineskip\halign{
\hfil#\hfil\cr
{\rm lim}\cr
$\stackrel{}{{}_{p_1,\ldots,p_k\to \infty}}$\cr
}} }\Biggl\Vert
K(t_1,\ldots,t_k)-
\sum_{j_1=0}^{p_1}\ldots \sum_{j_k=0}^{p_k}
C_{j_k\ldots j_1}\prod_{l=1}^{k} \phi_{j_l}(t_l)
\Biggr\Vert_{L_2([t,T]^k)}=0,
$$

\vspace{3mm}
\noindent
where
\begin{equation}
\label{ppppa}
C_{j_k\ldots j_1}=\int\limits_{[t,T]^k}
K(t_1,\ldots,t_k)\prod_{l=1}^{k}\phi_{j_l}(t_l)dt_1\ldots dt_k,
\end{equation}

$$
\left\Vert f\right\Vert_{L_2([t,T]^k)}=\left(\int\limits_{[t,T]^k}
f^2(t_1,\ldots,t_k)dt_1\ldots dt_k\right)^{1/2}.
$$

\vspace{6mm}

Consider the partition $\{\tau_j\}_{j=0}^N$ of $[t,T]$ such that

\begin{equation}
\label{1111}
t=\tau_0<\ldots <\tau_N=T,\ \ \
\Delta_N=
\hbox{\vtop{\offinterlineskip\halign{
\hfil#\hfil\cr
{\rm max}\cr
$\stackrel{}{{}_{0\le j\le N-1}}$\cr
}} }\Delta\tau_j\to 0\ \ \hbox{if}\ \ N\to \infty,\ \ \
\Delta\tau_j=\tau_{j+1}-\tau_j.
\end{equation}

\vspace{4mm}

{\bf Theorem 1}\ \cite{7} (2006), \cite{8}-\cite{19}, \cite{20}-\cite{33}. 
{\it Suppose that
every $\psi_l(\tau)$ $(l=1,\ldots, k)$ is a con\-ti\-nu\-ous nonrandom function on 
$[t, T]$ and
$\{\phi_j(x)\}_{j=0}^{\infty}$ is a complete orthonormal system  
of continuous func\-ti\-ons in the space $L_2([t,T]).$ Then

$$
J[\psi^{(k)}]_{T,t}\  =\ 
\hbox{\vtop{\offinterlineskip\halign{
\hfil#\hfil\cr
{\rm l.i.m.}\cr
$\stackrel{}{{}_{p_1,\ldots,p_k\to \infty}}$\cr
}} }\sum_{j_1=0}^{p_1}\ldots\sum_{j_k=0}^{p_k}
C_{j_k\ldots j_1}\Biggl(
\prod_{l=1}^k\zeta_{j_l}^{(i_l)}\ -
\Biggr.
$$

\vspace{2mm}
\begin{equation}
\label{tyyy}
-\ \Biggl.
\hbox{\vtop{\offinterlineskip\halign{
\hfil#\hfil\cr
{\rm l.i.m.}\cr
$\stackrel{}{{}_{N\to \infty}}$\cr
}} }\sum_{(l_1,\ldots,l_k)\in {\rm G}_k}
\phi_{j_{1}}(\tau_{l_1})
\Delta{\bf w}_{\tau_{l_1}}^{(i_1)}\ldots
\phi_{j_{k}}(\tau_{l_k})
\Delta{\bf w}_{\tau_{l_k}}^{(i_k)}\Biggr),
\end{equation}

\vspace{5mm}
\noindent
where $J[\psi^{(k)}]_{T,t}$ is defined by {\rm (\ref{ito}),}

\vspace{-1mm}
$$
{\rm G}_k={\rm H}_k\backslash{\rm L}_k,\ \ \
{\rm H}_k=\{(l_1,\ldots,l_k):\ l_1,\ldots,l_k=0,\ 1,\ldots,N-1\},
$$

\vspace{-1mm}
$$
{\rm L}_k=\{(l_1,\ldots,l_k):\ l_1,\ldots,l_k=0,\ 1,\ldots,N-1;\
l_g\ne l_r\ (g\ne r);\ g, r=1,\ldots,k\},
$$

\vspace{3mm}
\noindent
${\rm l.i.m.}$ is a limit in the mean-square sense$,$
$i_1,\ldots,i_k=0,1,\ldots,m,$

\vspace{-1mm}
\begin{equation}
\label{rr23}
\zeta_{j}^{(i)}=
\int\limits_t^T \phi_{j}(s) d{\bf w}_s^{(i)}
\end{equation} 

\vspace{2mm}
\noindent
are independent standard Gaussian random variables
for various
$i$ or $j$ {\rm(}if $i\ne 0${\rm),}
$C_{j_k\ldots j_1}$ is the Fourier coefficient {\rm(\ref{ppppa}),}
$\Delta{\bf w}_{\tau_{j}}^{(i)}=
{\bf w}_{\tau_{j+1}}^{(i)}-{\bf w}_{\tau_{j}}^{(i)}$
$(i=0, 1,\ldots,m),$
$\left\{\tau_{j}\right\}_{j=0}^{N}$ is a partition of
the interval $[t, T],$ which satisfies the condition {\rm (\ref{1111})}.
}

\vspace{2mm}

It was shown \cite{20xx} (Sect.~1.1.9, 1.11, 1.12), \cite{26a} (Sect.~6, 15, 16) that 
Theorem 1 is valid for convergence 
in the mean of degree $2n$ ($n\in \mathbb{N}$).
The convergence w.~p.~1 in Theorem 1
is proved in \cite{20xx} (Sect.~1.7.2), \cite{100-000-3}
for complete orthonormal systems of Legendre polynomials 
and trigonometric functions
in the space $L_2([t,T])$.
Moreover, the complete orthonormal systems of Haar and 
Rademacher--Walsh functions in $L_2([t,T])$ 
can also be applied in Theorem 1
\cite{9}-\cite{16}, \cite{19}, \cite{20}-\cite{20xxz}, \cite{26a}.
The modification of Theorem 1 for 
complete orthonormal with weigth $r(x)\ge 0$ systems
of functions in the space $L_2([t,T])$ can be found in 
\cite{20}-\cite{20xxz}, \cite{26b}.

In order to evaluate the significance of Theorem 1 for practice we will
demonstrate its transformed particular cases for 
$k=1,\ldots,6$ \cite{7}-\cite{19}, \cite{20}-\cite{33}

\begin{equation}
\label{a1}
J[\psi^{(1)}]_{T,t}
=\hbox{\vtop{\offinterlineskip\halign{
\hfil#\hfil\cr
{\rm l.i.m.}\cr
$\stackrel{}{{}_{p_1\to \infty}}$\cr
}} }\sum_{j_1=0}^{p_1}
C_{j_1}\zeta_{j_1}^{(i_1)},
\end{equation}

\vspace{4mm}
\begin{equation}
\label{a2}
J[\psi^{(2)}]_{T,t}
=\hbox{\vtop{\offinterlineskip\halign{
\hfil#\hfil\cr
{\rm l.i.m.}\cr
$\stackrel{}{{}_{p_1,p_2\to \infty}}$\cr
}} }\sum_{j_1=0}^{p_1}\sum_{j_2=0}^{p_2}
C_{j_2j_1}\Biggl(\zeta_{j_1}^{(i_1)}\zeta_{j_2}^{(i_2)}
-{\bf 1}_{\{i_1=i_2\ne 0\}}
{\bf 1}_{\{j_1=j_2\}}\Biggr),
\end{equation}

\vspace{7mm}
$$
J[\psi^{(3)}]_{T,t}=
\hbox{\vtop{\offinterlineskip\halign{
\hfil#\hfil\cr
{\rm l.i.m.}\cr
$\stackrel{}{{}_{p_1,\ldots,p_3\to \infty}}$\cr
}} }\sum_{j_1=0}^{p_1}\sum_{j_2=0}^{p_2}\sum_{j_3=0}^{p_3}
C_{j_3j_2j_1}\Biggl(
\zeta_{j_1}^{(i_1)}\zeta_{j_2}^{(i_2)}\zeta_{j_3}^{(i_3)}
-\Biggr.
$$

\begin{equation}
\label{a3}
\Biggl.-{\bf 1}_{\{i_1=i_2\ne 0\}}
{\bf 1}_{\{j_1=j_2\}}
\zeta_{j_3}^{(i_3)}
-{\bf 1}_{\{i_2=i_3\ne 0\}}
{\bf 1}_{\{j_2=j_3\}}
\zeta_{j_1}^{(i_1)}-
{\bf 1}_{\{i_1=i_3\ne 0\}}
{\bf 1}_{\{j_1=j_3\}}
\zeta_{j_2}^{(i_2)}\Biggr),
\end{equation}

\vspace{7mm}
$$
J[\psi^{(4)}]_{T,t}
=
\hbox{\vtop{\offinterlineskip\halign{
\hfil#\hfil\cr
{\rm l.i.m.}\cr
$\stackrel{}{{}_{p_1,\ldots,p_4\to \infty}}$\cr
}} }\sum_{j_1=0}^{p_1}\ldots\sum_{j_4=0}^{p_4}
C_{j_4\ldots j_1}\Biggl(
\prod_{l=1}^4\zeta_{j_l}^{(i_l)}
\Biggr.
-
$$
$$
-
{\bf 1}_{\{i_1=i_2\ne 0\}}
{\bf 1}_{\{j_1=j_2\}}
\zeta_{j_3}^{(i_3)}
\zeta_{j_4}^{(i_4)}
-
{\bf 1}_{\{i_1=i_3\ne 0\}}
{\bf 1}_{\{j_1=j_3\}}
\zeta_{j_2}^{(i_2)}
\zeta_{j_4}^{(i_4)}-
$$
$$
-
{\bf 1}_{\{i_1=i_4\ne 0\}}
{\bf 1}_{\{j_1=j_4\}}
\zeta_{j_2}^{(i_2)}
\zeta_{j_3}^{(i_3)}
-
{\bf 1}_{\{i_2=i_3\ne 0\}}
{\bf 1}_{\{j_2=j_3\}}
\zeta_{j_1}^{(i_1)}
\zeta_{j_4}^{(i_4)}-
$$
$$
-
{\bf 1}_{\{i_2=i_4\ne 0\}}
{\bf 1}_{\{j_2=j_4\}}
\zeta_{j_1}^{(i_1)}
\zeta_{j_3}^{(i_3)}
-
{\bf 1}_{\{i_3=i_4\ne 0\}}
{\bf 1}_{\{j_3=j_4\}}
\zeta_{j_1}^{(i_1)}
\zeta_{j_2}^{(i_2)}+
$$
$$
+
{\bf 1}_{\{i_1=i_2\ne 0\}}
{\bf 1}_{\{j_1=j_2\}}
{\bf 1}_{\{i_3=i_4\ne 0\}}
{\bf 1}_{\{j_3=j_4\}}
+
$$
$$
+
{\bf 1}_{\{i_1=i_3\ne 0\}}
{\bf 1}_{\{j_1=j_3\}}
{\bf 1}_{\{i_2=i_4\ne 0\}}
{\bf 1}_{\{j_2=j_4\}}+
$$
\begin{equation}
\label{a4}
+\Biggl.
{\bf 1}_{\{i_1=i_4\ne 0\}}
{\bf 1}_{\{j_1=j_4\}}
{\bf 1}_{\{i_2=i_3\ne 0\}}
{\bf 1}_{\{j_2=j_3\}}\Biggr),
\end{equation}

\vspace{9mm}

$$
J[\psi^{(5)}]_{T,t}
=\hbox{\vtop{\offinterlineskip\halign{
\hfil#\hfil\cr
{\rm l.i.m.}\cr
$\stackrel{}{{}_{p_1,\ldots,p_5\to \infty}}$\cr
}} }\sum_{j_1=0}^{p_1}\ldots\sum_{j_5=0}^{p_5}
C_{j_5\ldots j_1}\Biggl(
\prod_{l=1}^5\zeta_{j_l}^{(i_l)}
-\Biggr.
$$
$$
-
{\bf 1}_{\{i_1=i_2\ne 0\}}
{\bf 1}_{\{j_1=j_2\}}
\zeta_{j_3}^{(i_3)}
\zeta_{j_4}^{(i_4)}
\zeta_{j_5}^{(i_5)}-
{\bf 1}_{\{i_1=i_3\ne 0\}}
{\bf 1}_{\{j_1=j_3\}}
\zeta_{j_2}^{(i_2)}
\zeta_{j_4}^{(i_4)}
\zeta_{j_5}^{(i_5)}-
$$
$$
-
{\bf 1}_{\{i_1=i_4\ne 0\}}
{\bf 1}_{\{j_1=j_4\}}
\zeta_{j_2}^{(i_2)}
\zeta_{j_3}^{(i_3)}
\zeta_{j_5}^{(i_5)}-
{\bf 1}_{\{i_1=i_5\ne 0\}}
{\bf 1}_{\{j_1=j_5\}}
\zeta_{j_2}^{(i_2)}
\zeta_{j_3}^{(i_3)}
\zeta_{j_4}^{(i_4)}-
$$
$$
-
{\bf 1}_{\{i_2=i_3\ne 0\}}
{\bf 1}_{\{j_2=j_3\}}
\zeta_{j_1}^{(i_1)}
\zeta_{j_4}^{(i_4)}
\zeta_{j_5}^{(i_5)}-
{\bf 1}_{\{i_2=i_4\ne 0\}}
{\bf 1}_{\{j_2=j_4\}}
\zeta_{j_1}^{(i_1)}
\zeta_{j_3}^{(i_3)}
\zeta_{j_5}^{(i_5)}-
$$
$$
-
{\bf 1}_{\{i_2=i_5\ne 0\}}
{\bf 1}_{\{j_2=j_5\}}
\zeta_{j_1}^{(i_1)}
\zeta_{j_3}^{(i_3)}
\zeta_{j_4}^{(i_4)}
-{\bf 1}_{\{i_3=i_4\ne 0\}}
{\bf 1}_{\{j_3=j_4\}}
\zeta_{j_1}^{(i_1)}
\zeta_{j_2}^{(i_2)}
\zeta_{j_5}^{(i_5)}-
$$
$$
-
{\bf 1}_{\{i_3=i_5\ne 0\}}
{\bf 1}_{\{j_3=j_5\}}
\zeta_{j_1}^{(i_1)}
\zeta_{j_2}^{(i_2)}
\zeta_{j_4}^{(i_4)}
-{\bf 1}_{\{i_4=i_5\ne 0\}}
{\bf 1}_{\{j_4=j_5\}}
\zeta_{j_1}^{(i_1)}
\zeta_{j_2}^{(i_2)}
\zeta_{j_3}^{(i_3)}+
$$
$$
+
{\bf 1}_{\{i_1=i_2\ne 0\}}
{\bf 1}_{\{j_1=j_2\}}
{\bf 1}_{\{i_3=i_4\ne 0\}}
{\bf 1}_{\{j_3=j_4\}}\zeta_{j_5}^{(i_5)}+
{\bf 1}_{\{i_1=i_2\ne 0\}}
{\bf 1}_{\{j_1=j_2\}}
{\bf 1}_{\{i_3=i_5\ne 0\}}
{\bf 1}_{\{j_3=j_5\}}\zeta_{j_4}^{(i_4)}+
$$
$$
+
{\bf 1}_{\{i_1=i_2\ne 0\}}
{\bf 1}_{\{j_1=j_2\}}
{\bf 1}_{\{i_4=i_5\ne 0\}}
{\bf 1}_{\{j_4=j_5\}}\zeta_{j_3}^{(i_3)}+
{\bf 1}_{\{i_1=i_3\ne 0\}}
{\bf 1}_{\{j_1=j_3\}}
{\bf 1}_{\{i_2=i_4\ne 0\}}
{\bf 1}_{\{j_2=j_4\}}\zeta_{j_5}^{(i_5)}+
$$
$$
+
{\bf 1}_{\{i_1=i_3\ne 0\}}
{\bf 1}_{\{j_1=j_3\}}
{\bf 1}_{\{i_2=i_5\ne 0\}}
{\bf 1}_{\{j_2=j_5\}}\zeta_{j_4}^{(i_4)}+
{\bf 1}_{\{i_1=i_3\ne 0\}}
{\bf 1}_{\{j_1=j_3\}}
{\bf 1}_{\{i_4=i_5\ne 0\}}
{\bf 1}_{\{j_4=j_5\}}\zeta_{j_2}^{(i_2)}+
$$
$$
+
{\bf 1}_{\{i_1=i_4\ne 0\}}
{\bf 1}_{\{j_1=j_4\}}
{\bf 1}_{\{i_2=i_3\ne 0\}}
{\bf 1}_{\{j_2=j_3\}}\zeta_{j_5}^{(i_5)}+
{\bf 1}_{\{i_1=i_4\ne 0\}}
{\bf 1}_{\{j_1=j_4\}}
{\bf 1}_{\{i_2=i_5\ne 0\}}
{\bf 1}_{\{j_2=j_5\}}\zeta_{j_3}^{(i_3)}+
$$
$$
+
{\bf 1}_{\{i_1=i_4\ne 0\}}
{\bf 1}_{\{j_1=j_4\}}
{\bf 1}_{\{i_3=i_5\ne 0\}}
{\bf 1}_{\{j_3=j_5\}}\zeta_{j_2}^{(i_2)}+
{\bf 1}_{\{i_1=i_5\ne 0\}}
{\bf 1}_{\{j_1=j_5\}}
{\bf 1}_{\{i_2=i_3\ne 0\}}
{\bf 1}_{\{j_2=j_3\}}\zeta_{j_4}^{(i_4)}+
$$
$$
+
{\bf 1}_{\{i_1=i_5\ne 0\}}
{\bf 1}_{\{j_1=j_5\}}
{\bf 1}_{\{i_2=i_4\ne 0\}}
{\bf 1}_{\{j_2=j_4\}}\zeta_{j_3}^{(i_3)}+
{\bf 1}_{\{i_1=i_5\ne 0\}}
{\bf 1}_{\{j_1=j_5\}}
{\bf 1}_{\{i_3=i_4\ne 0\}}
{\bf 1}_{\{j_3=j_4\}}\zeta_{j_2}^{(i_2)}+
$$
$$
+
{\bf 1}_{\{i_2=i_3\ne 0\}}
{\bf 1}_{\{j_2=j_3\}}
{\bf 1}_{\{i_4=i_5\ne 0\}}
{\bf 1}_{\{j_4=j_5\}}\zeta_{j_1}^{(i_1)}+
{\bf 1}_{\{i_2=i_4\ne 0\}}
{\bf 1}_{\{j_2=j_4\}}
{\bf 1}_{\{i_3=i_5\ne 0\}}
{\bf 1}_{\{j_3=j_5\}}\zeta_{j_1}^{(i_1)}+
$$
\begin{equation}
\label{a5}
+\Biggl.
{\bf 1}_{\{i_2=i_5\ne 0\}}
{\bf 1}_{\{j_2=j_5\}}
{\bf 1}_{\{i_3=i_4\ne 0\}}
{\bf 1}_{\{j_3=j_4\}}\zeta_{j_1}^{(i_1)}\Biggr),
\end{equation}

\vspace{9mm}

$$
J[\psi^{(6)}]_{T,t}
=\hbox{\vtop{\offinterlineskip\halign{
\hfil#\hfil\cr
{\rm l.i.m.}\cr
$\stackrel{}{{}_{p_1,\ldots,p_6\to \infty}}$\cr
}} }\sum_{j_1=0}^{p_1}\ldots\sum_{j_6=0}^{p_6}
C_{j_6\ldots j_1}\Biggl(
\prod_{l=1}^6
\zeta_{j_l}^{(i_l)}
-\Biggr.
$$
$$
-
{\bf 1}_{\{i_1=i_6\ne 0\}}
{\bf 1}_{\{j_1=j_6\}}
\zeta_{j_2}^{(i_2)}
\zeta_{j_3}^{(i_3)}
\zeta_{j_4}^{(i_4)}
\zeta_{j_5}^{(i_5)}-
{\bf 1}_{\{i_2=i_6\ne 0\}}
{\bf 1}_{\{j_2=j_6\}}
\zeta_{j_1}^{(i_1)}
\zeta_{j_3}^{(i_3)}
\zeta_{j_4}^{(i_4)}
\zeta_{j_5}^{(i_5)}-
$$
$$
-
{\bf 1}_{\{i_3=i_6\ne 0\}}
{\bf 1}_{\{j_3=j_6\}}
\zeta_{j_1}^{(i_1)}
\zeta_{j_2}^{(i_2)}
\zeta_{j_4}^{(i_4)}
\zeta_{j_5}^{(i_5)}-
{\bf 1}_{\{i_4=i_6\ne 0\}}
{\bf 1}_{\{j_4=j_6\}}
\zeta_{j_1}^{(i_1)}
\zeta_{j_2}^{(i_2)}
\zeta_{j_3}^{(i_3)}
\zeta_{j_5}^{(i_5)}-
$$
$$
-
{\bf 1}_{\{i_5=i_6\ne 0\}}
{\bf 1}_{\{j_5=j_6\}}
\zeta_{j_1}^{(i_1)}
\zeta_{j_2}^{(i_2)}
\zeta_{j_3}^{(i_3)}
\zeta_{j_4}^{(i_4)}-
{\bf 1}_{\{i_1=i_2\ne 0\}}
{\bf 1}_{\{j_1=j_2\}}
\zeta_{j_3}^{(i_3)}
\zeta_{j_4}^{(i_4)}
\zeta_{j_5}^{(i_5)}
\zeta_{j_6}^{(i_6)}-
$$
$$
-
{\bf 1}_{\{i_1=i_3\ne 0\}}
{\bf 1}_{\{j_1=j_3\}}
\zeta_{j_2}^{(i_2)}
\zeta_{j_4}^{(i_4)}
\zeta_{j_5}^{(i_5)}
\zeta_{j_6}^{(i_6)}-
{\bf 1}_{\{i_1=i_4\ne 0\}}
{\bf 1}_{\{j_1=j_4\}}
\zeta_{j_2}^{(i_2)}
\zeta_{j_3}^{(i_3)}
\zeta_{j_5}^{(i_5)}
\zeta_{j_6}^{(i_6)}-
$$
$$
-
{\bf 1}_{\{i_1=i_5\ne 0\}}
{\bf 1}_{\{j_1=j_5\}}
\zeta_{j_2}^{(i_2)}
\zeta_{j_3}^{(i_3)}
\zeta_{j_4}^{(i_4)}
\zeta_{j_6}^{(i_6)}-
{\bf 1}_{\{i_2=i_3\ne 0\}}
{\bf 1}_{\{j_2=j_3\}}
\zeta_{j_1}^{(i_1)}
\zeta_{j_4}^{(i_4)}
\zeta_{j_5}^{(i_5)}
\zeta_{j_6}^{(i_6)}-
$$
$$
-
{\bf 1}_{\{i_2=i_4\ne 0\}}
{\bf 1}_{\{j_2=j_4\}}
\zeta_{j_1}^{(i_1)}
\zeta_{j_3}^{(i_3)}
\zeta_{j_5}^{(i_5)}
\zeta_{j_6}^{(i_6)}-
{\bf 1}_{\{i_2=i_5\ne 0\}}
{\bf 1}_{\{j_2=j_5\}}
\zeta_{j_1}^{(i_1)}
\zeta_{j_3}^{(i_3)}
\zeta_{j_4}^{(i_4)}
\zeta_{j_6}^{(i_6)}-
$$
$$
-
{\bf 1}_{\{i_3=i_4\ne 0\}}
{\bf 1}_{\{j_3=j_4\}}
\zeta_{j_1}^{(i_1)}
\zeta_{j_2}^{(i_2)}
\zeta_{j_5}^{(i_5)}
\zeta_{j_6}^{(i_6)}-
{\bf 1}_{\{i_3=i_5\ne 0\}}
{\bf 1}_{\{j_3=j_5\}}
\zeta_{j_1}^{(i_1)}
\zeta_{j_2}^{(i_2)}
\zeta_{j_4}^{(i_4)}
\zeta_{j_6}^{(i_6)}-
$$
$$
-
{\bf 1}_{\{i_4=i_5\ne 0\}}
{\bf 1}_{\{j_4=j_5\}}
\zeta_{j_1}^{(i_1)}
\zeta_{j_2}^{(i_2)}
\zeta_{j_3}^{(i_3)}
\zeta_{j_6}^{(i_6)}+
$$
$$
+
{\bf 1}_{\{i_1=i_2\ne 0\}}
{\bf 1}_{\{j_1=j_2\}}
{\bf 1}_{\{i_3=i_4\ne 0\}}
{\bf 1}_{\{j_3=j_4\}}
\zeta_{j_5}^{(i_5)}
\zeta_{j_6}^{(i_6)}+
{\bf 1}_{\{i_1=i_2\ne 0\}}
{\bf 1}_{\{j_1=j_2\}}
{\bf 1}_{\{i_3=i_5\ne 0\}}
{\bf 1}_{\{j_3=j_5\}}
\zeta_{j_4}^{(i_4)}
\zeta_{j_6}^{(i_6)}+
$$
$$
+
{\bf 1}_{\{i_1=i_2\ne 0\}}
{\bf 1}_{\{j_1=j_2\}}
{\bf 1}_{\{i_4=i_5\ne 0\}}
{\bf 1}_{\{j_4=j_5\}}
\zeta_{j_3}^{(i_3)}
\zeta_{j_6}^{(i_6)}
+
{\bf 1}_{\{i_1=i_3\ne 0\}}
{\bf 1}_{\{j_1=j_3\}}
{\bf 1}_{\{i_2=i_4\ne 0\}}
{\bf 1}_{\{j_2=j_4\}}
\zeta_{j_5}^{(i_5)}
\zeta_{j_6}^{(i_6)}+
$$
$$
+
{\bf 1}_{\{i_1=i_3\ne 0\}}
{\bf 1}_{\{j_1=j_3\}}
{\bf 1}_{\{i_2=i_5\ne 0\}}
{\bf 1}_{\{j_2=j_5\}}
\zeta_{j_4}^{(i_4)}
\zeta_{j_6}^{(i_6)}
+{\bf 1}_{\{i_1=i_3\ne 0\}}
{\bf 1}_{\{j_1=j_3\}}
{\bf 1}_{\{i_4=i_5\ne 0\}}
{\bf 1}_{\{j_4=j_5\}}
\zeta_{j_2}^{(i_2)}
\zeta_{j_6}^{(i_6)}+
$$
$$
+
{\bf 1}_{\{i_1=i_4\ne 0\}}
{\bf 1}_{\{j_1=j_4\}}
{\bf 1}_{\{i_2=i_3\ne 0\}}
{\bf 1}_{\{j_2=j_3\}}
\zeta_{j_5}^{(i_5)}
\zeta_{j_6}^{(i_6)}
+
{\bf 1}_{\{i_1=i_4\ne 0\}}
{\bf 1}_{\{j_1=j_4\}}
{\bf 1}_{\{i_2=i_5\ne 0\}}
{\bf 1}_{\{j_2=j_5\}}
\zeta_{j_3}^{(i_3)}
\zeta_{j_6}^{(i_6)}+
$$
$$
+
{\bf 1}_{\{i_1=i_4\ne 0\}}
{\bf 1}_{\{j_1=j_4\}}
{\bf 1}_{\{i_3=i_5\ne 0\}}
{\bf 1}_{\{j_3=j_5\}}
\zeta_{j_2}^{(i_2)}
\zeta_{j_6}^{(i_6)}
+
{\bf 1}_{\{i_1=i_5\ne 0\}}
{\bf 1}_{\{j_1=j_5\}}
{\bf 1}_{\{i_2=i_3\ne 0\}}
{\bf 1}_{\{j_2=j_3\}}
\zeta_{j_4}^{(i_4)}
\zeta_{j_6}^{(i_6)}+
$$
$$
+
{\bf 1}_{\{i_1=i_5\ne 0\}}
{\bf 1}_{\{j_1=j_5\}}
{\bf 1}_{\{i_2=i_4\ne 0\}}
{\bf 1}_{\{j_2=j_4\}}
\zeta_{j_3}^{(i_3)}
\zeta_{j_6}^{(i_6)}
+
{\bf 1}_{\{i_1=i_5\ne 0\}}
{\bf 1}_{\{j_1=j_5\}}
{\bf 1}_{\{i_3=i_4\ne 0\}}
{\bf 1}_{\{j_3=j_4\}}
\zeta_{j_2}^{(i_2)}
\zeta_{j_6}^{(i_6)}+
$$
$$
+
{\bf 1}_{\{i_2=i_3\ne 0\}}
{\bf 1}_{\{j_2=j_3\}}
{\bf 1}_{\{i_4=i_5\ne 0\}}
{\bf 1}_{\{j_4=j_5\}}
\zeta_{j_1}^{(i_1)}
\zeta_{j_6}^{(i_6)}
+
{\bf 1}_{\{i_2=i_4\ne 0\}}
{\bf 1}_{\{j_2=j_4\}}
{\bf 1}_{\{i_3=i_5\ne 0\}}
{\bf 1}_{\{j_3=j_5\}}
\zeta_{j_1}^{(i_1)}
\zeta_{j_6}^{(i_6)}+
$$
$$
+
{\bf 1}_{\{i_2=i_5\ne 0\}}
{\bf 1}_{\{j_2=j_5\}}
{\bf 1}_{\{i_3=i_4\ne 0\}}
{\bf 1}_{\{j_3=j_4\}}
\zeta_{j_1}^{(i_1)}
\zeta_{j_6}^{(i_6)}
+
{\bf 1}_{\{i_6=i_1\ne 0\}}
{\bf 1}_{\{j_6=j_1\}}
{\bf 1}_{\{i_3=i_4\ne 0\}}
{\bf 1}_{\{j_3=j_4\}}
\zeta_{j_2}^{(i_2)}
\zeta_{j_5}^{(i_5)}+
$$
$$
+
{\bf 1}_{\{i_6=i_1\ne 0\}}
{\bf 1}_{\{j_6=j_1\}}
{\bf 1}_{\{i_3=i_5\ne 0\}}
{\bf 1}_{\{j_3=j_5\}}
\zeta_{j_2}^{(i_2)}
\zeta_{j_4}^{(i_4)}
+
{\bf 1}_{\{i_6=i_1\ne 0\}}
{\bf 1}_{\{j_6=j_1\}}
{\bf 1}_{\{i_2=i_5\ne 0\}}
{\bf 1}_{\{j_2=j_5\}}
\zeta_{j_3}^{(i_3)}
\zeta_{j_4}^{(i_4)}+
$$
$$
+
{\bf 1}_{\{i_6=i_1\ne 0\}}
{\bf 1}_{\{j_6=j_1\}}
{\bf 1}_{\{i_2=i_4\ne 0\}}
{\bf 1}_{\{j_2=j_4\}}
\zeta_{j_3}^{(i_3)}
\zeta_{j_5}^{(i_5)}
+
{\bf 1}_{\{i_6=i_1\ne 0\}}
{\bf 1}_{\{j_6=j_1\}}
{\bf 1}_{\{i_4=i_5\ne 0\}}
{\bf 1}_{\{j_4=j_5\}}
\zeta_{j_2}^{(i_2)}
\zeta_{j_3}^{(i_3)}+
$$
$$
+
{\bf 1}_{\{i_6=i_1\ne 0\}}
{\bf 1}_{\{j_6=j_1\}}
{\bf 1}_{\{i_2=i_3\ne 0\}}
{\bf 1}_{\{j_2=j_3\}}
\zeta_{j_4}^{(i_4)}
\zeta_{j_5}^{(i_5)}
+
{\bf 1}_{\{i_6=i_2\ne 0\}}
{\bf 1}_{\{j_6=j_2\}}
{\bf 1}_{\{i_3=i_5\ne 0\}}
{\bf 1}_{\{j_3=j_5\}}
\zeta_{j_1}^{(i_1)}
\zeta_{j_4}^{(i_4)}+
$$
$$
+
{\bf 1}_{\{i_6=i_2\ne 0\}}
{\bf 1}_{\{j_6=j_2\}}
{\bf 1}_{\{i_4=i_5\ne 0\}}
{\bf 1}_{\{j_4=j_5\}}
\zeta_{j_1}^{(i_1)}
\zeta_{j_3}^{(i_3)}
+
{\bf 1}_{\{i_6=i_2\ne 0\}}
{\bf 1}_{\{j_6=j_2\}}
{\bf 1}_{\{i_3=i_4\ne 0\}}
{\bf 1}_{\{j_3=j_4\}}
\zeta_{j_1}^{(i_1)}
\zeta_{j_5}^{(i_5)}+
$$
$$
+
{\bf 1}_{\{i_6=i_2\ne 0\}}
{\bf 1}_{\{j_6=j_2\}}
{\bf 1}_{\{i_1=i_5\ne 0\}}
{\bf 1}_{\{j_1=j_5\}}
\zeta_{j_3}^{(i_3)}
\zeta_{j_4}^{(i_4)}
+
{\bf 1}_{\{i_6=i_2\ne 0\}}
{\bf 1}_{\{j_6=j_2\}}
{\bf 1}_{\{i_1=i_4\ne 0\}}
{\bf 1}_{\{j_1=j_4\}}
\zeta_{j_3}^{(i_3)}
\zeta_{j_5}^{(i_5)}+
$$
$$
+
{\bf 1}_{\{i_6=i_2\ne 0\}}
{\bf 1}_{\{j_6=j_2\}}
{\bf 1}_{\{i_1=i_3\ne 0\}}
{\bf 1}_{\{j_1=j_3\}}
\zeta_{j_4}^{(i_4)}
\zeta_{j_5}^{(i_5)}
+
{\bf 1}_{\{i_6=i_3\ne 0\}}
{\bf 1}_{\{j_6=j_3\}}
{\bf 1}_{\{i_2=i_5\ne 0\}}
{\bf 1}_{\{j_2=j_5\}}
\zeta_{j_1}^{(i_1)}
\zeta_{j_4}^{(i_4)}+
$$
$$
+
{\bf 1}_{\{i_6=i_3\ne 0\}}
{\bf 1}_{\{j_6=j_3\}}
{\bf 1}_{\{i_4=i_5\ne 0\}}
{\bf 1}_{\{j_4=j_5\}}
\zeta_{j_1}^{(i_1)}
\zeta_{j_2}^{(i_2)}
+
{\bf 1}_{\{i_6=i_3\ne 0\}}
{\bf 1}_{\{j_6=j_3\}}
{\bf 1}_{\{i_2=i_4\ne 0\}}
{\bf 1}_{\{j_2=j_4\}}
\zeta_{j_1}^{(i_1)}
\zeta_{j_5}^{(i_5)}+
$$
$$
+
{\bf 1}_{\{i_6=i_3\ne 0\}}
{\bf 1}_{\{j_6=j_3\}}
{\bf 1}_{\{i_1=i_5\ne 0\}}
{\bf 1}_{\{j_1=j_5\}}
\zeta_{j_2}^{(i_2)}
\zeta_{j_4}^{(i_4)}
+
{\bf 1}_{\{i_6=i_3\ne 0\}}
{\bf 1}_{\{j_6=j_3\}}
{\bf 1}_{\{i_1=i_4\ne 0\}}
{\bf 1}_{\{j_1=j_4\}}
\zeta_{j_2}^{(i_2)}
\zeta_{j_5}^{(i_5)}+
$$
$$
+
{\bf 1}_{\{i_6=i_3\ne 0\}}
{\bf 1}_{\{j_6=j_3\}}
{\bf 1}_{\{i_1=i_2\ne 0\}}
{\bf 1}_{\{j_1=j_2\}}
\zeta_{j_4}^{(i_4)}
\zeta_{j_5}^{(i_5)}
+
{\bf 1}_{\{i_6=i_4\ne 0\}}
{\bf 1}_{\{j_6=j_4\}}
{\bf 1}_{\{i_3=i_5\ne 0\}}
{\bf 1}_{\{j_3=j_5\}}
\zeta_{j_1}^{(i_1)}
\zeta_{j_2}^{(i_2)}+
$$
$$
+
{\bf 1}_{\{i_6=i_4\ne 0\}}
{\bf 1}_{\{j_6=j_4\}}
{\bf 1}_{\{i_2=i_5\ne 0\}}
{\bf 1}_{\{j_2=j_5\}}
\zeta_{j_1}^{(i_1)}
\zeta_{j_3}^{(i_3)}
+
{\bf 1}_{\{i_6=i_4\ne 0\}}
{\bf 1}_{\{j_6=j_4\}}
{\bf 1}_{\{i_2=i_3\ne 0\}}
{\bf 1}_{\{j_2=j_3\}}
\zeta_{j_1}^{(i_1)}
\zeta_{j_5}^{(i_5)}+
$$
$$
+
{\bf 1}_{\{i_6=i_4\ne 0\}}
{\bf 1}_{\{j_6=j_4\}}
{\bf 1}_{\{i_1=i_5\ne 0\}}
{\bf 1}_{\{j_1=j_5\}}
\zeta_{j_2}^{(i_2)}
\zeta_{j_3}^{(i_3)}
+
{\bf 1}_{\{i_6=i_4\ne 0\}}
{\bf 1}_{\{j_6=j_4\}}
{\bf 1}_{\{i_1=i_3\ne 0\}}
{\bf 1}_{\{j_1=j_3\}}
\zeta_{j_2}^{(i_2)}
\zeta_{j_5}^{(i_5)}+
$$
$$
+
{\bf 1}_{\{i_6=i_4\ne 0\}}
{\bf 1}_{\{j_6=j_4\}}
{\bf 1}_{\{i_1=i_2\ne 0\}}
{\bf 1}_{\{j_1=j_2\}}
\zeta_{j_3}^{(i_3)}
\zeta_{j_5}^{(i_5)}
+
{\bf 1}_{\{i_6=i_5\ne 0\}}
{\bf 1}_{\{j_6=j_5\}}
{\bf 1}_{\{i_3=i_4\ne 0\}}
{\bf 1}_{\{j_3=j_4\}}
\zeta_{j_1}^{(i_1)}
\zeta_{j_2}^{(i_2)}+
$$
$$
+
{\bf 1}_{\{i_6=i_5\ne 0\}}
{\bf 1}_{\{j_6=j_5\}}
{\bf 1}_{\{i_2=i_4\ne 0\}}
{\bf 1}_{\{j_2=j_4\}}
\zeta_{j_1}^{(i_1)}
\zeta_{j_3}^{(i_3)}
+
{\bf 1}_{\{i_6=i_5\ne 0\}}
{\bf 1}_{\{j_6=j_5\}}
{\bf 1}_{\{i_2=i_3\ne 0\}}
{\bf 1}_{\{j_2=j_3\}}
\zeta_{j_1}^{(i_1)}
\zeta_{j_4}^{(i_4)}+
$$
$$
+
{\bf 1}_{\{i_6=i_5\ne 0\}}
{\bf 1}_{\{j_6=j_5\}}
{\bf 1}_{\{i_1=i_4\ne 0\}}
{\bf 1}_{\{j_1=j_4\}}
\zeta_{j_2}^{(i_2)}
\zeta_{j_3}^{(i_3)}
+
{\bf 1}_{\{i_6=i_5\ne 0\}}
{\bf 1}_{\{j_6=j_5\}}
{\bf 1}_{\{i_1=i_3\ne 0\}}
{\bf 1}_{\{j_1=j_3\}}
\zeta_{j_2}^{(i_2)}
\zeta_{j_4}^{(i_4)}+
$$
$$
+
{\bf 1}_{\{i_6=i_5\ne 0\}}
{\bf 1}_{\{j_6=j_5\}}
{\bf 1}_{\{i_1=i_2\ne 0\}}
{\bf 1}_{\{j_1=j_2\}}
\zeta_{j_3}^{(i_3)}
\zeta_{j_4}^{(i_4)}-
$$
$$
-
{\bf 1}_{\{i_6=i_1\ne 0\}}
{\bf 1}_{\{j_6=j_1\}}
{\bf 1}_{\{i_2=i_5\ne 0\}}
{\bf 1}_{\{j_2=j_5\}}
{\bf 1}_{\{i_3=i_4\ne 0\}}
{\bf 1}_{\{j_3=j_4\}}-
$$
$$
-
{\bf 1}_{\{i_6=i_1\ne 0\}}
{\bf 1}_{\{j_6=j_1\}}
{\bf 1}_{\{i_2=i_4\ne 0\}}
{\bf 1}_{\{j_2=j_4\}}
{\bf 1}_{\{i_3=i_5\ne 0\}}
{\bf 1}_{\{j_3=j_5\}}-
$$
$$
-
{\bf 1}_{\{i_6=i_1\ne 0\}}
{\bf 1}_{\{j_6=j_1\}}
{\bf 1}_{\{i_2=i_3\ne 0\}}
{\bf 1}_{\{j_2=j_3\}}
{\bf 1}_{\{i_4=i_5\ne 0\}}
{\bf 1}_{\{j_4=j_5\}}-
$$
$$
-
{\bf 1}_{\{i_6=i_2\ne 0\}}
{\bf 1}_{\{j_6=j_2\}}
{\bf 1}_{\{i_1=i_5\ne 0\}}
{\bf 1}_{\{j_1=j_5\}}
{\bf 1}_{\{i_3=i_4\ne 0\}}
{\bf 1}_{\{j_3=j_4\}}-
$$
$$
-
{\bf 1}_{\{i_6=i_2\ne 0\}}
{\bf 1}_{\{j_6=j_2\}}
{\bf 1}_{\{i_1=i_4\ne 0\}}
{\bf 1}_{\{j_1=j_4\}}
{\bf 1}_{\{i_3=i_5\ne 0\}}
{\bf 1}_{\{j_3=j_5\}}-
$$
$$
-
{\bf 1}_{\{i_6=i_2\ne 0\}}
{\bf 1}_{\{j_6=j_2\}}
{\bf 1}_{\{i_1=i_3\ne 0\}}
{\bf 1}_{\{j_1=j_3\}}
{\bf 1}_{\{i_4=i_5\ne 0\}}
{\bf 1}_{\{j_4=j_5\}}-
$$
$$
-
{\bf 1}_{\{i_6=i_3\ne 0\}}
{\bf 1}_{\{j_6=j_3\}}
{\bf 1}_{\{i_1=i_5\ne 0\}}
{\bf 1}_{\{j_1=j_5\}}
{\bf 1}_{\{i_2=i_4\ne 0\}}
{\bf 1}_{\{j_2=j_4\}}-
$$
$$
-
{\bf 1}_{\{i_6=i_3\ne 0\}}
{\bf 1}_{\{j_6=j_3\}}
{\bf 1}_{\{i_1=i_4\ne 0\}}
{\bf 1}_{\{j_1=j_4\}}
{\bf 1}_{\{i_2=i_5\ne 0\}}
{\bf 1}_{\{j_2=j_5\}}-
$$
$$
-
{\bf 1}_{\{i_3=i_6\ne 0\}}
{\bf 1}_{\{j_3=j_6\}}
{\bf 1}_{\{i_1=i_2\ne 0\}}
{\bf 1}_{\{j_1=j_2\}}
{\bf 1}_{\{i_4=i_5\ne 0\}}
{\bf 1}_{\{j_4=j_5\}}-
$$
$$
-
{\bf 1}_{\{i_6=i_4\ne 0\}}
{\bf 1}_{\{j_6=j_4\}}
{\bf 1}_{\{i_1=i_5\ne 0\}}
{\bf 1}_{\{j_1=j_5\}}
{\bf 1}_{\{i_2=i_3\ne 0\}}
{\bf 1}_{\{j_2=j_3\}}-
$$
$$
-
{\bf 1}_{\{i_6=i_4\ne 0\}}
{\bf 1}_{\{j_6=j_4\}}
{\bf 1}_{\{i_1=i_3\ne 0\}}
{\bf 1}_{\{j_1=j_3\}}
{\bf 1}_{\{i_2=i_5\ne 0\}}
{\bf 1}_{\{j_2=j_5\}}-
$$
$$
-
{\bf 1}_{\{i_6=i_4\ne 0\}}
{\bf 1}_{\{j_6=j_4\}}
{\bf 1}_{\{i_1=i_2\ne 0\}}
{\bf 1}_{\{j_1=j_2\}}
{\bf 1}_{\{i_3=i_5\ne 0\}}
{\bf 1}_{\{j_3=j_5\}}-
$$
$$
-
{\bf 1}_{\{i_6=i_5\ne 0\}}
{\bf 1}_{\{j_6=j_5\}}
{\bf 1}_{\{i_1=i_4\ne 0\}}
{\bf 1}_{\{j_1=j_4\}}
{\bf 1}_{\{i_2=i_3\ne 0\}}
{\bf 1}_{\{j_2=j_3\}}-
$$
$$
-
{\bf 1}_{\{i_6=i_5\ne 0\}}
{\bf 1}_{\{j_6=j_5\}}
{\bf 1}_{\{i_1=i_2\ne 0\}}
{\bf 1}_{\{j_1=j_2\}}
{\bf 1}_{\{i_3=i_4\ne 0\}}
{\bf 1}_{\{j_3=j_4\}}-
$$
\begin{equation}
\label{a6}
\Biggl.-
{\bf 1}_{\{i_6=i_5\ne 0\}}
{\bf 1}_{\{j_6=j_5\}}
{\bf 1}_{\{i_1=i_3\ne 0\}}
{\bf 1}_{\{j_1=j_3\}}
{\bf 1}_{\{i_2=i_4\ne 0\}}
{\bf 1}_{\{j_2=j_4\}}\Biggr),
\end{equation}

\vspace{6mm}
\noindent
where ${\bf 1}_A$ is the indicator of the set $A$.

Thus, we obtain the following useful possibilities
of the method of generalized multiple Fourier series.

1. There is the explicit formula (see (\ref{ppppa})) for calculation 
of expansion coefficients 
of the iterated Ito stochastic integral (\ref{ito}) with any
fixed multiplicity $k$. 

2. We have new possibilities for exact calculation of the mean-square 
approximation error
for the iterated Ito stochastic integrals (\ref{ito})
(see \cite{17}, \cite{19}, \cite{20}-\cite{20xxz}, \cite{26}).

3. Since the used
multiple Fourier series is a generalized in the sense
that it is built using various complete orthonormal
systems of functions in the space $L_2([t, T])$, then we 
have new possibilities 
for approximation --- we can
use not only trigonometric functions as in \cite{KlPl2}-\cite{Mi3}
but Legendre polynomials.

4. As it turned out (see \cite{3}-\cite{19}, \cite{20}-\cite{33}), 
it is more convenient to work 
with the Legendre polynomials for constructing of approximations 
of the iterated Ito stochastic integrals (\ref{ito}). 
Ap\-pro\-xi\-ma\-ti\-ons based on the Legendre polynomials essentially simpler 
than their analogues based on the trigonometric functions
(see \cite{3}-\cite{19}, \cite{20}-\cite{33}).
Another advantages of the application of Legendre polynomials 
in the framework of the mentioned problem are considered
in \cite{20xx}-\cite{20xxz}, \cite{29}, \cite{30}.

5. The approach based on the Karhunen--Loeve expansion
of the Brownian bridge process \cite{Mi2} (also see \cite{rr})
leads to 
iterated application of the operation of limit
transition (the operation of limit transition 
is implemented only once in Theorem 1 and Theorem 2 (see below))
starting from the 
second multiplicity (in the general case) 
and third multiplicity (for the case
$\psi_1(s), \psi_2(s), \psi_3(s)\equiv 1;$ 
$i_1, i_2, i_3=1,\ldots,m$)
of iterated Ito stochastic integrals.
Multiple series (the operation of limit transition 
is implemented only once) are more convenient 
for approximation than the iterated ones
(iterated application of the operation of limit
transition), 
since partial sums of multiple series converge for any possible case of  
convergence to infinity of their upper limits of summation 
(let us denote them as $p_1,\ldots, p_k$). 
For example, 
when $p_1=\ldots=p_k=p\to\infty$. 
For iterated series, the condition $p_1=\ldots=p_k=p\to\infty$ obviously 
does not guarantee the convergence of this series.
However, in 
\cite{KlPl2}
(Sect.~5.8, pp.~202--204), \cite{KPS} (pp.~82-84),
\cite{KPW} (pp.~438-439),  
\cite{Zapad-9} (pp.~263-264) the authors use 
(without rigorous proof)
the condition $p_1=p_2=p_3=p\to\infty$
within the frames of the mentioned approach
based on the Karhunen--Loeve expansion of the Brownian bridge
process \cite{Mi2} together with the Wong--Zakai approximation
\cite{W-Z-1}-\cite{Watanabe} (see discussion
in Sect.~7 of this paper for details).

For further consideration, let us 
consider the generalization of formulas (\ref{a1})--(\ref{a6})                 
for the case of an arbitrary multiplicity $k$ $(k\in\mathbb{N})$ of 
the iterated Ito stochastic integral $J[\psi^{(k)}]_{T,t}$ defined by (\ref{ito}).
In order to do this, let us
introduce some notations. 
Consider the unordered
set $\{1, 2, \ldots, k\}$ 
and separate it into two parts:
the first part consists of $r$ unordered 
pairs (sequence order of these pairs is also unimportant) and the 
second one consists of the 
remaining $k-2r$ numbers.
So, we have

\begin{equation}
\label{leto5007}
(\{
\underbrace{\{g_1, g_2\}, \ldots, 
\{g_{2r-1}, g_{2r}\}}_{\small{\hbox{part 1}}}
\},
\{\underbrace{q_1, \ldots, q_{k-2r}}_{\small{\hbox{part 2}}}
\}),
\end{equation}

\vspace{4mm}
\noindent
where 

\vspace{-2mm}
$$
\{g_1, g_2, \ldots, 
g_{2r-1}, g_{2r}, q_1, \ldots, q_{k-2r}\}=\{1, 2, \ldots, k\},
$$

\vspace{4mm}
\noindent
braces   
mean an unordered 
set, and pa\-ren\-the\-ses mean an ordered set.

We will say that (\ref{leto5007}) is a partition 
and consider the sum with respect to all possible
partitions

\begin{equation}
\label{leto5008}
\sum_{\stackrel{(\{\{g_1, g_2\}, \ldots, 
\{g_{2r-1}, g_{2r}\}\}, \{q_1, \ldots, q_{k-2r}\})}
{{}_{\{g_1, g_2, \ldots, 
g_{2r-1}, g_{2r}, q_1, \ldots, q_{k-2r}\}=\{1, 2, \ldots, k\}}}}
a_{g_1 g_2, \ldots, 
g_{2r-1} g_{2r}, q_1 \ldots q_{k-2r}}.
\end{equation}

\vspace{4mm}

Below there are several examples of sums in the form (\ref{leto5008})

\vspace{2mm}
$$
\sum_{\stackrel{(\{g_1, g_2\})}{{}_{\{g_1, g_2\}=\{1, 2\}}}}
a_{g_1 g_2}=a_{12},
$$

\vspace{3mm}
$$
\sum_{\stackrel{(\{\{g_1, g_2\}, \{g_3, g_4\}\})}
{{}_{\{g_1, g_2, g_3, g_4\}=\{1, 2, 3, 4\}}}}
a_{g_1 g_2, g_3 g_4}=a_{12,34} + a_{13,24} + a_{23,14},
$$

\vspace{3mm}
$$
\sum_{\stackrel{(\{g_1, g_2\}, \{q_1, q_{2}\})}
{{}_{\{g_1, g_2, q_1, q_{2}\}=\{1, 2, 3, 4\}}}}
a_{g_1 g_2, q_1 q_{2}}=
$$

$$
=a_{12,34}+a_{13,24}+a_{14,23}
+a_{23,14}+a_{24,13}+a_{34,12},
$$

\vspace{3mm}
$$
\sum_{\stackrel{(\{g_1, g_2\}, \{q_1, q_{2}, q_3\})}
{{}_{\{g_1, g_2, q_1, q_{2}, q_3\}=\{1, 2, 3, 4, 5\}}}}
a_{g_1 g_2, q_1 q_{2}q_3}
=
$$

$$
=a_{12,345}+a_{13,245}+a_{14,235}
+a_{15,234}+a_{23,145}+a_{24,135}+
$$
$$
+a_{25,134}+a_{34,125}+a_{35,124}+a_{45,123},
$$

\vspace{4mm}
$$
\sum_{\stackrel{(\{\{g_1, g_2\}, \{g_3, g_{4}\}\}, \{q_1\})}
{{}_{\{g_1, g_2, g_3, g_{4}, q_1\}=\{1, 2, 3, 4, 5\}}}}
a_{g_1 g_2, g_3 g_{4},q_1}
=
$$

$$
=
a_{12,34,5}+a_{13,24,5}+a_{14,23,5}+
a_{12,35,4}+a_{13,25,4}+a_{15,23,4}+
$$
$$
+a_{12,54,3}+a_{15,24,3}+a_{14,25,3}+a_{15,34,2}+a_{13,54,2}+a_{14,53,2}+
$$
$$
+
a_{52,34,1}+a_{53,24,1}+a_{54,23,1}.
$$

\vspace{5mm}

Now we can write (\ref{tyyy}) as

\vspace{1mm}

$$
J[\psi^{(k)}]_{T,t}=
\hbox{\vtop{\offinterlineskip\halign{
\hfil#\hfil\cr
{\rm l.i.m.}\cr
$\stackrel{}{{}_{p_1,\ldots,p_k\to \infty}}$\cr
}} }
\sum\limits_{j_1=0}^{p_1}\ldots
\sum\limits_{j_k=0}^{p_k}
C_{j_k\ldots j_1}\Biggl(
\prod_{l=1}^k\zeta_{j_l}^{(i_l)}+\sum\limits_{r=1}^{[k/2]}
(-1)^r \times
\Biggr.
$$

\vspace{3mm}
\begin{equation}
\label{leto6000hh}
\times
\sum_{\stackrel{(\{\{g_1, g_2\}, \ldots, 
\{g_{2r-1}, g_{2r}\}\}, \{q_1, \ldots, q_{k-2r}\})}
{{}_{\{g_1, g_2, \ldots, 
g_{2r-1}, g_{2r}, q_1, \ldots, q_{k-2r}\}=\{1, 2, \ldots, k\}}}}
\prod\limits_{s=1}^r
{\bf 1}_{\{i_{g_{{}_{2s-1}}}=~i_{g_{{}_{2s}}}\ne 0\}}
\Biggl.{\bf 1}_{\{j_{g_{{}_{2s-1}}}=~j_{g_{{}_{2s}}}\}}
\prod_{l=1}^{k-2r}\zeta_{j_{q_l}}^{(i_{q_l})}\Biggr),
\end{equation}

\vspace{5mm}
\noindent
where $[x]$ is an integer part of a real number $x;$
another notations are the same as in Theorem {\bf 1}.

\vspace{2mm}

In particular, from (\ref{leto6000hh}) for $k=5$ we obtain

\vspace{3mm}

$$
J[\psi^{(5)}]_{T,t}=
\hbox{\vtop{\offinterlineskip\halign{
\hfil#\hfil\cr
{\rm l.i.m.}\cr
$\stackrel{}{{}_{p_1,\ldots,p_5\to \infty}}$\cr
}} }\sum_{j_1=0}^{p_1}\ldots\sum_{j_5=0}^{p_5}
C_{j_5\ldots j_1}\Biggl(
\prod_{l=1}^5\zeta_{j_l}^{(i_l)}-\Biggr.
$$

\vspace{2mm}
$$
-
\sum\limits_{\stackrel{(\{g_1, g_2\}, \{q_1, q_{2}, q_3\})}
{{}_{\{g_1, g_2, q_{1}, q_{2}, q_3\}=\{1, 2, 3, 4, 5\}}}}
{\bf 1}_{\{i_{g_{{}_{1}}}=~i_{g_{{}_{2}}}\ne 0\}}
{\bf 1}_{\{j_{g_{{}_{1}}}=~j_{g_{{}_{2}}}\}}
\prod_{l=1}^{3}\zeta_{j_{q_l}}^{(i_{q_l})}+
$$

\vspace{2mm}
$$
+
\sum_{\stackrel{(\{\{g_1, g_2\}, 
\{g_{3}, g_{4}\}\}, \{q_1\})}
{{}_{\{g_1, g_2, g_{3}, g_{4}, q_1\}=\{1, 2, 3, 4, 5\}}}}
{\bf 1}_{\{i_{g_{{}_{1}}}=~i_{g_{{}_{2}}}\ne 0\}}
{\bf 1}_{\{j_{g_{{}_{1}}}=~j_{g_{{}_{2}}}\}}
\Biggl.{\bf 1}_{\{i_{g_{{}_{3}}}=~i_{g_{{}_{4}}}\ne 0\}}
{\bf 1}_{\{j_{g_{{}_{3}}}=~j_{g_{{}_{4}}}\}}
\zeta_{j_{q_1}}^{(i_{q_1})}\Biggr).
$$

\vspace{5mm}
\noindent
The last equality obviously agrees with
(\ref{a5}).

Let us consider a generalization of Theorem 1 for the case
of an arbitrary complete orthonormal systems  
of functions in the space $L_2([t,T])$ 
and $\psi_1(\tau),\ldots,\psi_k(\tau)\in L_2([t, T]).$

\vspace{2mm}

{\bf Theorem~2}\ \cite{20xx} (Sect.~1.11), \cite{26a} (Sect.~15).
{\it Suppose that
$\psi_1(\tau),\ldots,\psi_k(\tau)\in L_2([t, T])$ and
$\{\phi_j(x)\}_{j=0}^{\infty}$ is an arbitrary complete orthonormal system  
of functions in the space $L_2([t,T]).$
Then the following expansion

\vspace{1mm}
$$
J[\psi^{(k)}]_{T,t}=
\hbox{\vtop{\offinterlineskip\halign{
\hfil#\hfil\cr
{\rm l.i.m.}\cr
$\stackrel{}{{}_{p_1,\ldots,p_k\to \infty}}$\cr
}} }
\sum\limits_{j_1=0}^{p_1}\ldots
\sum\limits_{j_k=0}^{p_k}
C_{j_k\ldots j_1}\Biggl(
\prod_{l=1}^k\zeta_{j_l}^{(i_l)}+\sum\limits_{r=1}^{[k/2]}
(-1)^r \times
\Biggr.
$$

\vspace{2mm}
\begin{equation}
\label{leto6000}
\times
\sum_{\stackrel{(\{\{g_1, g_2\}, \ldots, 
\{g_{2r-1}, g_{2r}\}\}, \{q_1, \ldots, q_{k-2r}\})}
{{}_{\{g_1, g_2, \ldots, 
g_{2r-1}, g_{2r}, q_1, \ldots, q_{k-2r}\}=\{1, 2, \ldots, k\}}}}
\prod\limits_{s=1}^r
{\bf 1}_{\{i_{g_{{}_{2s-1}}}=~i_{g_{{}_{2s}}}\ne 0\}}
\Biggl.{\bf 1}_{\{j_{g_{{}_{2s-1}}}=~j_{g_{{}_{2s}}}\}}
\prod_{l=1}^{k-2r}\zeta_{j_{q_l}}^{(i_{q_l})}\Biggr)
\end{equation}

\vspace{6mm}
\noindent
con\-verg\-ing in the mean-square sense is valid,
where $[x]$ is an integer part of a real number $x;$
another notations are the same as in Theorem~{\rm 1}.}

\vspace{2mm}

It should be noted that an analogue of Theorem 2 was considered 
in \cite{Rybakov1000}. 
Note that we use another notations 
\cite{20xx} (Sect.~1.11), \cite{26a} (Sect.~15)
in comparison with \cite{Rybakov1000}.
Moreover, the proof of an analogue of Theorem 2
from \cite{Rybakov1000} is somewhat different from the proof given in 
\cite{20xx} (Sect.~1.11), \cite{26a} (Sect.~15).

\vspace{5mm}

\section{Expansions of Iterated Stratonovich Stochastic Integrals 
of Multiplicities 2 to 4. Some Old Results}

\vspace{5mm}

As it turned out, Theorems 1, 2 can be adapted for the iterated
Stratonovich stochastic integrals (\ref{str}) at least
for multiplicities 1 to 6 (the case $k=1$ obviously corresponds to
(\ref{a1})). 
Expansions of the mentioned iterated Stratonovich 
stochastic integrals turned out
simpler than the appropriate expansions
for the iterated Ito stochastic integrals (\ref{ito}) based on Theorems 1, 2.
Let us formulate some theorems on expansions of the iterated
Stratonovich stochastic integrals (\ref{str}) of
multiplicities 2 to 4.

\vspace{2mm}

{\bf Theorem 3} \cite{14}-\cite{16}, \cite{19}, \cite{20}-\cite{20xxz}, 
\cite{32}.
{\it Suppose that 
$\{\phi_j(x)\}_{j=0}^{\infty}$ is a complete orthonormal system of 
Legendre polynomials or tri\-go\-no\-met\-ric functions in the space $L_2([t, T]).$
At the same time $\psi_2(\tau)$ is a continuously dif\-ferentiable 
function on $[t, T]$ and $\psi_1(\tau)$ is twice 
continuously differentiable function on $[t, T]$. 
Then, the iterated Stratonovich stochastic integral 
of second multiplicity

$$
J^{*}[\psi^{(2)}]_{T,t}={\int\limits_t^{*}}^T\psi_2(t_2)
{\int\limits_t^{*}}^{t_2}\psi_1(t_1)d{\bf w}_{t_1}^{(i_1)}
d{\bf w}_{t_2}^{(i_2)}\ \ \ (i_1, i_2=0, 1,\ldots,m)
$$

\vspace{3mm}
\noindent
is expanded into the 
multiple series

$$
J^{*}[\psi^{(2)}]_{T,t}=
\hbox{\vtop{\offinterlineskip\halign{
\hfil#\hfil\cr
{\rm l.i.m.}\cr
$\stackrel{}{{}_{p_1, p_2\to \infty}}$\cr
}} }
\sum\limits_{j_1=0}^{p_1}
\sum\limits_{j_2=0}^{p_2}
C_{j_2j_1}\zeta_{j_1}^{(i_1)}\zeta_{j_2}^{(i_2)}
$$

\vspace{3mm}
\noindent
converging in the mean-square sense, where 
$$
\zeta_{j}^{(i)}=
\int\limits_t^T \phi_{j}(\tau) d{\bf w}_{\tau}^{(i)}
$$ 

\vspace{2mm}
\noindent
are independent standard Gaussian random variables for various 
$i$\ or $j$\ {\rm (}if $i\ne 0${\rm ),}

$$
C_{j_2 j_1}=\int\limits_t^T\psi_2(t_2)\phi_{j_2}(t_2)
\int\limits_t^{t_2}\psi_1(t_1)\phi_{j_1}(t_1)dt_1dt_2
$$

\vspace{3mm}
\noindent
is the Fourier coefficient.}

\vspace{2mm}

Note that in \cite{20}-\cite{20xxz}, \cite{30a}, \cite{33} 
Theorem 3 is proved under
weaker conditions.

\vspace{2mm}

{\bf Theorem 4}\ \cite{20}-\cite{20xxz}, \cite{30a}, \cite{33}.
{\it Suppose that the following conditions are fulfilled{\rm :}

{\rm 1}. The functions $\psi_1(\tau)$ and $\psi_2(\tau)$ are continuously 
differentiable at the interval $[t, T]$.

{\rm 2}. $\{\phi_j(x)\}_{j=0}^{\infty}$ is a complete orthonormal system 
of Legendre polynomials or trigo\-no\-met\-ric functions
in the space $L_2([t, T])$.

Then, the iterated Stratonovich stochastic integral 
of second multiplicity

$$
J^{*}[\psi^{(2)}]_{T,t}={\int\limits_t^{*}}^T\psi_2(t_2)
{\int\limits_t^{*}}^{t_2}\psi_1(t_1)d{\bf w}_{t_1}^{(i_1)}
d{\bf w}_{t_2}^{(i_2)}\ \ \ (i_1, i_2=0, 1,\ldots,m)
$$

\vspace{3mm}
\noindent
is expanded into the 
multiple series

$$
J^{*}[\psi^{(2)}]_{T,t}=
\hbox{\vtop{\offinterlineskip\halign{
\hfil#\hfil\cr
{\rm l.i.m.}\cr
$\stackrel{}{{}_{p_1, p_2\to \infty}}$\cr
}} }
\sum\limits_{j_1=0}^{p_1}
\sum\limits_{j_2=0}^{p_2}
C_{j_2j_1}\zeta_{j_1}^{(i_1)}\zeta_{j_2}^{(i_2)}
$$

\vspace{3mm}
\noindent
converging in the mean-square sense, where 
$$
\zeta_{j}^{(i)}=
\int\limits_t^T \phi_{j}(\tau) d{\bf w}_{\tau}^{(i)}
$$ 

\vspace{2mm}
\noindent
are independent standard Gaussian random variables for various 
$i$\ or $j$\ {\rm (}if $i\ne 0${\rm ),}

$$
C_{j_2 j_1}=\int\limits_t^T\psi_2(t_2)\phi_{j_2}(t_2)
\int\limits_t^{t_2}\psi_1(t_1)\phi_{j_1}(t_1)dt_1dt_2
$$

\vspace{3mm}
\noindent
is the Fourier coefficient.}

\vspace{3mm}

{\bf Theorem 5}\ \cite{15}, \cite{16}, \cite{19}, \cite{20}-\cite{20xxz},
\cite{32}. 
{\it Suppose that
$\{\phi_j(x)\}_{j=0}^{\infty}$ is a complete orthonormal
system of Legendre polynomials or trigonometric functions
in the space $L_2([t, T])$. At the same time
$\psi_2(\tau)$ is a continuously
differentiable function at the interval $[t, T]$ and
$\psi_1(\tau), \psi_3(\tau)$ are twice continuously
differentiable functions at the interval $[t, T]$.
Then, for the iterated 
Stratonovich sto\-chas\-tic integral of third multiplicity

$$
J^{*}[\psi^{(3)}]_{T,t}={\int\limits_t^{*}}^T\psi_3(t_3)
{\int\limits_t^{*}}^{t_3}\psi_2(t_2)
{\int\limits_t^{*}}^{t_2}\psi_1(t_1)
d{\bf f}_{t_1}^{(i_1)}
d{\bf f}_{t_2}^{(i_2)}d{\bf f}_{t_3}^{(i_3)}\ \ \ (i_1, i_2, i_3=1,\ldots,m)
$$

\vspace{3mm}
\noindent
the following 
expansion 

\vspace{-1mm}
\begin{equation}
\label{newbegin96}
J^{*}[\psi^{(3)}]_{T,t}=
\hbox{\vtop{\offinterlineskip\halign{
\hfil#\hfil\cr
{\rm l.i.m.}\cr
$\stackrel{}{{}_{p\to \infty}}$\cr
}} }
\sum\limits_{j_1, j_2, j_3=0}^{p}
C_{j_3 j_2 j_1}\zeta_{j_1}^{(i_1)}\zeta_{j_2}^{(i_2)}\zeta_{j_3}^{(i_3)}
\end{equation}

\vspace{4mm}
\noindent
converging in the mean-square sense is valid, where

$$
C_{j_3 j_2 j_1}=\int\limits_t^T\psi_3(t_3)\phi_{j_3}(t_3)
\int\limits_t^{t_3}\psi_2(t_2)\phi_{j_2}(t_2)
\int\limits_t^{t_2}\psi_1(t_1)\phi_{j_1}(t_1)dt_1dt_2dt_3,
$$

\vspace{2mm}
\noindent
another notations are the same as in Theorems {\rm 1, 2}.}

\vspace{3mm}

{\bf Theorem 6}\ \cite{14}-\cite{16}, \cite{19}, \cite{20}-\cite{20xxz}, 
\cite{32}. 
{\it Suppose that
$\{\phi_j(x)\}_{j=0}^{\infty}$ is a complete orthonormal
system of Legendre polynomials or trigonometric functions
in the space $L_2([t, T])$.
Then, for the iterated 
Stratonovich stochastic integral of fourth multiplicity

$$
I_{T,t}^{*(i_1 i_2 i_3 i_4)}=
{\int\limits_t^{*}}^T
{\int\limits_t^{*}}^{t_4}
{\int\limits_t^{*}}^{t_3}
{\int\limits_t^{*}}^{t_2}
d{\bf w}_{t_1}^{(i_1)}
d{\bf w}_{t_2}^{(i_2)}d{\bf w}_{t_3}^{(i_3)}d{\bf w}_{t_4}^{(i_4)}\ \ \ 
(i_1, i_2, i_3, i_4=0, 1,\ldots,m)
$$

\vspace{3mm}
\noindent
the following 
expansion 

\vspace{-1mm}
$$
I_{T,t}^{*(i_1 i_2 i_3 i_4)}=
\hbox{\vtop{\offinterlineskip\halign{
\hfil#\hfil\cr
{\rm l.i.m.}\cr
$\stackrel{}{{}_{p\to \infty}}$\cr
}} }
\sum\limits_{j_1, j_2, j_3, j_4=0}^{p}
C_{j_4 j_3 j_2 j_1}\zeta_{j_1}^{(i_1)}\zeta_{j_2}^{(i_2)}\zeta_{j_3}^{(i_3)}
\zeta_{j_4}^{(i_4)}
$$

\vspace{4mm}
\noindent
converging in the mean-square sense is valid, where

$$
\zeta_{j}^{(i)}=
\int\limits_t^T \phi_{j}(s) d{\bf w}_s^{(i)}
$$ 

\vspace{3mm}
\noindent
are independent standard Gaussian random variables for various 
$i$ or $j$\ {\rm (}if $i\ne 0${\rm ),}

$$
C_{j_4 j_3 j_2 j_1}=\int\limits_t^T\phi_{j_4}(t_4)\int\limits_t^{t_4}
\phi_{j_3}(t_3)
\int\limits_t^{t_3}\phi_{j_2}(t_2)\int\limits_t^{t_2}\phi_{j_1}(t_1)
dt_1dt_2dt_3dt_4,
$$

\vspace{4mm}
\noindent
${\bf w}_{\tau}^{(i)}={\bf f}_{\tau}^{(i)}$ for
$i=1,\ldots,m$ and 
${\bf w}_{\tau}^{(0)}=\tau.$}

\vspace{2mm}

Note that 
in \cite{14}-\cite{16}, \cite{19}, \cite{20}-\cite{20xxz}, \cite{32}
the expansions (\ref{a1})--(\ref{a4}) 
have been applied 
for the proof
of Theorems 3--6. In this article, 
we will prove Theorems 4--6 by an another approach.
This approach will be called as the combined approach.
More precisely, we will use the scheme of the proof
of Theorem 1 from this paper (see 
\cite{7}--\cite{16}, \cite{19}, \cite{20}-\cite{20xxz}, 
\cite{26a} for details) 
for the iterated Stratonovich
stochastic integrals (\ref{str}) of multiplicities 2 to 4.
As a result, we will obtain 
two different parts of the expansion of 
iterated Stra\-to\-no\-vich stochastic integrals.
The mean-square convergence of the first part will be proved
on the base
of generalized multiple Fourier series converging 
in $L_2([t, T]^k)$ $(k=2, 3, 4)$. At the same time, the 
mean-square convergence
of the second part will be proved on the base of 
generalized iterated Fourier
series converging pointwise. At that, we do not use the iterated Ito
stochastic integrals (\ref{ito}) 
as a tool of the proof and directly consider
the iterated Stratonovich stochastic integrals (\ref{str}).

\vspace{5mm}

\section{Auxiliary Lemmas}

\vspace{5mm}

In this section, we collected several lemmas, which will be used
for the proof of Theorems 4--6.

Consider the partition $\{\tau_j\}_{j=0}^N$ of the interval $[t,T]$ such that

\vspace{-1mm}
\begin{equation}
\label{1111aaa}
t=\tau_0<\ldots <\tau_N=T,\ \ \
\Delta_N=
\hbox{\vtop{\offinterlineskip\halign{
\hfil#\hfil\cr
{\rm max}\cr
$\stackrel{}{{}_{0\le j\le N-1}}$\cr
}} }\Delta\tau_j\to 0\ \ \hbox{if}\ \ N\to\ \infty,\ \ \ 
\Delta\tau_j=\tau_{j+1}-\tau_j.
\end{equation}

\vspace{3mm}

{\bf Lemma 1}\ \cite{7}--\cite{16}, \cite{19}, \cite{20}-\cite{20xxz}, 
\cite{26a}. 
{\it Suppose that
every $\psi_l(\tau)$ $(l=1,\ldots, k)$ is a continuous 
nonrandom function at the interval
$[t, T]$. Then

\vspace{-1mm}
\begin{equation}
\label{30.30}
J[\psi^{(k)}]_{T,t}=
\hbox{\vtop{\offinterlineskip\halign{
\hfil#\hfil\cr
{\rm l.i.m.}\cr
$\stackrel{}{{}_{N\to \infty}}$\cr
}} }
\sum_{j_k=0}^{N-1}
\ldots \sum_{j_1=0}^{j_{2}-1}
\prod_{l=1}^k \psi_l(\tau_{j_l})\Delta{\bf w}_{\tau_
{j_l}}^{(i_l)}\ \ \ \hbox{w. p. {\rm 1}},
\end{equation}

\vspace{3mm}
\noindent
where $J[\psi^{(k)}]_{T,t}$ has the form {\rm (\ref{ito}),}\ 
$\Delta{\bf w}_{\tau_{j}}^{(i)}=
{\bf w}_{\tau_{j+1}}^{(i)}-{\bf w}_{\tau_{j}}^{(i)}$
$(i=0,\ 1,\ldots,m)$,\
$\left\{\tau_{j}\right\}_{j=0}^{N}$ is a partition 
of the interval $[t,T]$ satisfying the condition {\rm (\ref{1111aaa});}
hereinafter w.~p.~{\rm 1}  means  with probability {\rm 1}.
}

\vspace{2mm}

{\bf Remark 1.}\ \it It is easy to see that if
$\Delta{\bf w}_{\tau_{j_l}}^{(i_l)}$ in {\rm (\ref{30.30})}
for some $l\in\{1,\ldots,k\}$ is replaced with $\left(
\Delta{\bf w}_{\tau_{j_l}}^{(i_l)}\right)^p$ $(p=2,$
$i_l\ne 0),$ then
the differential $d{\bf w}_{t_{l}}^{(i_l)}$
in the integral $J[\psi^{(k)}]_{T,t}$
will be replaced with $dt_l$.\ 
If $p=3, 4,\ldots,$ then the
right-hand side
of the formula {\rm (\ref{30.30})}
will become zero w.~p.~{\rm 1}.
If we replace $\Delta{\bf w}_{\tau_{j_l}}^{(i_l)}$ in {\rm (\ref{30.30})}
for some $l\in\{1,\ldots,k\}$
with $\left(
\Delta \tau_{j_l}\right)^p$ $(p=2, 3,\ldots ),$
then the right-hand side of the formula
{\rm (\ref{30.30})} 
will also be equal to zero w.~p.~{\rm 1}.
\rm

\vspace{2mm}

Let us define the following
multiple stochastic integral

\begin{equation}
\label{30.34}
\hbox{\vtop{\offinterlineskip\halign{
\hfil#\hfil\cr
{\rm l.i.m.}\cr
$\stackrel{}{{}_{N\to \infty}}$\cr
}} }\sum_{j_1,\ldots,j_k=0}^{N-1}
\Phi\left(\tau_{j_1},\ldots,\tau_{j_k}\right)
\prod\limits_{l=1}^k\Delta{\bf w}_{\tau_{j_l}}^{(i_l)}
\stackrel{\rm def}{=}J[\Phi]_{T,t}^{(k)},
\end{equation}

\vspace{3mm}
\noindent
where $\Phi(t_1,\ldots,t_k):\ [t, T]^k\to\mathbb{R}$ is a nonrandom function 
(the properties of this function
will be specified further),
$\Delta{\bf w}_{\tau_{j}}^{(i)}=
{\bf w}_{\tau_{j+1}}^{(i)}-{\bf w}_{\tau_{j}}^{(i)}$
$(i=0,\ 1,\ldots,m)$,
$\left\{\tau_{j}\right\}_{j=0}^{N}$ is a partition 
of the interval $[t,T]$ satisfying the condition {\rm (\ref{1111aaa})}.

Denote
\begin{equation}
\label{dom1}
D_k=\{(t_1,\ldots,t_k):\ t\le t_1<\ldots <t_k\le T\}.
\end{equation}

\vspace{3mm}

We will use the same symbol $D_k$ to denote the open and closed 
domains corresponding to the domain $D_k$ defined by (\ref{dom1}).
However, we always specify what domain we consider (open or closed).

Also we will write $\Phi(t_1,\ldots,t_k)\in C(D_k)$
if 
$\Phi(t_1,\ldots,t_k)$ is a continuous nonrandom function of $k$ variables
in the closed domain $D_k$.

Let us consider the iterated Ito stochastic integral

\vspace{-1mm}
$$
I[\Phi]_{T,t}^{(k)}\stackrel{\rm def}{=}
\int\limits_t^T\ldots \int\limits_t^{t_2}
\Phi(t_1,\ldots,t_k)d{\bf w}_{t_1}^{(i_1)}\ldots
d{\bf w}_{t_k}^{(i_k)},
$$

\vspace{2mm}
\noindent
where $\Phi(t_1,\ldots,t_k)\in C(D_k).$

\vspace{2mm}

{\bf Lemma 2}\ \cite{7}--\cite{16}, \cite{19}, \cite{20}-\cite{20xxz}, 
\cite{26a}. 
{\it Suppose that $\Phi(t_1,\ldots,t_k)\in C(D_k)$ 
or $\Phi(t_1,\ldots,t_k)$ 
is a continuous nonrandom function in the open domain $D_k$ and bounded at its boundary.
Then

\begin{equation}
\label{30.52}
I[\Phi]_{T,t}^{(k)}=\hbox{\vtop{\offinterlineskip\halign{
\hfil#\hfil\cr
{\rm l.i.m.}\cr
$\stackrel{}{{}_{N\to \infty}}$\cr
}} }
\sum_{j_k=0}^{N-1}
\ldots \sum_{j_1=0}^{j_{2}-1}
\Phi(\tau_{j_1},\ldots,\tau_{j_k})
\prod\limits_{l=1}^k\Delta {\bf w}_{\tau_{j_l}}^{(i_l)}\ \ \ 
\hbox{w. p. {\rm 1}},
\end{equation}

\vspace{4mm}
\noindent
where $\Delta{\bf w}_{\tau_{j}}^{(i)}=
{\bf w}_{\tau_{j+1}}^{(i)}-{\bf w}_{\tau_{j}}^{(i)}$
$(i=0, 1,\ldots,m)$,\
$\left\{\tau_{j}\right\}_{j=0}^{N}$ is a partition 
of the interval $[t,T]$ satisfying the condition {\rm (\ref{1111aaa})}.
}

\vspace{3mm}

{\bf Lemma 3}\ \cite{7}--\cite{16}, \cite{19}, \cite{20}-\cite{20xxz}, 
\cite{26a}. 
{\it Suppose that every $\varphi_i(\tau)$
$(i=1,\ldots,k)$ is a continuous nonrandom function at the interval $[t, T]$.
Then

\vspace{-1mm}
\begin{equation}
\label{30.39}
\prod_{l=1}^k 
J[\varphi_l]_{T,t}=J[\Phi]_{T,t}^{(k)}\ \ \ \hbox{w.~p.~{\rm 1}},
\end{equation}

\vspace{4mm}
\noindent
where 
$$
J[\varphi_l]_{T,t}
=\int\limits_t^T \varphi_l(s) d{\bf w}_{s}^{(i_l)},\ \ \
\Phi(t_1,\ldots,t_k)=\prod\limits_{l=1}^k\varphi_l(t_l)
$$

\vspace{2mm}
\noindent
and the integral $J[\Phi]_{T,t}^{(k)}$ 
is defined
by the equality
{\rm (\ref{30.34})}.
}

\vspace{2mm}

Let us introduce the following notations

\vspace{1mm}
$$
J[\psi^{(k)}]_{T,t}^{s_l,\ldots,s_1}\ 
\stackrel{\rm def}{=}\  \prod_{p=1}^l {\bf 1}_{\{i_{s_p}=
i_{s_{p}+1}\ne 0\}}\ \times
$$

$$
\times\
\int\limits_t^T\psi_k(t_k)\ 
\ldots \int\limits_t^{t_{s_l+3}}\psi_{s_l+2}(t_{s_l+2})
\int\limits_t^{t_{s_l+2}}\psi_{s_l}(t_{s_l+1})
\psi_{s_l+1}(t_{s_l+1})\ \times
$$

$$
\times
\int\limits_t^{t_{s_l+1}}\psi_{s_l-1}(t_{s_l-1})\
\ldots
\int\limits_t^{t_{s_1+3}}\psi_{s_1+2}(t_{s_1+2})
\int\limits_t^{t_{s_1+2}}\psi_{s_1}(t_{s_1+1})
\psi_{s_1+1}(t_{s_1+1})\ \times
$$

$$
\times
\int\limits_t^{t_{s_1+1}}\psi_{s_1-1}(t_{s_1-1})\
\ldots \int\limits_t^{t_2}\psi_1(t_1)
d{\bf w}_{t_1}^{(i_1)}\ldots\ d{\bf w}_{t_{s_1-1}}^{(i_{s_1-1})}
dt_{s_1+1}d{\bf w}_{t_{s_1+2}}^{(i_{s_1+2})}\ \ldots
$$

\begin{equation}
\label{30.1}
\ldots\
d{\bf w}_{t_{s_l-1}}^{(i_{s_l-1})}
dt_{s_l+1}d{\bf w}_{t_{s_l+2}}^{(i_{s_l+2})}\ldots\ d{\bf w}_{t_k}^{(i_k)},
\end{equation}

\vspace{2mm}
\noindent
where 
\begin{equation}
\label{30.5550001}
{\rm A}_{k,l}=\biggl\{(s_l,\ldots,s_1):\
s_l>s_{l-1}+1,\ldots,s_2>s_1+1,\ s_l,\ldots,s_1=1,\ldots,k-1\biggr\},
\end{equation}

\vspace{5mm}
\noindent
$(s_l,\ldots,s_1)\in{\rm A}_{k,l},$\ \
$l=1,\ldots,\left[k/2\right],$\ \
$i_s=0, 1,\ldots,m,$\ \
$s=1,\ldots,k,$\ \ $[x]$ is an
integer
part of a real number $x,$\ \
${\bf 1}_A$ is the indicator of the set $A$.

\vspace{2mm}

{\bf Lemma 4}\ \cite{3} (1997), \cite{4}, \cite{7}-\cite{16}, \cite{19}, 
\cite{20}-\cite{20xxz}, \cite{23}. {\it Suppose that
every $\psi_l(\tau) (l=1,\ldots,k)$ is a continuous
nonrandom
function at the interval $[t, T]$.
Then, the following relation between iterated
Stra\-to\-no\-vich and Ito stochastic integrals is correct 

\begin{equation}
\label{30.4}
J^{*}[\psi^{(k)}]_{T,t}=J[\psi^{(k)}]_{T,t}+
\sum_{r=1}^{\left[k/2\right]}\frac{1}{2^r}
\sum_{(s_r,\ldots,s_1)\in {\rm A}_{k,r}}
J[\psi^{(k)}]_{T,t}^{s_r,\ldots,s_1}\ \ \ \hbox{w.~p.~{\rm 1}},
\end{equation}

\vspace{3mm}
\noindent
where $\sum\limits_{\emptyset}$ is supposed to be equal to zero.}

\vspace{3mm}

Let us define the function $K^{*}(t_1,\ldots,t_k)$ on the hypercube 
$[t,T]^k$ $(k\ge 2)$ by the following relation

\vspace{2mm}
$$
K^{*}(t_1,\ldots,t_k)=\prod\limits_{l=1}^k\psi_l(t_l)
\prod_{l=1}^{k-1}\Biggl({\bf 1}_{\{t_l<t_{l+1}\}}+
\frac{1}{2}{\bf 1}_{\{t_l=t_{l+1}\}}\Biggr)=
$$

\vspace{2mm}
\begin{equation}
\label{1999.1}
=\prod_{l=1}^k \psi_l(t_l)\left(\prod_{l=1}^{k-1}
{\bf 1}_{\{t_l<t_{l+1}\}}+
\sum_{r=1}^{k-1}\frac{1}{2^r}
\sum_{\stackrel{s_r,\ldots,s_1=1}{{}_{s_r>\ldots>s_1}}}^{k-1}
\prod_{l=1}^r {\bf 1}_{\{t_{s_l}=t_{s_l+1}\}}
\prod_{\stackrel{l=1}{{}_{l\ne s_1,\ldots, s_r}}}^{k-1}
{\bf 1}_{\{t_{l}<t_{l+1}\}}\right),
\end{equation}

\vspace{6mm}
\noindent
where ${\bf 1}_A$ is the indicator of the set $A$.

\vspace{2mm}

{\bf Lemma 5}\ \cite{3}, \cite{4}, \cite{11}-\cite{16}, \cite{19}, 
\cite{20}-\cite{20xxz}, \cite{23}.
{\it Under the conditions of Lemma {\rm 4}
the following relation is correct

\vspace{-3mm}
\begin{equation}
\label{30.36}
J[{K^{*}}]_{T,t}^{(k)}=
J^{*}[\psi^{(k)}]_{T,t}\ \ \ \hbox{w.~p.~{\rm 1}},
\end{equation}

\vspace{4mm}
\noindent
where $J[{K^{*}}]_{T,t}^{(k)}$ is defined by the equality 
{\rm (\ref{30.34})}.}

\vspace{3mm}

{\bf Proof.} Substituting
(\ref{1999.1}) into (\ref{30.34}) and using Lemmas 1, 2, 4 with Remark 1,
it is easy to notice that w.~p.~1

\vspace{-2mm}
\begin{equation}
\label{30.37}
J[{K^{*}}]_{T,t}^{(k)}
=J[\psi^{(k)}]_{T,t}+
\sum_{r=1}^{\left[k/2\right]}\frac{1}{2^r}
\sum_{(s_r,\ldots,s_1)\in {\rm A}_{k,r}}
J[\psi^{(k)}]_{T,t}^{s_r,\ldots,s_1}=J^{*}[\psi^{(k)}]_{T,t}.
\end{equation}

\vspace{5mm}

Let us consider the following generalized multiple Fourier sum

$$
\sum_{j_1=0}^{p_1}\ldots\sum_{j_k=0}^{p_k}
C_{j_k\ldots j_1} \prod_{l=1}^{k} \phi_{j_l}(t_l),
$$

\vspace{4mm}
\noindent
where $C_{j_k\ldots j_1}$ is the Fourier coefficient of the form

\vspace{-1mm}
\begin{equation}
\label{1}
C_{j_k\ldots j_1}=\int\limits_{[t,T]^k}
K^{*}(t_1,\ldots,t_k)\prod_{l=1}^{k}\phi_{j_l}(t_l)dt_1\ldots dt_k.
\end{equation}

\vspace{3mm}

Let us subsitute the relation

$$
K^{*}(t_1,\ldots,t_k)=
\sum_{j_1=0}^{p_1}\ldots\sum_{j_k=0}^{p_k}
C_{j_k\ldots j_1} \prod_{l=1}^{k} \phi_{j_l}(t_l)
+K^{*}(t_1,\ldots,t_k)-\sum_{j_1=0}^{p_1}\ldots\sum_{j_k=0}^{p_k}
C_{j_k\ldots j_1} \prod_{l=1}^{k} \phi_{j_l}(t_l)
$$

\vspace{3mm}
\noindent
into $J[{K^{*}}]_{T,t}^{(k)}$ (here $p_1,\ldots,p_k<\infty).$

Then, using Lemma 3, we obtain

\begin{equation}
\label{proof1}
J^{*}[\psi^{(k)}]_{T,t}=
\sum_{j_1=0}^{p_1}\ldots\sum_{j_k=0}^{p_k}
C_{j_k\ldots j_1}
\prod_{l=1}^k \zeta_{j_l}^{(i_l)}+
J[R_{p_1\ldots p_k}]_{T,t}^{(k)}\ \ \ \hbox{w. p. {\rm 1}},
\end{equation}

\vspace{3mm}
\noindent
where the stochastic integral
$J[R_{p_1\ldots p_k}]_{T,t}^{(k)}$
is defined in accordance with (\ref{30.34}) and

\begin{equation}
\label{30.46}
R_{p_1\ldots p_k}(t_1,\ldots,t_k)=
K^{*}(t_1,\ldots,t_k)-
\sum_{j_1=0}^{p_1}\ldots\sum_{j_k=0}^{p_k}
C_{j_k\ldots j_1} \prod_{l=1}^{k} \phi_{j_l}(t_l),
\end{equation}

\vspace{3mm}
\noindent
where
$$
\zeta_{j}^{(i)}=\int\limits_t^T \phi_{j}(s) d{\bf w}_s^{(i)}.
$$

\vspace{5mm}

\section{Proof of Theorem 4 Using the Combined Approach}

\vspace{5mm}

From (\ref{proof1}) we obtain

\vspace{-1mm}
\begin{equation}
\label{leto9000ddd}
J^{*}[\psi^{(2)}]_{T,t}=
\sum_{j_1=0}^{p_1}\sum_{j_2=0}^{p_2}C_{j_2 j_1}
\zeta_{j_1}^{(i_1)}\zeta_{j_2}^{(i_2)}+
J[R_{p_1p_2}]_{T,t}^{(2)}\ \ \ \hbox{w. p. {\rm 1}},
\end{equation}

\vspace{2mm}
\noindent
where

$$
J[R_{p_1p_2}]_{T,t}^{(2)}=\int\limits_t^T\int\limits_t^{t_2}
R_{p_1p_2}(t_1,t_2)d{\bf w}_{t_1}^{(i_1)}d{\bf w}_{t_2}^{(i_2)}
+\int\limits_t^T\int\limits_t^{t_1}
R_{p_1p_2}(t_1,t_2)d{\bf w}_{t_2}^{(i_2)}d{\bf w}_{t_1}^{(i_1)}+
$$

$$
+{\bf 1}_{\{i_1=i_2\ne 0\}}
\int\limits_t^T R_{p_1p_2}(t_1,t_1)dt_1,
$$

\vspace{2mm}
$$
R_{p_1 p_2}(t_1,t_2)=
K^{*}(t_1,t_2)-
\sum_{j_1=0}^{p_1}\sum_{j_2=0}^{p_2}C_{j_2 j_1}
\phi_{j_1}(t_1)\phi_{j_2}(t_2)\ \ \ (p_1, p_2<\infty),
$$

\vspace{4mm}

$$
K^{*}(t_1,t_2)=K(t_1,t_2)+\frac{1}{2}{\bf 1}_{\{t_1=t_2\}}
\psi_1(t_1)\psi_2(t_1),
$$

\vspace{6mm}
\noindent
where 
$$
K(t_1,t_2)=
\begin{cases}
\psi_1(t_1)\psi_2(t_2),\ &t_1<t_2\cr\cr
0,\ &\hbox{\rm otherwise}
\end{cases},\ \ \ t_1, t_2\in[t, T].
$$

\vspace{4mm}

Let us consider the case $i_1, i_2\ne 0$ (another cases can be considered
absolutely analogously).
Using standard estimates for 
moments of stochastic integrals \cite{1}, we get

\vspace{3mm}
$$
{\sf M}\left\{\left(J[R_{p_1p_2}]_{T,t}^{(2)}\right)^{2}
\right\}=
$$

\vspace{2mm}
$$
={\sf M}\left\{\left(\int\limits_t^T\int\limits_t^{t_2}
R_{p_1p_2}(t_1,t_2)d{\bf w}_{t_1}^{(i_1)}d{\bf w}_{t_2}^{(i_2)}
+\int\limits_t^T\int\limits_t^{t_1}
R_{p_1p_2}(t_1,t_2)d{\bf w}_{t_2}^{(i_2)}d{\bf w}_{t_1}^{(i_1)}
\right)^2\right\}+
$$

\vspace{2mm}
$$
+{\bf 1}_{\{i_1=i_2\ne 0\}}
\left(\int\limits_t^T R_{p_1p_2}(t_1,t_1)dt_1\right)^2\le
$$

\vspace{2mm}
$$
\le 
2\left(\int\limits_t^T\int\limits_t^{t_2}
\left(R_{p_1p_2}(t_1,t_2)\right)^{2}dt_1 dt_2
+
\int\limits_t^T\int\limits_t^{t_1}
\left(R_{p_1p_2}(t_1,t_2)\right)^{2}dt_2 dt_1\right)+
$$

\vspace{2mm}
$$
+
{\bf 1}_{\{i_1=i_2\ne 0\}}
\left(\int\limits_t^T R_{p_1p_2}(t_1,t_1)dt_1\right)^2=
$$

\vspace{2mm}
\begin{equation}
\label{newbegin1}
=
2\int\limits_{[t, T]^2}
\left(R_{p_1p_2}(t_1,t_2)\right)^{2}dt_1 dt_2
+
{\bf 1}_{\{i_1=i_2\ne 0\}}
\left(\int\limits_t^T R_{p_1p_2}(t_1,t_1)dt_1\right)^2.
\end{equation}

\vspace{6mm}

Moreover, we have

\vspace{-2mm}
$$
\int\limits_{[t, T]^2}
\left(R_{p_1p_2}(t_1,t_2)\right)^{2}dt_1 dt_2=
$$

\vspace{2mm}
$$
=
\int\limits_{[t, T]^2}
\Biggl(
K^{*}(t_1,t_2)-
\sum_{j_1=0}^{p_1}\sum_{j_2=0}^{p_2}C_{j_2 j_1}
\phi_{j_1}(t_1)\phi_{j_2}(t_2)\Biggr)^2 dt_1 dt_2=
$$

\vspace{2mm}
$$
=\int\limits_{[t, T]^2}
\Biggl(
K(t_1,t_2)-
\sum_{j_1=0}^{p_1}\sum_{j_2=0}^{p_2}C_{j_2 j_1}
\phi_{j_1}(t_1)\phi_{j_2}(t_2)\Biggr)^2 dt_1 dt_2.
$$

\vspace{5mm}

The function $K(t_1,t_2)$ is piecewise continuous in the 
square $[t, T]^2$.
At this situation it is well known that the
generalized multiple Fourier series 
of the function $K(t_1,t_2)\in L_2([t, T]^2)$ is converging 
to this function in the square $[t, T]^2$ in the mean-square sense, i.e.

\vspace{2mm}
$$
\hbox{\vtop{\offinterlineskip\halign{
\hfil#\hfil\cr
{\rm lim}\cr
$\stackrel{}{{}_{p_1,p_2\to \infty}}$\cr
}} }\Biggl\Vert
K(t_1,t_2)-
\sum_{j_1=0}^{p_1}\sum_{j_2=0}^{p_k}
C_{j_2 j_1}\prod_{l=1}^{2} \phi_{j_l}(t_l)\Biggr\Vert_{L_2([t,T]^2)}=0,
$$

\vspace{3mm}
\noindent
where
$$
\left\Vert f\right\Vert_{L_2([t,T]^2)}=\left(\int\limits_{[t,T]^2}
f^2(t_1,t_2)dt_1dt_2\right)^{1/2}.
$$

\vspace{2mm}

So, we obtain

\vspace{-1mm}
\begin{equation}
\label{newbegin2}
\hbox{\vtop{\offinterlineskip\halign{
\hfil#\hfil\cr
{\rm lim}\cr
$\stackrel{}{{}_{p_1,p_2\to \infty}}$\cr
}} }
\int\limits_{[t, T]^2}
\left(R_{p_1p_2}(t_1,t_2)\right)^{2}dt_1 dt_2=0.
\end{equation}

\vspace{3mm}

Note that

$$
\int\limits_t^T R_{p_1p_2}(t_1,t_1)dt_1=
$$

\vspace{2mm}
$$
=
\int\limits_t^T
\Biggl(
\frac{1}{2}\psi_1(t_1)\psi_2(t_1) - 
\sum_{j_1=0}^{p_1}\sum_{j_2=0}^{p_2}C_{j_2 j_1}
\phi_{j_1}(t_1)\phi_{j_2}(t_1)\Biggr) dt_1=
$$

\vspace{2mm}
$$
=
\frac{1}{2}\int\limits_t^T
\psi_1(t_1)\psi_2(t_1)dt_1 -
\sum_{j_1=0}^{p_1}\sum_{j_2=0}^{p_2}C_{j_2 j_1}
\int\limits_t^T\phi_{j_1}(t_1)\phi_{j_2}(t_1)dt_1=
$$

\vspace{2mm}
$$
=
\frac{1}{2}\int\limits_t^T
\psi_1(t_1)\psi_2(t_1)dt_1 -
\sum_{j_1=0}^{p_1}\sum_{j_2=0}^{p_2}C_{j_2 j_1}
{\bf 1}_{\{j_1=j_2\}}=
$$

\vspace{2mm}
\begin{equation}
\label{newbegin3}
=
\frac{1}{2}\int\limits_t^T
\psi_1(t_1)\psi_2(t_1)dt_1 -
\sum_{j_1=0}^{{\rm min}\{p_1,p_2\}}C_{j_1 j_1}.
\end{equation}

\vspace{6mm}

In \cite{16} (Theorem 3, p.~A.59), 
\cite{20} (Theorem 5.3, p.~A.294), \cite{20xx}-\cite{20xxz} (Theorems 2.1, 2.2), 
\cite{32} (Theorem 2), \cite{33} (Theorem 6) the following equality 

\vspace{-1mm}
\begin{equation}
\label{newbegin4}
\frac{1}{2}\int\limits_t^T
\psi_1(t_1)\psi_2(t_1)dt_1 =
\sum_{j_1=0}^{\infty}C_{j_1 j_1}
\end{equation}

\vspace{3mm}
\noindent
is proved. Note that the existence of the limit on the right-hand side of (\ref{newbegin4})
is proved in \cite{20xx}-\cite{20xxz}, \cite{32} for the polynomial and trigonometric cases.

From (\ref{newbegin1})--(\ref{newbegin4}) it follows that 

$$
\hbox{\vtop{\offinterlineskip\halign{
\hfil#\hfil\cr
{\rm lim}\cr
$\stackrel{}{{}_{p_1,p_2\to \infty}}$\cr
}} }
{\sf M}\left\{\left(J[R_{p_1p_2}]_{T,t}^{(2)}\right)^{2}
\right\}=0.
$$

\vspace{4mm}

Theorem 4 is proved.

\vspace{5mm}

\section{Proof of Theorem 5 Using the Combined Approach}

\vspace{5mm}

Let us consider (\ref{proof1}) for
$k=3$ and $p_1=p_2=p_3=p$

\vspace{2mm}
\begin{equation}
\label{newbegin97}
J^{*}[\psi^{(3)}]_{T,t}=
\sum_{j_1=0}^{p}\sum_{j_2=0}^{p}\sum_{j_3=0}^{p}C_{j_3j_2 j_1}
\zeta_{j_1}^{(i_1)}\zeta_{j_2}^{(i_2)}\zeta_{j_3}^{(i_3)}
+J[R_{ppp}]_{T,t}^{(3)}\ \ \ \hbox{w.\ p.\ 1,}
\end{equation}

\vspace{4mm}
\noindent
where

$$
J[R_{ppp}]_{T,t}^{(3)}=
\hbox{\vtop{\offinterlineskip\halign{
\hfil#\hfil\cr
{\rm l.i.m.}\cr
$\stackrel{}{{}_{N\to \infty}}$\cr
}} }\sum_{l_3=0}^{N-1}\sum_{l_2=0}^{N-1}
\sum_{l_1=0}^{N-1}
R_{ppp}(\tau_{l_1},\tau_{l_2},\tau_{l_3})
\Delta{\bf f}_{\tau_{l_1}}^{(i_1)}
\Delta{\bf f}_{\tau_{l_2}}^{(i_2)}
\Delta{\bf f}_{\tau_{l_3}}^{(i_3)},
$$

\vspace{3mm}

$$
R_{ppp}(t_1,t_2,t_3)
\stackrel{\small{\sf def}}{=}K^{*}(t_1,t_2,t_3)-
\sum_{j_1=0}^{p}\sum_{j_2=0}^{p}\sum_{j_3=0}^{p}
C_{j_3j_2 j_1}\phi_{j_1}(t_1)\phi_{j_2}(t_2)\phi_{j_3}(t_3),
$$

\vspace{5mm}

$$
K^{*}(t_1,t_2,t_3)
=\prod_{l=1}^{3}\psi_l(t_l)\Biggl(
{\bf 1}_{\{t_1<t_2\}}{\bf 1}_{\{t_2<t_3\}}+
\frac{1}{2}{\bf 1}_{\{t_1=t_2\}}{\bf 1}_{\{t_2<t_3\}}+\Biggr.
$$

\vspace{1mm}
$$
\Biggl.
+\frac{1}{2}{\bf 1}_{\{t_1<t_2\}}{\bf 1}_{\{t_2=t_3\}}+
\frac{1}{4}{\bf 1}_{\{t_1=t_2\}}{\bf 1}_{\{t_2=t_3\}}\Biggr).
$$

\vspace{6mm}

Furthermore, we have w.~p.~1

\vspace{2mm}
$$
J[R_{ppp}]_{T,t}^{(3)}=
\hbox{\vtop{\offinterlineskip\halign{
\hfil#\hfil\cr
{\rm l.i.m.}\cr
$\stackrel{}{{}_{N\to \infty}}$\cr
}} }\sum_{l_3=0}^{N-1}\sum_{l_2=0}^{N-1}
\sum_{l_1=0}^{N-1}
R_{ppp}(\tau_{l_1},\tau_{l_2},\tau_{l_3})\Delta{\bf f}_{\tau_{l_1}}^{(i_1)}
\Delta{\bf f}_{\tau_{l_2}}^{(i_2)}
\Delta{\bf f}_{\tau_{l_3}}^{(i_3)}=
$$

\vspace{2mm}
$$
=\hbox{\vtop{\offinterlineskip\halign{
\hfil#\hfil\cr
{\rm l.i.m.}\cr
$\stackrel{}{{}_{N\to \infty}}$\cr
}} }\sum_{l_3=0}^{N-1}\sum_{l_2=0}^{l_3-1}
\sum_{l_1=0}^{l_2-1}\Biggl(
R_{ppp}(\tau_{l_1},\tau_{l_2},\tau_{l_3})
\Delta{\bf f}_{\tau_{l_1}}^{(i_1)}
\Delta{\bf f}_{\tau_{l_2}}^{(i_2)}
\Delta{\bf f}_{\tau_{l_3}}^{(i_3)}+\Biggr.
$$

\vspace{2mm}
$$
+R_{ppp}(\tau_{l_1},\tau_{l_3},\tau_{l_2})
\Delta{\bf f}_{\tau_{l_1}}^{(i_1)}
\Delta{\bf f}_{\tau_{l_3}}^{(i_2)}
\Delta{\bf f}_{\tau_{l_2}}^{(i_3)}+
$$

\vspace{2mm}
$$
+
R_{ppp}(\tau_{l_2},\tau_{l_1},\tau_{l_3})
\Delta{\bf f}_{\tau_{l_2}}^{(i_1)}
\Delta{\bf f}_{\tau_{l_1}}^{(i_2)}
\Delta{\bf f}_{\tau_{l_3}}^{(i_3)}+
$$

\vspace{2mm}
$$
+R_{ppp}(\tau_{l_2},\tau_{l_3},\tau_{l_1})
\Delta{\bf f}_{\tau_{l_2}}^{(i_1)}
\Delta{\bf f}_{\tau_{l_3}}^{(i_2)}
\Delta{\bf f}_{\tau_{l_1}}^{(i_3)}+
$$

\vspace{2mm}
$$
+
R_{ppp}(\tau_{l_3},\tau_{l_2},\tau_{l_1})
\Delta{\bf f}_{\tau_{l_3}}^{(i_1)}
\Delta{\bf f}_{\tau_{l_2}}^{(i_2)}
\Delta{\bf f}_{\tau_{l_1}}^{(i_3)}+
$$

\vspace{2mm}
$$
\Biggl.
+R_{ppp}(\tau_{l_3},\tau_{l_1},\tau_{l_2})\Delta{\bf f}_{\tau_{l_3}}^{(i_1)}
\Delta{\bf f}_{\tau_{l_1}}^{(i_2)}
\Delta{\bf f}_{\tau_{l_2}}^{(i_3)}\Biggr)+
$$

\vspace{2mm}
$$
+
\hbox{\vtop{\offinterlineskip\halign{
\hfil#\hfil\cr
{\rm l.i.m.}\cr
$\stackrel{}{{}_{N\to \infty}}$\cr
}} }\sum_{l_3=0}^{N-1}\sum_{l_2=0}^{l_3-1}\Biggl(
R_{ppp}(\tau_{l_2},\tau_{l_2},\tau_{l_3})
\Delta{\bf f}_{\tau_{l_2}}^{(i_1)}
\Delta{\bf f}_{\tau_{l_2}}^{(i_2)}
\Delta{\bf f}_{\tau_{l_3}}^{(i_3)}+\Biggr.
$$

\vspace{2mm}
$$
+R_{ppp}(\tau_{l_2},\tau_{l_3},\tau_{l_2})
\Delta{\bf f}_{\tau_{l_2}}^{(i_1)}
\Delta{\bf f}_{\tau_{l_3}}^{(i_2)}
\Delta{\bf f}_{\tau_{l_2}}^{(i_3)}+
$$

\vspace{2mm}
$$
\Biggl.+R_{ppp}(\tau_{l_3},\tau_{l_2},\tau_{l_2})
\Delta{\bf f}_{\tau_{l_3}}^{(i_1)}
\Delta{\bf f}_{\tau_{l_2}}^{(i_2)}
\Delta{\bf f}_{\tau_{l_2}}^{(i_3)}\Biggr)+
$$

\vspace{2mm}
$$
+
\hbox{\vtop{\offinterlineskip\halign{
\hfil#\hfil\cr
{\rm l.i.m.}\cr
$\stackrel{}{{}_{N\to \infty}}$\cr
}} }\sum_{l_3=0}^{N-1}\sum_{l_1=0}^{l_3-1}\Biggl(
R_{ppp}(\tau_{l_1},\tau_{l_3},\tau_{l_3})
\Delta{\bf f}_{\tau_{l_1}}^{(i_1)}
\Delta{\bf f}_{\tau_{l_3}}^{(i_2)}
\Delta{\bf f}_{\tau_{l_3}}^{(i_3)}+\Biggr.
$$

\vspace{2mm}
$$
+R_{ppp}(\tau_{l_3},\tau_{l_1},\tau_{l_3})
\Delta{\bf f}_{\tau_{l_3}}^{(i_1)}
\Delta{\bf f}_{\tau_{l_1}}^{(i_2)}
\Delta{\bf f}_{\tau_{l_3}}^{(i_3)}+
$$

\vspace{2mm}
$$
\Biggl.+R_{ppp}(\tau_{l_3},\tau_{l_3},\tau_{l_1})
\Delta{\bf f}_{\tau_{l_3}}^{(i_1)}
\Delta{\bf f}_{\tau_{l_3}}^{(i_2)}
\Delta{\bf f}_{\tau_{l_1}}^{(i_3)}\Biggr)+
$$

\vspace{2mm}
$$
+\hbox{\vtop{\offinterlineskip\halign{
\hfil#\hfil\cr
{\rm l.i.m.}\cr
$\stackrel{}{{}_{N\to \infty}}$\cr
}} }\sum_{l_3=0}^{N-1}
R_{ppp}(\tau_{l_3},\tau_{l_3},\tau_{l_3})
\Delta{\bf f}_{\tau_{l_3}}^{(i_1)}
\Delta{\bf f}_{\tau_{l_3}}^{(i_2)}
\Delta{\bf f}_{\tau_{l_3}}^{(i_3)}=
$$

\vspace{2mm}
$$
=R_{T,t}^{(1)ppp}+R_{T,t}^{(2)ppp},
$$

\vspace{8mm}
\noindent
where

$$
R_{T,t}^{(1)ppp}=
$$

\vspace{2mm}
$$
=\int\limits_t^T\int\limits_t^{t_3}\int\limits_t^{t_2}
R_{ppp}(t_1,t_2,t_3)
d{\bf f}_{t_1}^{(i_1)}
d{\bf f}_{t_2}^{(i_2)}
d{\bf f}_{t_3}^{(i_3)}+
\int\limits_t^T\int\limits_t^{t_3}\int\limits_t^{t_2}
R_{ppp}(t_1,t_3,t_2)
d{\bf f}_{t_1}^{(i_1)}
d{\bf f}_{t_2}^{(i_3)}
d{\bf f}_{t_3}^{(i_2)}+
$$

\vspace{2mm}
$$
+
\int\limits_t^T\int\limits_t^{t_3}\int\limits_t^{t_2}
R_{ppp}(t_2,t_1,t_3)
d{\bf f}_{t_1}^{(i_2)}
d{\bf f}_{t_2}^{(i_1)}
d{\bf f}_{t_3}^{(i_3)}+
\int\limits_t^T\int\limits_t^{t_3}\int\limits_t^{t_2}
R_{ppp}(t_2,t_3,t_1)
d{\bf f}_{t_1}^{(i_3)}
d{\bf f}_{t_2}^{(i_1)}
d{\bf f}_{t_3}^{(i_2)}+
$$

\vspace{2mm}
$$
+
\int\limits_t^T\int\limits_t^{t_3}\int\limits_t^{t_2}
R_{ppp}(t_3,t_2,t_1)
d{\bf f}_{t_1}^{(i_3)}
d{\bf f}_{t_2}^{(i_2)}
d{\bf f}_{t_3}^{(i_1)}+
\int\limits_t^T\int\limits_t^{t_3}\int\limits_t^{t_2}
R_{ppp}(t_3,t_1,t_2)
d{\bf f}_{t_1}^{(i_2)}
d{\bf f}_{t_2}^{(i_3)}
d{\bf f}_{t_3}^{(i_1)},
$$

\vspace{7mm}

$$
R_{T,t}^{(2)ppp}=
$$

\vspace{1mm}
$$
={\bf 1}_{\{i_1=i_2\ne 0\}}
\int\limits_t^T\int\limits_t^{t_3}
R_{ppp}(t_2,t_2,t_3)
dt_2
d{\bf f}_{t_3}^{(i_3)}+
{\bf 1}_{\{i_1=i_3\ne 0\}}
\int\limits_t^T\int\limits_t^{t_3}
R_{ppp}(t_2,t_3,t_2)
dt_2
d{\bf f}_{t_3}^{(i_2)}+
$$

\vspace{2mm}
$$
+{\bf 1}_{\{i_2=i_3\ne 0\}}
\int\limits_t^T\int\limits_t^{t_3}
R_{ppp}(t_3,t_2,t_2)
dt_2
d{\bf f}_{t_3}^{(i_1)}+
{\bf 1}_{\{i_2=i_3\ne 0\}}
\int\limits_t^T\int\limits_t^{t_3}
R_{ppp}(t_1,t_3,t_3)
d{\bf f}_{t_1}^{(i_1)}dt_3+
$$

\vspace{2mm}
$$
+{\bf 1}_{\{i_1=i_3\ne 0\}}
\int\limits_t^T\int\limits_t^{t_3}
R_{ppp}(t_3,t_1,t_3)
d{\bf f}_{t_1}^{(i_2)}dt_3+
{\bf 1}_{\{i_1=i_2\ne 0\}}
\int\limits_t^T\int\limits_t^{t_3}
R_{ppp}(t_3,t_3,t_1)
d{\bf f}_{t_1}^{(i_3)}dt_3.
$$

\vspace{8mm}

Moreover, we obtain

\begin{equation}
\label{newbegin98}
{\sf M}\left\{\left(J[R_{ppp}]_{T,t}^{(3)}\right)^2\right\}\le
2{\sf M}\left\{\left(R_{T,t}^{(1)ppp}\right)^2\right\}
+2{\sf M}\left\{\left(R_{T,t}^{(2)ppp}\right)^2\right\}.
\end{equation}

\vspace{6mm}

Now, using standard estimates for moments of stochastic 
integrals \cite{1}, we obtain the following inequality

\vspace{1mm}
$$
{\sf M}\left\{\left(R_{T,t}^{(1)ppp}\right)^{2}
\right\}\le 
$$

\vspace{2mm}
$$
\le 6 \int\limits_t^T\int\limits_t^{t_3}\int\limits_t^{t_2}
\Biggl(
\left(R_{p_1 p_2 p_3}(t_1,t_2,t_3)\right)^{2}+
\left(R_{p_1 p_2 p_3}(t_1,t_3,t_2)\right)^{2}+\Biggr.
\left(R_{p_1 p_2 p_3}(t_2,t_1,t_3)\right)^{2}+
$$

\vspace{2mm}
$$
+\left(R_{p_1 p_2 p_3}(t_2,t_3,t_1)\right)^{2}+
\left(R_{p_1 p_2 p_3}(t_3,t_2,t_1)\right)^{2}+
\Biggl.\left(R_{p_1 p_2 p_3}(t_3,t_1,t_2)\right)^{2}\Biggr)dt_1dt_2dt_3=
$$

\vspace{4mm}
$$
=6\int\limits_{[t, T]^3}
\left(R_{ppp}(t_1,t_2,t_3)\right)^{2}dt_1 dt_2 dt_3.
$$

\vspace{6mm}

We have

$$
\int\limits_{[t, T]^3}
\left(R_{ppp}(t_1,t_2,t_3)\right)^{2}dt_1 dt_2 dt_3=
$$

\vspace{2mm}
$$
=
\int\limits_{[t, T]^3}
\Biggl(
K^{*}(t_1,t_2,t_3)-
\sum_{j_1=0}^{p}\sum_{j_2=0}^{p}\sum_{j_3=0}^{p}C_{j_3 j_2 j_1}
\phi_{j_1}(t_1)\phi_{j_2}(t_2)\phi_{j_3}(t_3)\Biggr)^2 
dt_1 dt_2 dt_3=
$$

\vspace{2mm}
$$
=
\int\limits_{[t, T]^3}
\Biggl(
K(t_1,t_2,t_3)-
\sum_{j_1=0}^{p}\sum_{j_2=0}^{p}\sum_{j_3=0}^{p}C_{j_3 j_2 j_1}
\phi_{j_1}(t_1)\phi_{j_2}(t_2)\phi_{j_3}(t_3)\Biggr)^2 dt_1 dt_2 dt_3,
$$

\vspace{7mm}
\noindent
where

\vspace{-1mm}
$$
K(t_1,t_2,t_3)=
\begin{cases}
\psi_1(t_1)\psi_2(t_2)\psi_3(t_3),\ &t_1<t_2<t_3\cr\cr
0,\ &\hbox{\rm otherwise}
\end{cases},\ \ \ t_1, t_2, t_3\in[t, T].
$$

\vspace{7mm}

So, we get

\vspace{-1mm}
$$
\hbox{\vtop{\offinterlineskip\halign{
\hfil#\hfil\cr
{\rm lim}\cr
$\stackrel{}{{}_{p\to \infty}}$\cr
}} }
{\sf M}\left\{\left(R_{T,t}^{(1)ppp}\right)^{2}\right\}\le
$$

\vspace{2mm}
\begin{equation}
\label{newbegin99}
\le 6\hbox{\vtop{\offinterlineskip\halign{
\hfil#\hfil\cr
{\rm lim}\cr
$\stackrel{}{{}_{p\to \infty}}$\cr
}} }
\int\limits_{[t, T]^3}
\Biggl(
K(t_1,t_2,t_3)-
\sum_{j_1=0}^{p}\sum_{j_2=0}^{p}\sum_{j_3=0}^{p}C_{j_3 j_2 j_1}
\phi_{j_1}(t_1)\phi_{j_2}(t_2)\phi_{j_3}(t_3)\Biggr)^2 dt_1 dt_2 dt_3=0,
\end{equation}

\vspace{5mm}
\noindent
where
$K(t_1,t_2,t_3) \in L_2([t, T]^3)$.

After the integration order replacement 
in iterated Ito stochastic integrals \cite{th} (also see
\cite{16}, \cite{20} or Chapter 3 in \cite{20xx}-\cite{20xxz})
from $R_{T,t}^{(2)ppp}$ we obtain w.~p.~1

\vspace{2mm}
$$
R_{T,t}^{(2)ppp}=
$$

\vspace{2mm}
$$
={\bf 1}_{\{i_1=i_2\ne 0\}}\left(
\int\limits_t^T\int\limits_t^{t_3}
R_{ppp}(t_2,t_2,t_3)
dt_2
d{\bf f}_{t_3}^{(i_3)}+\Biggr.
\Biggl.
\int\limits_t^T\int\limits_t^{t_3}
R_{ppp}(t_3,t_3,t_1)
d{\bf f}_{t_1}^{(i_3)}dt_3\right)+
$$

\vspace{3mm}
$$
+{\bf 1}_{\{i_2=i_3\ne 0\}}\left(
\int\limits_t^T\int\limits_t^{t_3}
R_{ppp}(t_3,t_2,t_2)
dt_2
d{\bf f}_{t_3}^{(i_1)}+\Biggr.
\Biggl.
\int\limits_t^T\int\limits_t^{t_3}
R_{ppp}(t_1,t_3,t_3)
d{\bf f}_{t_1}^{(i_1)}dt_3\right)+
$$

\vspace{3mm}
$$
+{\bf 1}_{\{i_1=i_3\ne 0\}}\left(
\int\limits_t^T\int\limits_t^{t_3}
R_{ppp}(t_2,t_3,t_2)dt_2
d{\bf f}_{t_3}^{(i_2)}+\Biggr.
\Biggl.
\int\limits_t^T\int\limits_t^{t_3}
R_{ppp}(t_3,t_1,t_3)
d{\bf f}_{t_1}^{(i_2)}dt_3\right)=
$$

\vspace{3mm}
$$
={\bf 1}_{\{i_1=i_2\ne 0\}}\left(
\int\limits_t^T\int\limits_t^{t_1}
R_{ppp}(t_2,t_2,t_1)
dt_2
d{\bf f}_{t_1}^{(i_3)}+\Biggr.
\Biggl.
\int\limits_t^T\int\limits_{t_1}^{T}
R_{ppp}(t_2,t_2,t_1)
dt_2d{\bf f}_{t_1}^{(i_3)}\right)+
$$

\vspace{3mm}
$$
+{\bf 1}_{\{i_2=i_3\ne 0\}}\left(
\int\limits_t^T\int\limits_t^{t_1}
R_{ppp}(t_1,t_2,t_2)
dt_2
d{\bf f}_{t_1}^{(i_1)}+\Biggr.
\Biggl.
\int\limits_t^T\int\limits_{t_1}^{T}
R_{ppp}(t_1,t_2,t_2)
dt_2d{\bf f}_{t_1}^{(i_1)}\right)+
$$

\vspace{3mm}
$$
+{\bf 1}_{\{i_1=i_3\ne 0\}}\left(
\int\limits_t^T\int\limits_t^{t_1}
R_{ppp}(t_2,t_1,t_2)dt_2
d{\bf f}_{t_1}^{(i_2)}+\Biggr.
\Biggl.
\int\limits_t^T\int\limits_{t_1}^{T}
R_{ppp}(t_2,t_1,t_2)
dt_2d{\bf f}_{t_1}^{(i_2)}\right)=
$$

\vspace{3mm}
$$
={\bf 1}_{\{i_1=i_2\ne 0\}}
\int\limits_t^T\left(\int\limits_t^{T}
R_{ppp}(t_2,t_2,t_3)dt_2\right)
d{\bf f}_{t_3}^{(i_3)}+
$$

\vspace{3mm}
$$
+
{\bf 1}_{\{i_2=i_3\ne 0\}}
\int\limits_t^T\left(\int\limits_t^{T}
R_{ppp}(t_1,t_2,t_2)dt_2\right)
d{\bf f}_{t_1}^{(i_1)}+
$$

\vspace{3mm}
$$
+{\bf 1}_{\{i_1=i_3\ne 0\}}
\int\limits_t^T\left(\int\limits_t^{T}
R_{ppp}(t_3,t_2,t_3)dt_3\right)
d{\bf f}_{t_2}^{(i_2)}=
$$

\vspace{6mm}
$$
={\bf 1}_{\{i_1=i_2\ne 0\}}
\int\limits_t^T\int\limits_t^{T}\Biggl(\Biggl(
\frac{1}{2}{\bf 1}_{\{t_2<t_3\}} + \frac{1}{4}{\bf 1}_{\{t_2=t_3\}}
\Biggr)
\psi_1(t_2)\psi_2(t_2)\psi_3(t_3)-\Biggr.
$$

\vspace{3mm}
$$
\Biggl.
-\sum_{j_1=0}^{p}\sum_{j_2=0}^{p}\sum_{j_3=0}^{p}
C_{j_3 j_2 j_1}\phi_{j_1}(t_2)\phi_{j_2}(t_2)\phi_{j_3}(t_3)
\Biggr)dt_2 d{\bf f}_{t_3}^{(i_3)}+
$$

\vspace{3mm}
$$
+{\bf 1}_{\{i_2=i_3\ne 0\}}
\int\limits_t^T\int\limits_t^{T}\Biggl(\Biggl(
\frac{1}{2}{\bf 1}_{\{t_1<t_2\}} + \frac{1}{4}{\bf 1}_{\{t_1=t_2\}}
\Biggr)
\psi_1(t_1)\psi_2(t_2)\psi_3(t_2)-\Biggr.
$$

\vspace{3mm}
$$
\Biggl.
-\sum_{j_1=0}^{p}\sum_{j_2=0}^{p}\sum_{j_3=0}^{p}
C_{j_3 j_2 j_1}\phi_{j_1}(t_1)\phi_{j_2}(t_2)\phi_{j_3}(t_2)
\Biggr)dt_2 d{\bf f}_{t_1}^{(i_1)}+
$$

\vspace{3mm}
$$
+{\bf 1}_{\{i_1=i_3\ne 0\}}
\int\limits_t^T\int\limits_t^{T}\Biggl(
\frac{1}{4}{\bf 1}_{\{t_2=t_3\}}
\psi_1(t_3)\psi_2(t_2)\psi_3(t_3)-\Biggr.
$$

\vspace{3mm}
$$
\Biggl.
-\sum_{j_1=0}^{p}\sum_{j_2=0}^{p}\sum_{j_3=0}^{p}
C_{j_3 j_2 j_1}\phi_{j_1}(t_3)\phi_{j_2}(t_2)\phi_{j_3}(t_3)
\Biggr)dt_3 d{\bf f}_{t_2}^{(i_2)}=
$$

\vspace{6mm}
$$
={\bf 1}_{\{i_1=i_2\ne 0\}}
\int\limits_t^T\left(\frac{1}{2}\psi_3(t_3)\int\limits_t^{t_3}
\psi_1(t_2)\psi_2(t_2)dt_2-\Biggr.
\Biggl.
\sum_{j_1=0}^{p}\sum_{j_3=0}^{p}
C_{j_3 j_1 j_1}\phi_{j_3}(t_3)
\right)d{\bf f}_{t_3}^{(i_3)}+
$$

\vspace{3mm}
$$
+{\bf 1}_{\{i_2=i_3\ne 0\}}
\int\limits_t^T\left(\frac{1}{2}\psi_1(t_1)\int\limits_{t_1}^{T}
\psi_2(t_2)\psi_3(t_2)dt_2-\Biggr.
\Biggl.
\sum_{j_1=0}^{p}\sum_{j_3=0}^{p}
C_{j_3 j_3 j_1}\phi_{j_1}(t_1)
\right)d{\bf f}_{t_1}^{(i_1)}+
$$

\vspace{3mm}
$$
+{\bf 1}_{\{i_1=i_3\ne 0\}}
\int\limits_t^T (-1)\sum_{j_1=0}^{p}\sum_{j_2=0}^{p}
C_{j_1 j_2 j_1}\phi_{j_2}(t_2)
d{\bf f}_{t_2}^{(i_2)}=
$$

\vspace{6mm}
$$
={\bf 1}_{\{i_1=i_2\ne 0\}}\left(
\frac{1}{2}\int\limits_t^T\psi_3(t_3)\int\limits_t^{t_3}
\psi_1(t_2)\psi_2(t_2)dt_2d{\bf f}_{t_3}^{(i_3)}-
\sum_{j_1=0}^{p}\sum_{j_3=0}^{p}
C_{j_3 j_1 j_1}\zeta_{j_3}^{(i_3)}\right)+
$$

\vspace{3mm}
$$
+{\bf 1}_{\{i_2=i_3\ne 0\}}\left(
\frac{1}{2}\int\limits_t^T\psi_1(t_1)\int\limits_{t_1}^{T}
\psi_2(t_2)\psi_3(t_2)dt_2d{\bf f}_{t_1}^{(i_1)}-\Biggr.
\Biggl.
\sum_{j_1=0}^{p}\sum_{j_3=0}^{p}
C_{j_3 j_3 j_1}\zeta_{j_1}^{(i_1)}\right)-
$$

\vspace{3mm}
$$
-{\bf 1}_{\{i_1=i_3\ne 0\}}
\sum_{j_1=0}^{p}\sum_{j_3=0}^{p}
C_{j_1 j_3 j_1}\zeta_{j_3}^{(i_2)}.
$$

\vspace{8mm}

From \cite{16} (Theorem 6, pp.~A.116--A.117), 
\cite{20} (Theorem 5.5', p.~A.371), \cite{20xx}-\cite{20xxz} (Chapter 2),
\cite{32} (Theorem 3) we obtain

$$
{\sf M}\left\{\left(R_{T,t}^{(2)ppp}\right)^2\right\}\le
$$

\vspace{2mm}
$$
\le
3\left({\bf 1}_{\{i_1=i_2\ne 0\}}{\sf M}\left\{\left(
\frac{1}{2}\int\limits_t^T\psi_3(t_3)\int\limits_t^{t_3}
\psi_1(t_2)\psi_2(t_2)dt_2d{\bf f}_{t_3}^{(i_3)}-
\sum_{j_1=0}^{p}\sum_{j_3=0}^{p}
C_{j_3 j_1 j_1}\zeta_{j_3}^{(i_3)}\right)^2\right\}+\right.
$$

\vspace{3mm}
$$
+{\bf 1}_{\{i_1=i_3\ne 0\}}
{\sf M}\left\{\left(\sum_{j_1=0}^{p}\sum_{j_3=0}^{p}
C_{j_1 j_3 j_1}\zeta_{j_3}^{(i_2)}\right)^2\right\}+
$$

\vspace{3mm}
\begin{equation}
\label{newbegin100}
\left.+{\bf 1}_{\{i_2=i_3\ne 0\}}{\sf M}\left\{\left(
\frac{1}{2}\int\limits_t^T\psi_1(t_1)\int\limits_{t_1}^{T}
\psi_2(t_2)\psi_3(t_2)dt_2d{\bf f}_{t_1}^{(i_1)}-
\sum_{j_1=0}^{p}\sum_{j_3=0}^{p}
C_{j_3 j_3 j_1}\zeta_{j_1}^{(i_1)}\right)^2\right\}\right)\ \to 0\
\end{equation}

\vspace{6mm}
\noindent
if $p \to \infty$. Using (\ref{newbegin97})--(\ref{newbegin100}),
we obtain the expansion (\ref{newbegin96}). 
Theorem 5 is proved.

\vspace{7mm}

\section{Proof of Theorem 6 Using the Combined Approach}

\vspace{7mm}

Let us consider (\ref{proof1}) for the case $k=4,$
$p_1=p_2=p_3=p_4=p$, and $\psi_1(\tau),$ $\psi_2(\tau),$ $\psi_3(\tau),$
$\psi_4(\tau)\equiv 1$

\vspace{2mm}
$$
{\int\limits_t^{*}}^T
{\int\limits_t^{*}}^{t_4}
{\int\limits_t^{*}}^{t_3}
{\int\limits_t^{*}}^{t_2}
d{\bf w}_{t_1}^{(i_1)}
d{\bf w}_{t_2}^{(i_2)}d{\bf w}_{t_3}^{(i_3)}d{\bf w}_{t_4}^{(i_4)}
=
$$

\vspace{2mm}

\begin{equation}
\label{proof2}
=
\sum_{j_1=0}^{p}\sum_{j_2=0}^{p}\sum_{j_3=0}^{p}\sum_{j_4=0}^{p}
C_{j_4 j_3 j_2 j_1}
\zeta_{j_1}^{(i_1)}\zeta_{j_2}^{(i_2)}\zeta_{j_3}^{(i_3)}\zeta_{j_4}^{(i_4)}+
J[R_{pppp}]_{T,t}^{(4)}\ \ \ \hbox{w. p. 1},
\end{equation}

\vspace{4mm}
\noindent
where 

\vspace{-3mm}
$$
J[R_{pppp}]_{T,t}^{(4)}=
\hbox{\vtop{\offinterlineskip\halign{
\hfil#\hfil\cr
{\rm l.i.m.}\cr
$\stackrel{}{{}_{N\to \infty}}$\cr
}} }\sum_{l_4=0}^{N-1}\sum_{l_3=0}^{N-1}\sum_{l_2=0}^{N-1}
\sum_{l_1=0}^{N-1}
R_{pppp}(\tau_{l_1},\tau_{l_2},\tau_{l_3},\tau_{l_4})
\Delta{\bf w}_{\tau_{l_1}}^{(i_1)}
\Delta{\bf w}_{\tau_{l_2}}^{(i_2)}
\Delta{\bf w}_{\tau_{l_3}}^{(i_3)}
\Delta{\bf w}_{\tau_{l_4}}^{(i_4)},
$$

\vspace{4mm}

$$
R_{pppp}(t_1,t_2,t_3,t_4)
\stackrel{\small{\sf def}}{=}
K^{*}(t_1,t_2,t_3,t_4)-
$$

\vspace{1mm}
$$
-
\sum_{j_1=0}^{p}\sum_{j_2=0}^{p}
\sum_{j_3=0}^{p}\sum_{j_4=0}^{p}
C_{j_4j_3j_2 j_1}\phi_{j_1}(t_1)\phi_{j_2}(t_2)\phi_{j_3}(t_3)\phi_{j_4}(t_4),
$$

\vspace{5mm}

$$
K^{*}(t_1,t_2,t_3,t_4)\stackrel{\small{\sf def}}{=}
\prod_{l=1}^3\Biggl({\bf 1}_{\{t_l<t_{l+1}\}}+
\frac{1}{2}{\bf 1}_{\{t_l=t_{l+1}\}}\Biggr)=
$$

\vspace{1mm}
$$
={\bf 1}_{\{t_1<t_2<t_3<t_4\}}+
\frac{1}{2}{\bf 1}_{\{t_1=t_2<t_3<t_4\}}+
\frac{1}{2}{\bf 1}_{\{t_1<t_2=t_3<t_4\}}+
$$

\vspace{1mm}
$$
+
\frac{1}{4}{\bf 1}_{\{t_1=t_2=t_3<t_4\}}+
\frac{1}{2}{\bf 1}_{\{t_1<t_2<t_3=t_4\}}+
\frac{1}{4}{\bf 1}_{\{t_1=t_2<t_3=t_4\}}+
$$

\vspace{1mm}
$$
+
\frac{1}{4}{\bf 1}_{\{t_1<t_2=t_3=t_4\}}+
\frac{1}{8}{\bf 1}_{\{t_1=t_2=t_3=t_4\}}.
$$

\vspace{5mm}

Moreover, we have

\vspace{-1mm}
\begin{equation}
\label{klor1}
J[R_{pppp}]_{T,t}^{(4)}=\sum\limits_{i=0}^7 R_{T,t}^{(i)pppp}\ \ \ 
\hbox{w.~p.~1},
\end{equation}

\vspace{2mm}
\noindent
where

\vspace{1mm}
$$
R_{T,t}^{(0)pppp}=
\hbox{\vtop{\offinterlineskip\halign{
\hfil#\hfil\cr
{\rm l.i.m.}\cr
$\stackrel{}{{}_{N\to \infty}}$\cr
}} }\sum_{l_4=0}^{N-1}\sum_{l_3=0}^{l_4-1}\sum_{l_2=0}^{l_3-1}
\sum_{l_1=0}^{l_2-1}
\sum\limits_{(l_1,l_2,l_3,l_4)}
\Biggl(
R_{pppp}(\tau_{l_1},\tau_{l_2},\tau_{l_3},\tau_{l_4})\times\Biggr.
$$

\vspace{2mm}
$$
\Biggl.
\times
\Delta{\bf w}_{\tau_{l_1}}^{(i_1)}
\Delta{\bf w}_{\tau_{l_2}}^{(i_2)}
\Delta{\bf w}_{\tau_{l_3}}^{(i_3)}
\Delta{\bf w}_{\tau_{l_4}}^{(i_4)}\Biggr),
$$

\vspace{4mm}

\noindent
where 
permutations
$(l_1,l_2,l_3,l_4)$ when summing are
performed only in the expression, which is enclosed in parentheses,

\vspace{2mm}
$$
R_{T,t}^{(1)pppp}={\bf 1}_{\{i_1=i_2\ne 0\}}
\hbox{\vtop{\offinterlineskip\halign{
\hfil#\hfil\cr
{\rm l.i.m.}\cr
$\stackrel{}{{}_{N\to \infty}}$\cr
}} }
\sum_{\stackrel{l_4,l_3,l_1=0}{{}_{l_1\ne l_3, l_1\ne l_4,
l_3\ne l_4}}}^{N-1}
R_{pppp}(\tau_{l_1},\tau_{l_1},\tau_{l_3},\tau_{l_4})
\Delta\tau_{l_1}
\Delta{\bf w}_{\tau_{l_3}}^{(i_3)}
\Delta{\bf w}_{\tau_{l_4}}^{(i_4)},
$$

\vspace{2mm}
$$
R_{T,t}^{(2)pppp}={\bf 1}_{\{i_1=i_3\ne 0\}}
\hbox{\vtop{\offinterlineskip\halign{
\hfil#\hfil\cr
{\rm l.i.m.}\cr
$\stackrel{}{{}_{N\to \infty}}$\cr
}} }
\sum_{\stackrel{l_4,l_2,l_1=0}{{}_{l_1\ne l_2, l_1\ne l_4,
l_2\ne l_4}}}^{N-1}
R_{pppp}(\tau_{l_1},\tau_{l_2},\tau_{l_1},\tau_{l_4})
\Delta\tau_{l_1}
\Delta{\bf w}_{\tau_{l_2}}^{(i_2)}
\Delta{\bf w}_{\tau_{l_4}}^{(i_4)},
$$

\vspace{2mm}
$$
R_{T,t}^{(3)pppp}={\bf 1}_{\{i_1=i_4\ne 0\}}
\hbox{\vtop{\offinterlineskip\halign{
\hfil#\hfil\cr
{\rm l.i.m.}\cr
$\stackrel{}{{}_{N\to \infty}}$\cr
}} }
\sum_{\stackrel{l_3,l_2,l_1=0}{{}_{l_1\ne l_2, l_1\ne l_3,
l_2\ne l_3}}}^{N-1}
R_{pppp}(\tau_{l_1},\tau_{l_2},\tau_{l_3},\tau_{l_1})
\Delta\tau_{l_1}
\Delta{\bf w}_{\tau_{l_2}}^{(i_2)}
\Delta{\bf w}_{\tau_{l_3}}^{(i_3)},
$$

\vspace{2mm}
$$
R_{T,t}^{(4)pppp}={\bf 1}_{\{i_2=i_3\ne 0\}}
\hbox{\vtop{\offinterlineskip\halign{
\hfil#\hfil\cr
{\rm l.i.m.}\cr
$\stackrel{}{{}_{N\to \infty}}$\cr
}} }
\sum_{\stackrel{l_4,l_2,l_1=0}{{}_{l_1\ne l_2, l_1\ne l_4,
l_2\ne l_4}}}^{N-1}
R_{pppp}(\tau_{l_1},\tau_{l_2},\tau_{l_2},\tau_{l_4})
\Delta{\bf w}_{\tau_{l_1}}^{(i_1)}
\Delta\tau_{l_2}
\Delta{\bf w}_{\tau_{l_4}}^{(i_4)},
$$

\vspace{2mm}
$$
R_{T,t}^{(5)pppp}={\bf 1}_{\{i_2=i_4\ne 0\}}
\hbox{\vtop{\offinterlineskip\halign{
\hfil#\hfil\cr
{\rm l.i.m.}\cr
$\stackrel{}{{}_{N\to \infty}}$\cr
}} }
\sum_{\stackrel{l_3,l_2,l_1=0}{{}_{l_1\ne l_2, l_1\ne l_3,
l_2\ne l_3}}}^{N-1}
R_{pppp}(\tau_{l_1},\tau_{l_2},\tau_{l_3},\tau_{l_2})
\Delta{\bf w}_{\tau_{l_1}}^{(i_1)}
\Delta\tau_{l_2}
\Delta{\bf w}_{\tau_{l_3}}^{(i_3)},
$$

\vspace{2mm}
$$
R_{T,t}^{(6)pppp}={\bf 1}_{\{i_3=i_4\ne 0\}}
\hbox{\vtop{\offinterlineskip\halign{
\hfil#\hfil\cr
{\rm l.i.m.}\cr
$\stackrel{}{{}_{N\to \infty}}$\cr
}} }
\sum_{\stackrel{l_3,l_2,l_1=0}{{}_{l_1\ne l_2, l_1\ne l_3,
l_2\ne l_3}}}^{N-1}
R_{pppp}(\tau_{l_1},\tau_{l_2},\tau_{l_3},\tau_{l_3})
\Delta{\bf w}_{\tau_{l_1}}^{(i_1)}
\Delta{\bf w}_{\tau_{l_2}}^{(i_2)}
\Delta\tau_{l_3},
$$

\vspace{5mm}

$$
R_{T,t}^{(7)pppp}={\bf 1}_{\{i_1=i_2\ne 0\}}{\bf 1}_{\{i_3=i_4\ne 0\}}
\hbox{\vtop{\offinterlineskip\halign{
\hfil#\hfil\cr
{\rm lim}\cr
$\stackrel{}{{}_{N\to \infty}}$\cr
}} }
\sum_{\stackrel{l_4,l_2=0}{{}_{l_2\ne l_4}}}^{N-1}
R_{pppp}(\tau_{l_2},\tau_{l_2},\tau_{l_4},\tau_{l_4})
\Delta\tau_{l_2}
\Delta\tau_{l_4}+
$$

\vspace{2mm}
$$
+{\bf 1}_{\{i_1=i_3\ne 0\}}{\bf 1}_{\{i_2=i_4\ne 0\}}
\hbox{\vtop{\offinterlineskip\halign{
\hfil#\hfil\cr
{\rm lim}\cr
$\stackrel{}{{}_{N\to \infty}}$\cr
}} }
\sum_{\stackrel{l_4,l_2=0}{{}_{l_2\ne l_4}}}^{N-1}
R_{pppp}(\tau_{l_2},\tau_{l_4},\tau_{l_2},\tau_{l_4})
\Delta\tau_{l_2}
\Delta\tau_{l_4}+
$$

\vspace{2mm}
$$
+{\bf 1}_{\{i_1=i_4\ne 0\}}{\bf 1}_{\{i_2=i_3\ne 0\}}
\hbox{\vtop{\offinterlineskip\halign{
\hfil#\hfil\cr
{\rm lim}\cr
$\stackrel{}{{}_{N\to \infty}}$\cr
}} }
\sum_{\stackrel{l_4,l_2=0}{{}_{l_2\ne l_4}}}^{N-1}
R_{pppp}(\tau_{l_2},\tau_{l_4},\tau_{l_4},\tau_{l_2})
\Delta\tau_{l_2}
\Delta\tau_{l_4}.
$$

\vspace{8mm}

The relations  (\ref{proof2}) and (\ref{klor1}) imply that Theorem 6
will be proved if 

\vspace{2mm}
$$
\hbox{\vtop{\offinterlineskip\halign{
\hfil#\hfil\cr
{\rm lim}\cr
$\stackrel{}{{}_{p\to \infty}}$\cr
}} }
{\sf M}\left\{\left(R_{T,t}^{(i)pppp}\right)^2\right\}=0,\ \ \ 
i=0, 1, \ldots,7.
$$

\vspace{4mm}

We have (see (\ref{30.52}))

\vspace{2mm}
$$
R_{T,t}^{(0)pppp}=
\int\limits_t^T
\int\limits_t^{t_4}\int\limits_t^{t_3}\int\limits_t^{t_2}
\sum\limits_{(t_1,t_2,t_3,t_4)}
\Biggl(
R_{pppp}(t_1,t_2,t_3,t_4)
d{\bf w}_{t_1}^{(i_1)}
d{\bf w}_{t_2}^{(i_2)}
d{\bf w}_{t_3}^{(i_3)}
d{\bf w}_{t_4}^{(i_4)}\Biggr),
$$

\vspace{5mm}
\noindent
where 
permutations
$(t_1,t_2,t_3,t_4)$ when summing
are performed only in the expression, which is enclosed in parentheses.

From the other hand \cite{16}, \cite{20}-\cite{20xxz}, \cite{26a}

\vspace{2mm}
$$
R_{T,t}^{(0)pppp}=
\sum\limits_{(t_1,t_2,t_3,t_4)}\int\limits_t^T
\int\limits_t^{t_4}\int\limits_t^{t_3}\int\limits_t^{t_2}
R_{pppp}(t_1,t_2,t_3,t_4)
d{\bf w}_{t_1}^{(i_1)}
d{\bf w}_{t_2}^{(i_2)}
d{\bf w}_{t_3}^{(i_3)}
d{\bf w}_{t_4}^{(i_4)},
$$

\vspace{4mm}
\noindent
where permutations
$(t_1,\ldots,t_4)$ when summing
are performed only in 
the values
$d{\bf w}_{t_1}^{(i_1)}
\ldots $
$d{\bf w}_{t_4}^{(i_4)}$. At the same time, the indexes near upper 
limits of integration in the iterated stochastic integrals are changed 
correspondently and if $t_r$ swapped with $t_q$ in the  
permutation $(t_1,\ldots,t_4)$, then $i_r$ swapped with $i_q$ in 
the permutations $(i_1,\ldots,i_4)$.

So, we obtain

\vspace{1mm}
$$
{\sf M}\left\{\left(R_{T,t}^{(0)pppp}\right)^2\right\}\le
24\sum\limits_{(t_1,t_2,t_3,t_4)}\int\limits_t^T
\int\limits_t^{t_4}\int\limits_t^{t_3}\int\limits_t^{t_2}
\left(R_{pppp}(t_1,t_2,t_3,t_4)\right)^2dt_1dt_2dt_3dt_4=
$$

\vspace{4mm}
$$
= 24\int\limits_{[t, T]^4}
\left(R_{pppp}(t_1,t_2,t_3,t_4))\right)^2
dt_1dt_2dt_3dt_4\ \to 0
$$

\vspace{5mm}
\noindent
if $p\to\infty,$ where $K^{*}(t_1,t_2,t_3,t_4)\in L_2([t, T]^4)$.

Let us consider $R_{T,t}^{(1)pppp}$

\vspace{3mm}
$$
R_{T,t}^{(1)pppp}={\bf 1}_{\{i_1=i_2\ne 0\}}
\hbox{\vtop{\offinterlineskip\halign{
\hfil#\hfil\cr
{\rm l.i.m.}\cr
$\stackrel{}{{}_{N\to \infty}}$\cr
}} }
\sum_{\stackrel{l_4,l_3,l_1=0}{{}_{l_1\ne l_3, l_1\ne l_4,
l_3\ne l_4}}}^{N-1}
R_{pppp}(\tau_{l_1},\tau_{l_1},\tau_{l_3},\tau_{l_4})
\Delta\tau_{l_1}
\Delta{\bf w}_{\tau_{l_3}}^{(i_3)}
\Delta{\bf w}_{\tau_{l_4}}^{(i_4)}=
$$

\vspace{2mm}
$$
={\bf 1}_{\{i_1=i_2\ne 0\}}
\hbox{\vtop{\offinterlineskip\halign{
\hfil#\hfil\cr
{\rm l.i.m.}\cr
$\stackrel{}{{}_{N\to \infty}}$\cr
}} }
\sum_{\stackrel{l_4,l_3,l_1=0}{{}_{l_3\ne l_4}}}^{N-1}
R_{pppp}(\tau_{l_1},\tau_{l_1},\tau_{l_3},\tau_{l_4})
\Delta\tau_{l_1}
\Delta{\bf w}_{\tau_{l_3}}^{(i_3)}
\Delta{\bf w}_{\tau_{l_4}}^{(i_4)}=
$$

\vspace{2mm}
$$
={\bf 1}_{\{i_1=i_2\ne 0\}}
\hbox{\vtop{\offinterlineskip\halign{
\hfil#\hfil\cr
{\rm l.i.m.}\cr
$\stackrel{}{{}_{N\to \infty}}$\cr
}} }
\sum_{\stackrel{l_4,l_3,l_1=0}{{}_{l_3\ne l_4}}}^{N-1}
\Biggl(
\frac{1}{2}{\bf 1}_{\{\tau_{l_1}<\tau_{l_3}<\tau_{l_4}\}}
+ 
\Biggr.
$$

\vspace{2mm}
$$
+\frac{1}{4}{\bf 1}_{\{\tau_{l_1}=\tau_{l_3}<\tau_{l_4}\}}+
\frac{1}{4}{\bf 1}_{\{\tau_{l_1}<\tau_{l_3}=\tau_{l_4}\}}+
\frac{1}{8}{\bf 1}_{\{\tau_{l_1}=\tau_{l_3}=\tau_{l_4}\}}-
$$

\vspace{2mm}
$$
\Biggl.
-\sum\limits_{j_4,j_3,j_2,j_1=0}^p C_{j_4j_3j_2j_1}
\phi_{j_1}(\tau_{l_1})
\phi_{j_2}(\tau_{l_1})
\phi_{j_3}(\tau_{l_3})
\phi_{j_4}(\tau_{l_4})\Biggr)
\Delta\tau_{l_1}
\Delta{\bf w}_{\tau_{l_3}}^{(i_3)}
\Delta{\bf w}_{\tau_{l_4}}^{(i_4)}=
$$

\vspace{2mm}
$$
={\bf 1}_{\{i_1=i_2\ne 0\}}
\hbox{\vtop{\offinterlineskip\halign{
\hfil#\hfil\cr
{\rm l.i.m.}\cr
$\stackrel{}{{}_{N\to \infty}}$\cr
}} }
\sum_{\stackrel{l_4,l_3,l_1=0}{{}_{l_3\ne l_4}}}^{N-1}
\Biggl(
\frac{1}{2}{\bf 1}_{\{\tau_{l_1}<\tau_{l_3}<\tau_{l_4}\}}
- 
\Biggr.
$$

\vspace{2mm}
$$
\Biggl.
-\sum\limits_{j_4,j_3,j_2,j_1=0}^p C_{j_4j_3j_2j_1}
\phi_{j_1}(\tau_{l_1})
\phi_{j_2}(\tau_{l_1})
\phi_{j_3}(\tau_{l_3})
\phi_{j_4}(\tau_{l_4})\Biggr)
\Delta\tau_{l_1}
\Delta{\bf w}_{\tau_{l_3}}^{(i_3)}
\Delta{\bf w}_{\tau_{l_4}}^{(i_4)}=
$$

\vspace{2mm}
$$
={\bf 1}_{\{i_1=i_2\ne 0\}}
\hbox{\vtop{\offinterlineskip\halign{
\hfil#\hfil\cr
{\rm l.i.m.}\cr
$\stackrel{}{{}_{N\to \infty}}$\cr
}} }
\sum\limits_{l_4=0}^{N-1}\sum\limits_{l_3=0}^{N-1}
\sum\limits_{l_1=0}^{N-1}
\Biggl(
\frac{1}{2}{\bf 1}_{\{\tau_{l_1}<\tau_{l_3}<\tau_{l_4}\}}
- 
\Biggr.
$$

\vspace{2mm}
$$
\Biggl.
-\sum\limits_{j_4,j_3,j_2,j_1=0}^p C_{j_4j_3j_2j_1}
\phi_{j_1}(\tau_{l_1})
\phi_{j_2}(\tau_{l_1})
\phi_{j_3}(\tau_{l_3})
\phi_{j_4}(\tau_{l_4})\Biggr)
\Delta\tau_{l_1}
\Delta{\bf w}_{\tau_{l_3}}^{(i_3)}
\Delta{\bf w}_{\tau_{l_4}}^{(i_4)}-
$$

\vspace{2mm}
$$
-{\bf 1}_{\{i_1=i_2\ne 0\}}{\bf 1}_{\{i_3=i_4\ne 0\}}
\hbox{\vtop{\offinterlineskip\halign{
\hfil#\hfil\cr
{\rm l.i.m.}\cr
$\stackrel{}{{}_{N\to \infty}}$\cr
}} }
\sum\limits_{l_4=0}^{N-1}
\sum\limits_{l_1=0}^{N-1}
\Biggl(
0
- 
\Biggr.
$$

\vspace{2mm}
$$
\Biggl.
-\sum\limits_{j_4,j_3,j_2,j_1=0}^p C_{j_4j_3j_2j_1}
\phi_{j_1}(\tau_{l_1})
\phi_{j_2}(\tau_{l_1})
\phi_{j_3}(\tau_{l_4})
\phi_{j_4}(\tau_{l_4})\Biggr)
\Delta\tau_{l_1}
\Delta\tau_{l_4}=
$$

\vspace{2mm}
$$
=
{\bf 1}_{\{i_1=i_2\ne 0\}}\left(\frac{1}{2}
\int\limits_t^T
\int\limits_t^{t_4}
\int\limits_t^{t_3}dt_1d{\bf w}_{t_3}^{(i_3)}
d{\bf w}_{t_4}^{(i_4)} - 
\sum\limits_{j_4,j_3,j_1=0}^p C_{j_4j_3j_1j_1}
\zeta_{j_3}^{(i_3)}
\zeta_{j_4}^{(i_4)}\right)+
$$

\vspace{2mm}
$$
+{\bf 1}_{\{i_1=i_2\ne 0\}}{\bf 1}_{\{i_3=i_4\ne 0\}}
\sum\limits_{j_4,j_1=0}^p C_{j_4j_4j_1j_1}\ \ \ \hbox{w. p. 1.}
$$

\vspace{6mm}

In \cite{16} (see the proof of Theorem 8, p.~A.135),
\cite{20} (see the proof of Theorem 5.7, p.~A.388),
\cite{20xx}-\cite{20xxz} (see Chapter 2),
\cite{32} (see the proof of Theorem 4) we have proved that

$$
\hbox{\vtop{\offinterlineskip\halign{
\hfil#\hfil\cr
{\rm lim}\cr
$\stackrel{}{{}_{p\to \infty}}$\cr
}} }
\sum\limits_{j_4,j_1=0}^{p} C_{j_4j_4j_1j_1}=
\frac{1}{4}
\int\limits_t^T
\int\limits_t^{t_2}dt_1dt_2,
$$

\vspace{2mm}
$$
\hbox{\vtop{\offinterlineskip\halign{
\hfil#\hfil\cr
{\rm l.i.m.}\cr
$\stackrel{}{{}_{p\to \infty}}$\cr
}} }
\sum\limits_{j_4,j_3,j_1=0}^p C_{j_4j_3j_1j_1}
\zeta_{j_3}^{(i_3)}
\zeta_{j_4}^{(i_4)}=
\frac{1}{2}
\int\limits_t^T
\int\limits_t^{t_4}
\int\limits_t^{t_3}dt_1d{\bf w}_{t_3}^{(i_3)}
d{\bf w}_{t_4}^{(i_4)}+
$$

\vspace{1mm}
$$
+{\bf 1}_{\{i_3=i_4\ne 0\}}
\frac{1}{4}
\int\limits_t^T
\int\limits_t^{t_2}dt_1dt_2\ \ \ \hbox{w. p. 1.}
$$

\vspace{3mm}

Then

\vspace{-2mm}
$$
\hbox{\vtop{\offinterlineskip\halign{
\hfil#\hfil\cr
{\rm lim}\cr
$\stackrel{}{{}_{p\to \infty}}$\cr
}} }
{\sf M}\left\{\left(R_{T,t}^{(1)pppp}\right)^2\right\}=0.
$$

\vspace{6mm}

Let us consider $R_{T,t}^{(2)pppp}$

\vspace{3mm}
$$
R_{T,t}^{(2)pppp}={\bf 1}_{\{i_1=i_3\ne 0\}}
\hbox{\vtop{\offinterlineskip\halign{
\hfil#\hfil\cr
{\rm l.i.m.}\cr
$\stackrel{}{{}_{N\to \infty}}$\cr
}} }
\sum_{\stackrel{l_4,l_2,l_1=0}{{}_{l_1\ne l_2, l_1\ne l_4,
l_2\ne l_4}}}^{N-1}
G_{pppp}(\tau_{l_1},\tau_{l_2},\tau_{l_1},\tau_{l_4})
\Delta\tau_{l_1}
\Delta{\bf w}_{\tau_{l_2}}^{(i_2)}
\Delta{\bf w}_{\tau_{l_4}}^{(i_4)}=
$$

\vspace{2mm}
$$
={\bf 1}_{\{i_1=i_3\ne 0\}}
\hbox{\vtop{\offinterlineskip\halign{
\hfil#\hfil\cr
{\rm l.i.m.}\cr
$\stackrel{}{{}_{N\to \infty}}$\cr
}} }
\sum_{\stackrel{l_4,l_2,l_1=0}{{}_{l_2\ne l_4}}}^{N-1}
G_{pppp}(\tau_{l_1},\tau_{l_2},\tau_{l_1},\tau_{l_4})
\Delta\tau_{l_1}
\Delta{\bf w}_{\tau_{l_2}}^{(i_2)}
\Delta{\bf w}_{\tau_{l_4}}^{(i_4)}=
$$

\vspace{2mm}
$$
={\bf 1}_{\{i_1=i_3\ne 0\}}
\hbox{\vtop{\offinterlineskip\halign{
\hfil#\hfil\cr
{\rm l.i.m.}\cr
$\stackrel{}{{}_{N\to \infty}}$\cr
}} }
\sum_{\stackrel{l_4,l_2,l_1=0}{{}_{l_2\ne l_4}}}^{N-1}
\Biggl(
\frac{1}{4}{\bf 1}_{\{\tau_{l_1}=\tau_{l_2}<\tau_{l_4}\}}
+ 
\Biggr.
\frac{1}{8}{\bf 1}_{\{\tau_{l_1}=\tau_{l_2}=\tau_{l_4}\}}-
$$

\vspace{2mm}
$$
\Biggl.
-\sum\limits_{j_4,j_3,j_2,j_1=0}^p C_{j_4j_3j_2j_1}
\phi_{j_1}(\tau_{l_1})
\phi_{j_2}(\tau_{l_2})
\phi_{j_3}(\tau_{l_1})
\phi_{j_4}(\tau_{l_4})\Biggr)
\Delta\tau_{l_1}
\Delta{\bf w}_{\tau_{l_2}}^{(i_2)}
\Delta{\bf w}_{\tau_{l_4}}^{(i_4)}=
$$

\vspace{2mm}
$$
={\bf 1}_{\{i_1=i_3\ne 0\}}
\hbox{\vtop{\offinterlineskip\halign{
\hfil#\hfil\cr
{\rm l.i.m.}\cr
$\stackrel{}{{}_{N\to \infty}}$\cr
}} }
\sum\limits_{l_4=0}^{N-1}\sum\limits_{l_2=0}^{N-1}
\sum\limits_{l_1=0}^{N-1}
(-1)\sum\limits_{j_4,j_3,j_2,j_1=0}^p C_{j_4j_3j_2j_1}\times
$$

\vspace{4mm}
$$
\times\phi_{j_1}(\tau_{l_1})
\phi_{j_2}(\tau_{l_2})
\phi_{j_3}(\tau_{l_1})
\phi_{j_4}(\tau_{l_4})
\Delta\tau_{l_1}
\Delta{\bf w}_{\tau_{l_2}}^{(i_2)}
\Delta{\bf w}_{\tau_{l_4}}^{(i_4)}-
$$

\vspace{2mm}
$$
-{\bf 1}_{\{i_1=i_3\ne 0\}}{\bf 1}_{\{i_2=i_4\ne 0\}}
\hbox{\vtop{\offinterlineskip\halign{
\hfil#\hfil\cr
{\rm l.i.m.}\cr
$\stackrel{}{{}_{N\to \infty}}$\cr
}} }
\sum\limits_{l_4=0}^{N-1}
\sum\limits_{l_1=0}^{N-1}
(-1)\sum\limits_{j_4,j_3,j_2,j_1=0}^p C_{j_4j_3j_2j_1}\times
$$

\vspace{4mm}
$$
\times
\phi_{j_1}(\tau_{l_1})
\phi_{j_2}(\tau_{l_4})
\phi_{j_3}(\tau_{l_1})
\phi_{j_4}(\tau_{l_4})
\Delta\tau_{l_1}
\Delta\tau_{l_4}=
$$

\vspace{2mm}
$$
=
-{\bf 1}_{\{i_1=i_3\ne 0\}}
\sum\limits_{j_4,j_2,j_1=0}^p C_{j_4j_1j_2j_1}
\zeta_{j_2}^{(i_2)}
\zeta_{j_4}^{(i_4)}+
$$

\vspace{2mm}
$$
+{\bf 1}_{\{i_1=i_3\ne 0\}}{\bf 1}_{\{i_2=i_4\ne 0\}}
\sum\limits_{j_4,j_1=0}^p C_{j_4j_1j_4j_1}\ \ \ \hbox{w. p. 1.}
$$

\vspace{6mm}

In \cite{16} (see the proof of Theorem 8, p.~A.135),
\cite{20} (see the proof of Theorem 5.7, p.~A.388),
\cite{20xx}-\cite{20xxz} (see Chapter 2),
\cite{32} (see the proof of Theorem 4) we have proved that

$$
\hbox{\vtop{\offinterlineskip\halign{
\hfil#\hfil\cr
{\rm l.i.m.}\cr
$\stackrel{}{{}_{p\to \infty}}$\cr
}} }
\sum\limits_{j_4,j_2,j_1=0}^p C_{j_4j_1j_2j_1}
\zeta_{j_2}^{(i_2)}
\zeta_{j_4}^{(i_4)}=0\ \ \ \hbox{w. p. 1,}
$$

\vspace{2mm}
$$
\hbox{\vtop{\offinterlineskip\halign{
\hfil#\hfil\cr
{\rm lim}\cr
$\stackrel{}{{}_{p\to \infty}}$\cr
}} }
\sum\limits_{j_4,j_1=0}^p C_{j_4j_1j_4j_1}=0.
$$

\vspace{4mm}

Then

\vspace{-2mm}
$$
\hbox{\vtop{\offinterlineskip\halign{
\hfil#\hfil\cr
{\rm lim}\cr
$\stackrel{}{{}_{p\to \infty}}$\cr
}} }
{\sf M}\left\{\left(R_{T,t}^{(2)pppp}\right)^2\right\}=0.
$$

\vspace{6mm}

Let us consider $R_{T,t}^{(3)pppp}$

\vspace{2mm}

$$
R_{T,t}^{(3)pppp}={\bf 1}_{\{i_1=i_4\ne 0\}}
\hbox{\vtop{\offinterlineskip\halign{
\hfil#\hfil\cr
{\rm l.i.m.}\cr
$\stackrel{}{{}_{N\to \infty}}$\cr
}} }
\sum_{\stackrel{l_3,l_2,l_1=0}{{}_{l_1\ne l_2, l_1\ne l_3,
l_2\ne l_3}}}^{N-1}
G_{pppp}(\tau_{l_1},\tau_{l_2},\tau_{l_3},\tau_{l_1})
\Delta\tau_{l_1}
\Delta{\bf w}_{\tau_{l_2}}^{(i_2)}
\Delta{\bf w}_{\tau_{l_3}}^{(i_3)}=
$$

\vspace{2mm}
$$
={\bf 1}_{\{i_1=i_4\ne 0\}}
\hbox{\vtop{\offinterlineskip\halign{
\hfil#\hfil\cr
{\rm l.i.m.}\cr
$\stackrel{}{{}_{N\to \infty}}$\cr
}} }
\sum_{\stackrel{l_3,l_2,l_1=0}{{}_{l_2\ne l_3}}}^{N-1}
G_{pppp}(\tau_{l_1},\tau_{l_2},\tau_{l_3},\tau_{l_1})
\Delta\tau_{l_1}
\Delta{\bf w}_{\tau_{l_2}}^{(i_2)}
\Delta{\bf w}_{\tau_{l_3}}^{(i_3)}=
$$

\vspace{2mm}
$$
={\bf 1}_{\{i_1=i_4\ne 0\}}
\hbox{\vtop{\offinterlineskip\halign{
\hfil#\hfil\cr
{\rm l.i.m.}\cr
$\stackrel{}{{}_{N\to \infty}}$\cr
}} }
\sum_{\stackrel{l_3,l_2,l_1=0}{{}_{l_2\ne l_3}}}^{N-1}
\Biggl(
\frac{1}{8}{\bf 1}_{\{\tau_{l_1}=\tau_{l_2}=\tau_{l_3}\}}-
$$

\vspace{2mm}
$$
\Biggl.
-\sum\limits_{j_4,j_3,j_2,j_1=0}^p C_{j_4j_3j_2j_1}
\phi_{j_1}(\tau_{l_1})
\phi_{j_2}(\tau_{l_2})
\phi_{j_3}(\tau_{l_3})
\phi_{j_4}(\tau_{l_1})\Biggr)
\Delta\tau_{l_1}
\Delta{\bf w}_{\tau_{l_2}}^{(i_2)}
\Delta{\bf w}_{\tau_{l_3}}^{(i_3)}=
$$

\vspace{2mm}
$$
={\bf 1}_{\{i_1=i_4\ne 0\}}
\hbox{\vtop{\offinterlineskip\halign{
\hfil#\hfil\cr
{\rm l.i.m.}\cr
$\stackrel{}{{}_{N\to \infty}}$\cr
}} }
\sum\limits_{l_3=0}^{N-1}\sum\limits_{l_2=0}^{N-1}
\sum\limits_{l_1=0}^{N-1}
(-1)\sum\limits_{j_4,j_3,j_2,j_1=0}^p C_{j_4j_3j_2j_1}\times
$$

\vspace{4mm}
$$
\times\phi_{j_1}(\tau_{l_1})
\phi_{j_2}(\tau_{l_2})
\phi_{j_3}(\tau_{l_3})
\phi_{j_4}(\tau_{l_1})
\Delta\tau_{l_1}
\Delta{\bf w}_{\tau_{l_2}}^{(i_2)}
\Delta{\bf w}_{\tau_{l_3}}^{(i_3)}-
$$

\vspace{2mm}
$$
-{\bf 1}_{\{i_1=i_4\ne 0\}}{\bf 1}_{\{i_2=i_3\ne 0\}}
\hbox{\vtop{\offinterlineskip\halign{
\hfil#\hfil\cr
{\rm l.i.m.}\cr
$\stackrel{}{{}_{N\to \infty}}$\cr
}} }
\sum\limits_{l_3=0}^{N-1}
\sum\limits_{l_1=0}^{N-1}
(-1)\sum\limits_{j_4,j_3,j_2,j_1=0}^p C_{j_4j_3j_2j_1}\times
$$

\vspace{4mm}
$$
\times
\phi_{j_1}(\tau_{l_1})
\phi_{j_2}(\tau_{l_3})
\phi_{j_3}(\tau_{l_3})
\phi_{j_4}(\tau_{l_1})
\Delta\tau_{l_1}
\Delta\tau_{l_3}=
$$

\vspace{2mm}
$$
=
-{\bf 1}_{\{i_1=i_4\ne 0\}}
\sum\limits_{j_4,j_3,j_2=0}^p C_{j_4j_3j_2j_4}
\zeta_{j_2}^{(i_2)}
\zeta_{j_3}^{(i_3)}+
$$

\vspace{2mm}
$$
+{\bf 1}_{\{i_1=i_4\ne 0\}}{\bf 1}_{\{i_2=i_3\ne 0\}}
\sum\limits_{j_4,j_2=0}^p C_{j_4j_2j_2j_4}\ \ \ \hbox{w. p. 1.}
$$

\vspace{6mm}

In \cite{16} (see the proof of Theorem 8, p.~A.135),
\cite{20} (see the proof of Theorem 5.7, p.~A.388),
\cite{20xx}-\cite{20xxz} (see Chapter 2),
\cite{32} (see the proof of Theorem 4) we have proved that

$$
\hbox{\vtop{\offinterlineskip\halign{
\hfil#\hfil\cr
{\rm l.i.m.}\cr
$\stackrel{}{{}_{p\to \infty}}$\cr
}} }
\sum\limits_{j_4,j_3,j_2=0}^p C_{j_4j_3j_2j_4}
\zeta_{j_2}^{(i_2)}
\zeta_{j_3}^{(i_3)}=0\ \ \ \hbox{w. p. 1,}
$$

\vspace{2mm}

$$
\hbox{\vtop{\offinterlineskip\halign{
\hfil#\hfil\cr
{\rm lim}\cr
$\stackrel{}{{}_{p\to \infty}}$\cr
}} }
\sum\limits_{j_4,j_2=0}^p C_{j_4j_2j_2j_4}=0.
$$

\vspace{4mm}

Then

\vspace{-2mm}

$$
\hbox{\vtop{\offinterlineskip\halign{
\hfil#\hfil\cr
{\rm lim}\cr
$\stackrel{}{{}_{p\to \infty}}$\cr
}} }
{\sf M}\left\{\left(R_{T,t}^{(3)pppp}\right)^2\right\}=0.
$$

\vspace{6mm}

Let us consider $R_{T,t}^{(4)pppp}$

\vspace{2mm}

$$
R_{T,t}^{(4)pppp}={\bf 1}_{\{i_2=i_3\ne 0\}}
\hbox{\vtop{\offinterlineskip\halign{
\hfil#\hfil\cr
{\rm l.i.m.}\cr
$\stackrel{}{{}_{N\to \infty}}$\cr
}} }
\sum_{\stackrel{l_4,l_2,l_1=0}{{}_{l_1\ne l_2, l_1\ne l_4,
l_2\ne l_4}}}^{N-1}
G_{pppp}(\tau_{l_1},\tau_{l_2},\tau_{l_2},\tau_{l_4})
\Delta{\bf w}_{\tau_{l_1}}^{(i_1)}
\Delta\tau_{l_2}
\Delta{\bf w}_{\tau_{l_4}}^{(i_4)}=
$$

\vspace{2mm}
$$
={\bf 1}_{\{i_2=i_3\ne 0\}}
\hbox{\vtop{\offinterlineskip\halign{
\hfil#\hfil\cr
{\rm l.i.m.}\cr
$\stackrel{}{{}_{N\to \infty}}$\cr
}} }
\sum_{\stackrel{l_4,l_2,l_1=0}{{}_{l_1\ne l_4}}}^{N-1}
G_{pppp}(\tau_{l_1},\tau_{l_2},\tau_{l_2},\tau_{l_4})
\Delta{\bf w}_{\tau_{l_1}}^{(i_1)}
\Delta\tau_{l_2}
\Delta{\bf w}_{\tau_{l_4}}^{(i_4)}=
$$

\vspace{2mm}
$$
={\bf 1}_{\{i_2=i_3\ne 0\}}
\hbox{\vtop{\offinterlineskip\halign{
\hfil#\hfil\cr
{\rm l.i.m.}\cr
$\stackrel{}{{}_{N\to \infty}}$\cr
}} }
\sum_{\stackrel{l_4,l_2,l_1=0}{{}_{l_1\ne l_4}}}^{N-1}
\Biggl(
\frac{1}{2}{\bf 1}_{\{\tau_{l_1}<\tau_{l_2}<\tau_{l_4}\}}
+ 
\Biggr.
$$

\vspace{2mm}
$$
+\frac{1}{4}{\bf 1}_{\{\tau_{l_1}=\tau_{l_2}<\tau_{l_4}\}}+
\frac{1}{4}{\bf 1}_{\{\tau_{l_1}<\tau_{l_2}=\tau_{l_4}\}}+
\frac{1}{8}{\bf 1}_{\{\tau_{l_1}=\tau_{l_2}=\tau_{l_4}\}}-
$$

\vspace{2mm}
$$
\Biggl.
-\sum\limits_{j_4,j_3,j_2,j_1=0}^p C_{j_4j_3j_2j_1}
\phi_{j_1}(\tau_{l_1})
\phi_{j_2}(\tau_{l_2})
\phi_{j_3}(\tau_{l_2})
\phi_{j_4}(\tau_{l_4})\Biggr)
\Delta{\bf w}_{\tau_{l_1}}^{(i_1)}
\Delta\tau_{l_2}
\Delta{\bf w}_{\tau_{l_4}}^{(i_4)}=
$$

\vspace{2mm}
$$
={\bf 1}_{\{i_2=i_3\ne 0\}}
\hbox{\vtop{\offinterlineskip\halign{
\hfil#\hfil\cr
{\rm l.i.m.}\cr
$\stackrel{}{{}_{N\to \infty}}$\cr
}} }
\sum_{\stackrel{l_4,l_2,l_1=0}{{}_{l_1\ne l_4}}}^{N-1}
\Biggl(
\frac{1}{2}{\bf 1}_{\{\tau_{l_1}<\tau_{l_2}<\tau_{l_4}\}}
- 
\Biggr.
$$

\vspace{2mm}
$$
\Biggl.
-\sum\limits_{j_4,j_3,j_2,j_1=0}^p C_{j_4j_3j_2j_1}
\phi_{j_1}(\tau_{l_1})
\phi_{j_2}(\tau_{l_2})
\phi_{j_3}(\tau_{l_2})
\phi_{j_4}(\tau_{l_4})\Biggr)
\Delta{\bf w}_{\tau_{l_1}}^{(i_1)}
\Delta\tau_{l_2}
\Delta{\bf w}_{\tau_{l_4}}^{(i_4)}=
$$

\vspace{2mm}
$$
={\bf 1}_{\{i_2=i_3\ne 0\}}
\hbox{\vtop{\offinterlineskip\halign{
\hfil#\hfil\cr
{\rm l.i.m.}\cr
$\stackrel{}{{}_{N\to \infty}}$\cr
}} }
\sum\limits_{l_4=0}^{N-1}\sum\limits_{l_2=0}^{N-1}
\sum\limits_{l_1=0}^{N-1}
\Biggl(
\frac{1}{2}{\bf 1}_{\{\tau_{l_1}<\tau_{l_2}<\tau_{l_4}\}}
- 
\Biggr.
$$

\vspace{2mm}
$$
\Biggl.
-\sum\limits_{j_4,j_3,j_2,j_1=0}^p C_{j_4j_3j_2j_1}
\phi_{j_1}(\tau_{l_1})
\phi_{j_2}(\tau_{l_2})
\phi_{j_3}(\tau_{l_2})
\phi_{j_4}(\tau_{l_4})\Biggr)
\Delta{\bf w}_{\tau_{l_1}}^{(i_1)}
\Delta\tau_{l_2}
\Delta{\bf w}_{\tau_{l_4}}^{(i_4)}-
$$

\vspace{2mm}
$$
-{\bf 1}_{\{i_2=i_3\ne 0\}}{\bf 1}_{\{i_1=i_4\ne 0\}}
\hbox{\vtop{\offinterlineskip\halign{
\hfil#\hfil\cr
{\rm l.i.m.}\cr
$\stackrel{}{{}_{N\to \infty}}$\cr
}} }
\sum\limits_{l_4=0}^{N-1}
\sum\limits_{l_2=0}^{N-1}
(-1)
\sum\limits_{j_4,j_3,j_2,j_1=0}^p C_{j_4j_3j_2j_1}\times
$$

\vspace{4mm}
$$
\times
\phi_{j_1}(\tau_{l_4})
\phi_{j_2}(\tau_{l_2})
\phi_{j_3}(\tau_{l_2})
\phi_{j_4}(\tau_{l_4})
\Delta\tau_{l_2}
\Delta\tau_{l_4}=
$$

\vspace{2mm}
$$
=
{\bf 1}_{\{i_2=i_3\ne 0\}}\left(\frac{1}{2}
\int\limits_t^T
\int\limits_t^{t_4}
\int\limits_t^{t_2}d{\bf w}_{t_1}^{(i_1)}dt_2
d{\bf w}_{t_4}^{(i_4)} - 
\sum\limits_{j_4,j_2,j_1=0}^p C_{j_4j_2j_2j_1}
\zeta_{j_1}^{(i_1)}
\zeta_{j_4}^{(i_4)}\right)+
$$

\vspace{2mm}
$$
+{\bf 1}_{\{i_2=i_3\ne 0\}}{\bf 1}_{\{i_1=i_4\ne 0\}}
\sum\limits_{j_4,j_2=0}^p C_{j_4j_2j_2j_4}\ \ \ \hbox{w. p. 1.}
$$

\vspace{6mm}

In \cite{16} (see the proof of Theorem 8, p.~A.135),
\cite{20} (see the proof of Theorem 5.7, p.~A.388),
\cite{20xx}-\cite{20xxz} (see Chapter 2),
\cite{32} (see the proof of Theorem 4) we have proved that

$$
\hbox{\vtop{\offinterlineskip\halign{
\hfil#\hfil\cr
{\rm lim}\cr
$\stackrel{}{{}_{p\to \infty}}$\cr
}} }
\sum\limits_{j_4,j_2=0}^{p} C_{j_4j_2j_2j_4}=0,
$$

\vspace{2mm}
$$
\hbox{\vtop{\offinterlineskip\halign{
\hfil#\hfil\cr
{\rm l.i.m.}\cr
$\stackrel{}{{}_{p\to \infty}}$\cr
}} }
\sum\limits_{j_4,j_2,j_1=0}^p C_{j_4j_2j_2j_1}
\zeta_{j_1}^{(i_1)}
\zeta_{j_4}^{(i_4)}=
\frac{1}{2}
\int\limits_t^T
\int\limits_t^{t_4}
\int\limits_t^{t_2}d{\bf w}_{t_1}^{(i_1)}dt_2
d{\bf w}_{t_4}^{(i_4)}\ \ \ \hbox{w. p. 1.}
$$

\vspace{7mm}

Then

\vspace{-4mm}

$$
\hbox{\vtop{\offinterlineskip\halign{
\hfil#\hfil\cr
{\rm lim}\cr
$\stackrel{}{{}_{p\to \infty}}$\cr
}} }
{\sf M}\left\{\left(R_{T,t}^{(4)pppp}\right)^2\right\}=0.
$$

\vspace{6mm}

Let us consider $R_{T,t}^{(5)pppp}$

\vspace{2mm}

$$
R_{T,t}^{(5)pppp}={\bf 1}_{\{i_2=i_4\ne 0\}}
\hbox{\vtop{\offinterlineskip\halign{
\hfil#\hfil\cr
{\rm l.i.m.}\cr
$\stackrel{}{{}_{N\to \infty}}$\cr
}} }
\sum_{\stackrel{l_3,l_2,l_1=0}{{}_{l_1\ne l_2, l_1\ne l_3,
l_2\ne l_3}}}^{N-1}
G_{pppp}(\tau_{l_1},\tau_{l_2},\tau_{l_3},\tau_{l_2})
\Delta{\bf w}_{\tau_{l_1}}^{(i_1)}
\Delta\tau_{l_2}
\Delta{\bf w}_{\tau_{l_3}}^{(i_3)}=
$$

\vspace{2mm}
$$
={\bf 1}_{\{i_2=i_4\ne 0\}}
\hbox{\vtop{\offinterlineskip\halign{
\hfil#\hfil\cr
{\rm l.i.m.}\cr
$\stackrel{}{{}_{N\to \infty}}$\cr
}} }
\sum_{\stackrel{l_3,l_2,l_1=0}{{}_{l_1\ne l_3}}}^{N-1}
G_{pppp}(\tau_{l_1},\tau_{l_2},\tau_{l_3},\tau_{l_2})
\Delta{\bf w}_{\tau_{l_1}}^{(i_1)}
\Delta\tau_{l_2}
\Delta{\bf w}_{\tau_{l_3}}^{(i_3)}=
$$

\vspace{2mm}
$$
={\bf 1}_{\{i_2=i_4\ne 0\}}
\hbox{\vtop{\offinterlineskip\halign{
\hfil#\hfil\cr
{\rm l.i.m.}\cr
$\stackrel{}{{}_{N\to \infty}}$\cr
}} }
\sum_{\stackrel{l_3,l_2,l_1=0}{{}_{l_1\ne l_3}}}^{N-1}
\Biggl(
\frac{1}{4}{\bf 1}_{\{\tau_{l_1}<\tau_{l_2}=\tau_{l_3}\}}
+ 
\Biggr.
\frac{1}{8}{\bf 1}_{\{\tau_{l_1}=\tau_{l_2}=\tau_{l_3}\}}-
$$

\vspace{2mm}
$$
\Biggl.
-\sum\limits_{j_4,j_3,j_2,j_1=0}^p C_{j_4j_3j_2j_1}
\phi_{j_1}(\tau_{l_1})
\phi_{j_2}(\tau_{l_2})
\phi_{j_3}(\tau_{l_3})
\phi_{j_4}(\tau_{l_2})\Biggr)
\Delta{\bf w}_{\tau_{l_1}}^{(i_1)}
\Delta\tau_{l_2}
\Delta{\bf w}_{\tau_{l_3}}^{(i_3)}=
$$

\vspace{2mm}
$$
={\bf 1}_{\{i_2=i_4\ne 0\}}
\hbox{\vtop{\offinterlineskip\halign{
\hfil#\hfil\cr
{\rm l.i.m.}\cr
$\stackrel{}{{}_{N\to \infty}}$\cr
}} }
\sum_{\stackrel{l_3,l_2,l_1=0}{{}_{l_1\ne l_3}}}^{N-1}
(-1)\sum\limits_{j_4,j_3,j_2,j_1=0}^p C_{j_4j_3j_2j_1}\times
$$

\vspace{4mm}
$$
\times
\phi_{j_1}(\tau_{l_1})
\phi_{j_2}(\tau_{l_2})
\phi_{j_3}(\tau_{l_3})
\phi_{j_4}(\tau_{l_2})
\Delta{\bf w}_{\tau_{l_1}}^{(i_1)}
\Delta\tau_{l_2}
\Delta{\bf w}_{\tau_{l_3}}^{(i_3)}=
$$

\vspace{2mm}
$$
=
-{\bf 1}_{\{i_2=i_4\ne 0\}}
\sum\limits_{j_4,j_3,j_1=0}^p C_{j_4j_3j_4j_1}
\zeta_{j_1}^{(i_1)}
\zeta_{j_3}^{(i_3)}-
$$

\vspace{2mm}
$$
-{\bf 1}_{\{i_2=i_4\ne 0\}}{\bf 1}_{\{i_1=i_3\ne 0\}}
\hbox{\vtop{\offinterlineskip\halign{
\hfil#\hfil\cr
{\rm l.i.m.}\cr
$\stackrel{}{{}_{N\to \infty}}$\cr
}} }
\sum\limits_{l_3=0}^{N-1}
\sum\limits_{l_2=0}^{N-1}
(-1)\sum\limits_{j_4,j_3,j_2,j_1=0}^p C_{j_4j_3j_2j_1}\times
$$

\vspace{4mm}
$$
\times
\phi_{j_1}(\tau_{l_3})
\phi_{j_2}(\tau_{l_2})
\phi_{j_3}(\tau_{l_3})
\phi_{j_4}(\tau_{l_2})
\Delta\tau_{l_2}
\Delta\tau_{l_3}=
$$

\vspace{2mm}
$$
=
-{\bf 1}_{\{i_2=i_4\ne 0\}}
\sum\limits_{j_4,j_3,j_1=0}^p C_{j_4j_3j_4j_1}
\zeta_{j_1}^{(i_1)}
\zeta_{j_3}^{(i_3)}+
$$

\vspace{2mm}
$$
+{\bf 1}_{\{i_2=i_4\ne 0\}}{\bf 1}_{\{i_1=i_3\ne 0\}}
\sum\limits_{j_4,j_1=0}^p C_{j_4j_1j_4j_1}\ \ \ \hbox{w. p. 1.}
$$

\vspace{6mm}

In \cite{16} (see the proof of Theorem 8, p.~A.135),
\cite{20} (see the proof of Theorem 5.7, p.~A.388),
\cite{20xx}-\cite{20xxz} (see Chapter 2),
\cite{32} (see the proof of Theorem 4) we have proved that

$$
\hbox{\vtop{\offinterlineskip\halign{
\hfil#\hfil\cr
{\rm l.i.m.}\cr
$\stackrel{}{{}_{p\to \infty}}$\cr
}} }
\sum\limits_{j_4,j_3,j_1=0}^p C_{j_4j_3j_4j_1}
\zeta_{j_1}^{(i_1)}
\zeta_{j_3}^{(i_3)}=0\ \ \ \hbox{w. p. 1,}
$$

\vspace{2mm}
$$
\hbox{\vtop{\offinterlineskip\halign{
\hfil#\hfil\cr
{\rm lim}\cr
$\stackrel{}{{}_{p\to \infty}}$\cr
}} }
\sum\limits_{j_4,j_1=0}^p C_{j_4j_1j_4j_1}=0.
$$

\vspace{4mm}

Then

\vspace{-2mm}

$$
\hbox{\vtop{\offinterlineskip\halign{
\hfil#\hfil\cr
{\rm lim}\cr
$\stackrel{}{{}_{p\to \infty}}$\cr
}} }
{\sf M}\left\{\left(R_{T,t}^{(5)pppp}\right)^2\right\}=0.
$$

\vspace{6mm}

Let us consider $R_{T,t}^{(6)pppp}$

\vspace{2mm}

$$
R_{T,t}^{(6)pppp}={\bf 1}_{\{i_3=i_4\ne 0\}}
\hbox{\vtop{\offinterlineskip\halign{
\hfil#\hfil\cr
{\rm l.i.m.}\cr
$\stackrel{}{{}_{N\to \infty}}$\cr
}} }
\sum_{\stackrel{l_3,l_2,l_1=0}{{}_{l_1\ne l_2, l_1\ne l_3,
l_2\ne l_3}}}^{N-1}
G_{pppp}(\tau_{l_1},\tau_{l_2},\tau_{l_3},\tau_{l_3})
\Delta{\bf w}_{\tau_{l_1}}^{(i_1)}
\Delta{\bf w}_{\tau_{l_2}}^{(i_2)}
\Delta\tau_{l_3}=
$$

\vspace{2mm}
$$
={\bf 1}_{\{i_3=i_4\ne 0\}}
\hbox{\vtop{\offinterlineskip\halign{
\hfil#\hfil\cr
{\rm l.i.m.}\cr
$\stackrel{}{{}_{N\to \infty}}$\cr
}} }
\sum_{\stackrel{l_3,l_2,l_1=0}{{}_{l_1\ne l_2}}}^{N-1}
G_{pppp}(\tau_{l_1},\tau_{l_2},\tau_{l_3},\tau_{l_3})
\Delta{\bf w}_{\tau_{l_1}}^{(i_1)}
\Delta{\bf w}_{\tau_{l_2}}^{(i_2)}
\Delta\tau_{l_3}=
$$

\vspace{2mm}
$$
={\bf 1}_{\{i_3=i_4\ne 0\}}
\hbox{\vtop{\offinterlineskip\halign{
\hfil#\hfil\cr
{\rm l.i.m.}\cr
$\stackrel{}{{}_{N\to \infty}}$\cr
}} }
\sum_{\stackrel{l_3,l_2,l_1=0}{{}_{l_1\ne l_2}}}^{N-1}
\Biggl(
\frac{1}{2}{\bf 1}_{\{\tau_{l_1}<\tau_{l_2}<\tau_{l_3}\}}
+ 
\Biggr.
$$

\vspace{2mm}
$$
+\frac{1}{4}{\bf 1}_{\{\tau_{l_1}=\tau_{l_2}<\tau_{l_3}\}}+
\frac{1}{4}{\bf 1}_{\{\tau_{l_1}<\tau_{l_2}=\tau_{l_3}\}}+
\frac{1}{8}{\bf 1}_{\{\tau_{l_1}=\tau_{l_2}=\tau_{l_3}\}}-
$$

\vspace{2mm}
$$
\Biggl.
-\sum\limits_{j_4,j_3,j_2,j_1=0}^p C_{j_4j_3j_2j_1}
\phi_{j_1}(\tau_{l_1})
\phi_{j_2}(\tau_{l_2})
\phi_{j_3}(\tau_{l_3})
\phi_{j_4}(\tau_{l_3})\Biggr)
\Delta{\bf w}_{\tau_{l_1}}^{(i_1)}
\Delta{\bf w}_{\tau_{l_2}}^{(i_2)}
\Delta\tau_{l_3}=
$$

\vspace{2mm}
$$
={\bf 1}_{\{i_3=i_4\ne 0\}}
\hbox{\vtop{\offinterlineskip\halign{
\hfil#\hfil\cr
{\rm l.i.m.}\cr
$\stackrel{}{{}_{N\to \infty}}$\cr
}} }
\sum_{\stackrel{l_3,l_2,l_1=0}{{}_{l_1\ne l_2}}}^{N-1}
\Biggl(
\frac{1}{2}{\bf 1}_{\{\tau_{l_1}<\tau_{l_2}<\tau_{l_3}\}}
- 
\Biggr.
$$

\vspace{2mm}
$$
\Biggl.
-\sum\limits_{j_4,j_3,j_2,j_1=0}^p C_{j_4j_3j_2j_1}
\phi_{j_1}(\tau_{l_1})
\phi_{j_2}(\tau_{l_2})
\phi_{j_3}(\tau_{l_3})
\phi_{j_4}(\tau_{l_3})\Biggr)
\Delta{\bf w}_{\tau_{l_1}}^{(i_1)}
\Delta{\bf w}_{\tau_{l_2}}^{(i_2)}
\Delta\tau_{l_3}=
$$

\vspace{2mm}
$$
={\bf 1}_{\{i_3=i_4\ne 0\}}
\left(\frac{1}{2}
\int\limits_t^T
\int\limits_t^{t_3}
\int\limits_t^{t_2}d{\bf w}_{t_1}^{(i_1)}
d{\bf w}_{t_2}^{(i_2)}dt_3 - 
\sum\limits_{j_4,j_2,j_1=0}^p C_{j_4j_4j_2j_1}
\zeta_{j_1}^{(i_1)}
\zeta_{j_2}^{(i_2)}\right)-
$$

\vspace{2mm}
$$
-{\bf 1}_{\{i_3=i_4\ne 0\}}{\bf 1}_{\{i_1=i_2\ne 0\}}
\hbox{\vtop{\offinterlineskip\halign{
\hfil#\hfil\cr
{\rm l.i.m.}\cr
$\stackrel{}{{}_{N\to \infty}}$\cr
}} }
\sum\limits_{l_3=0}^{N-1}
\sum\limits_{l_1=0}^{N-1}
(-1)
\sum\limits_{j_4,j_3,j_2,j_1=0}^p C_{j_4j_3j_2j_1}\times
$$

\vspace{4mm}
$$
\times
\phi_{j_1}(\tau_{l_1})
\phi_{j_2}(\tau_{l_1})
\phi_{j_3}(\tau_{l_3})
\phi_{j_4}(\tau_{l_3})
\Delta\tau_{l_1}
\Delta\tau_{l_3}=
$$

\vspace{2mm}
$$
={\bf 1}_{\{i_3=i_4\ne 0\}}
\left(\frac{1}{2}
\int\limits_t^T
\int\limits_t^{t_3}
\int\limits_t^{t_2}d{\bf w}_{t_1}^{(i_1)}
d{\bf w}_{t_2}^{(i_2)}dt_3 - 
\sum\limits_{j_4,j_2,j_1=0}^p C_{j_4j_4j_2j_1}
\zeta_{j_1}^{(i_1)}
\zeta_{j_2}^{(i_2)}\right)+
$$

\vspace{2mm}
$$
+{\bf 1}_{\{i_1=i_2\ne 0\}}{\bf 1}_{\{i_3=i_4\ne 0\}}
\sum\limits_{j_4,j_1=0}^p C_{j_4j_4j_1j_1}=
$$

\vspace{2mm}
$$
={\bf 1}_{\{i_3=i_4\ne 0\}}
\left(\frac{1}{2}
\int\limits_t^T
\int\limits_t^{t_3}
\int\limits_t^{t_2}d{\bf w}_{t_1}^{(i_1)}
d{\bf w}_{t_2}^{(i_2)}dt_3 
+\frac{1}{4}{\bf 1}_{\{i_1=i_2\ne 0\}}
\int\limits_t^T\int\limits_t^{t_3}dt_1dt_3-
\right.
$$

\vspace{2mm}
$$
\left.
-\sum\limits_{j_4,j_2,j_1=0}^p C_{j_4j_4j_2j_1}
\zeta_{j_1}^{(i_1)}
\zeta_{j_2}^{(i_2)}\right)+
$$

\vspace{2mm}
$$
+{\bf 1}_{\{i_1=i_2\ne 0\}}{\bf 1}_{\{i_3=i_4\ne 0\}}
\left(\sum\limits_{j_4,j_1=0}^p C_{j_4j_4j_1j_1}-
\frac{1}{4}
\int\limits_t^T\int\limits_t^{t_3}dt_1dt_3\right)\ \ \ \hbox{w. p. 1.}
$$

\vspace{6mm}

In \cite{16} (see the proof of Theorem 8, p.~A.135),
\cite{20} (see the proof of Theorem 5.7, p.~A.388),
\cite{20xx}-\cite{20xxz} (see Chapter 2),
\cite{32} (see the proof of Theorem 4) we have proved that

$$
\hbox{\vtop{\offinterlineskip\halign{
\hfil#\hfil\cr
{\rm lim}\cr
$\stackrel{}{{}_{p\to \infty}}$\cr
}} }
\sum\limits_{j_4,j_1=0}^{p} C_{j_4j_4j_1j_1}=
\frac{1}{4}
\int\limits_t^T\int\limits_t^{t_3}dt_1dt_3,
$$

\vspace{2mm}
$$
\hbox{\vtop{\offinterlineskip\halign{
\hfil#\hfil\cr
{\rm l.i.m.}\cr
$\stackrel{}{{}_{p\to \infty}}$\cr
}} }
\sum\limits_{j_4,j_2,j_1=0}^p C_{j_4j_4j_2j_1}
\zeta_{j_1}^{(i_1)}
\zeta_{j_2}^{(i_2)}=
\frac{1}{2}
\int\limits_t^T
\int\limits_t^{t_3}
\int\limits_t^{t_2}{\bf w}_{t_1}^{(i_1)}
d{\bf w}_{t_2}^{(i_2)}dt_3+
$$

$$
+
{\bf 1}_{\{i_1=i_2\ne 0\}}\frac{1}{4}
\int\limits_t^T\int\limits_t^{t_3}dt_1dt_3\ \ \ \hbox{w. p. 1.}
$$

\vspace{4mm}

Then

\vspace{-2mm}

$$
\hbox{\vtop{\offinterlineskip\halign{
\hfil#\hfil\cr
{\rm lim}\cr
$\stackrel{}{{}_{p\to \infty}}$\cr
}} }
{\sf M}\left\{\left(R_{T,t}^{(6)pppp}\right)^2\right\}=0.
$$

\vspace{6mm}

Let us consider $R_{T,t}^{(7)pppp}$

\vspace{2mm}

$$
R_{T,t}^{(7)pppp}={\bf 1}_{\{i_1=i_2\ne 0\}}{\bf 1}_{\{i_3=i_4\ne 0\}}
\hbox{\vtop{\offinterlineskip\halign{
\hfil#\hfil\cr
{\rm l.i.m.}\cr
$\stackrel{}{{}_{N\to \infty}}$\cr
}} }
\sum_{\stackrel{l_4,l_2=0}{{}_{l_2\ne l_4}}}^{N-1}
G_{pppp}(\tau_{l_2},\tau_{l_2},\tau_{l_4},\tau_{l_4})
\Delta\tau_{l_2}
\Delta\tau_{l_4}
+
$$

\vspace{2mm}
$$
+{\bf 1}_{\{i_1=i_3\ne 0\}}{\bf 1}_{\{i_2=i_4\ne 0\}}
\hbox{\vtop{\offinterlineskip\halign{
\hfil#\hfil\cr
{\rm l.i.m.}\cr
$\stackrel{}{{}_{N\to \infty}}$\cr
}} }
\sum_{\stackrel{l_4,l_2=0}{{}_{l_2\ne l_4}}}^{N-1}
G_{pppp}(\tau_{l_2},\tau_{l_4},\tau_{l_2},\tau_{l_4})
\Delta\tau_{l_2}
\Delta\tau_{l_4}+
$$

\vspace{2mm}
$$
+{\bf 1}_{\{i_1=i_4\ne 0\}}{\bf 1}_{\{i_2=i_3\ne 0\}}
\hbox{\vtop{\offinterlineskip\halign{
\hfil#\hfil\cr
{\rm l.i.m.}\cr
$\stackrel{}{{}_{N\to \infty}}$\cr
}} }
\sum_{\stackrel{l_4,l_2=0}{{}_{l_2\ne l_4}}}^{N-1}
G_{pppp}(\tau_{l_2},\tau_{l_4},\tau_{l_4},\tau_{l_2})
\Delta\tau_{l_2}
\Delta\tau_{l_4}=
$$

\vspace{2mm}
$$
={\bf 1}_{\{i_1=i_2\ne 0\}}{\bf 1}_{\{i_3=i_4\ne 0\}}
\hbox{\vtop{\offinterlineskip\halign{
\hfil#\hfil\cr
{\rm l.i.m.}\cr
$\stackrel{}{{}_{N\to \infty}}$\cr
}} }
\sum\limits_{l_4=0}^{N-1}
\sum\limits_{l_2=0}^{N-1} 
G_{pppp}(\tau_{l_2},\tau_{l_2},\tau_{l_4},\tau_{l_4})
\Delta\tau_{l_2}
\Delta\tau_{l_4}+
$$

\vspace{2mm}
$$
+{\bf 1}_{\{i_1=i_3\ne 0\}}{\bf 1}_{\{i_2=i_4\ne 0\}}
\hbox{\vtop{\offinterlineskip\halign{
\hfil#\hfil\cr
{\rm l.i.m.}\cr
$\stackrel{}{{}_{N\to \infty}}$\cr
}} }
\sum\limits_{l_4=0}^{N-1}
\sum\limits_{l_2=0}^{N-1} 
G_{pppp}(\tau_{l_2},\tau_{l_4},\tau_{l_2},\tau_{l_4})
\Delta\tau_{l_2}
\Delta\tau_{l_4}+
$$

\vspace{2mm}
$$
+{\bf 1}_{\{i_1=i_4\ne 0\}}{\bf 1}_{\{i_2=i_3\ne 0\}}
\hbox{\vtop{\offinterlineskip\halign{
\hfil#\hfil\cr
{\rm l.i.m.}\cr
$\stackrel{}{{}_{N\to \infty}}$\cr
}} }
\sum\limits_{l_4=0}^{N-1}
\sum\limits_{l_2=0}^{N-1} 
G_{pppp}(\tau_{l_2},\tau_{l_4},\tau_{l_4},\tau_{l_2})
\Delta\tau_{l_2}
\Delta\tau_{l_4}=
$$

\vspace{2mm}
$$
={\bf 1}_{\{i_1=i_2\ne 0\}}{\bf 1}_{\{i_3=i_4\ne 0\}}
\hbox{\vtop{\offinterlineskip\halign{
\hfil#\hfil\cr
{\rm l.i.m.}\cr
$\stackrel{}{{}_{N\to \infty}}$\cr
}} }
\sum\limits_{l_4=0}^{N-1}
\sum\limits_{l_2=0}^{N-1} 
\Biggl(\frac{1}{4}{\bf 1}_{\{\tau_{l_2}<\tau_{l_4}\}}
+\frac{1}{8}{\bf 1}_{\{\tau_{l_2}=\tau_{l_4}\}}-\Biggr.
$$

\vspace{2mm}
$$
\Biggl.
-\sum\limits_{j_4,j_3,j_2,j_1=0}^p C_{j_4j_3j_2j_1}
\phi_{j_1}(\tau_{l_2})
\phi_{j_2}(\tau_{l_2})
\phi_{j_3}(\tau_{l_4})
\phi_{j_4}(\tau_{l_4})\Biggr)
\Delta\tau_{l_2}
\Delta\tau_{l_4}+
$$

\vspace{2mm}
$$
+{\bf 1}_{\{i_1=i_3\ne 0\}}{\bf 1}_{\{i_2=i_4\ne 0\}}
\hbox{\vtop{\offinterlineskip\halign{
\hfil#\hfil\cr
{\rm l.i.m.}\cr
$\stackrel{}{{}_{N\to \infty}}$\cr
}} }
\sum\limits_{l_4=0}^{N-1}
\sum\limits_{l_2=0}^{N-1}\Biggl(
\frac{1}{8}{\bf 1}_{\{\tau_{l_2}=\tau_{l_4}\}}-
\sum\limits_{j_4,j_3,j_2,j_1=0}^p C_{j_4j_3j_2j_1}\times\Biggr.
$$

\vspace{2mm}
$$
\Biggl.\times
\phi_{j_1}(\tau_{l_2})
\phi_{j_2}(\tau_{l_4})
\phi_{j_3}(\tau_{l_2})
\phi_{j_4}(\tau_{l_4})\Biggr)
\Delta\tau_{l_2}
\Delta\tau_{l_4}+
$$

\vspace{2mm}
$$
+{\bf 1}_{\{i_1=i_4\ne 0\}}{\bf 1}_{\{i_2=i_3\ne 0\}}
\hbox{\vtop{\offinterlineskip\halign{
\hfil#\hfil\cr
{\rm l.i.m.}\cr
$\stackrel{}{{}_{N\to \infty}}$\cr
}} }
\sum\limits_{l_4=0}^{N-1}
\sum\limits_{l_2=0}^{N-1} 
\Biggl(\frac{1}{8}{\bf 1}_{\{\tau_{l_2}=\tau_{l_4}\}}-
\sum\limits_{j_4,j_3,j_2,j_1=0}^p C_{j_4j_3j_2j_1}\times\Biggr.
$$

\vspace{2mm}
$$
\Biggl.\times
\phi_{j_1}(\tau_{l_2})
\phi_{j_2}(\tau_{l_4})
\phi_{j_3}(\tau_{l_4})
\phi_{j_4}(\tau_{l_2})\Biggr)
\Delta\tau_{l_2}
\Delta\tau_{l_4}=
$$

\vspace{2mm}
$$
=
{\bf 1}_{\{i_1=i_2\ne 0\}}{\bf 1}_{\{i_3=i_4\ne 0\}}
\left(
\frac{1}{4}
\int\limits_t^T\int\limits_t^{t_4}dt_2dt_4
-\sum\limits_{j_4,j_1=0}^p C_{j_4j_4j_1j_1}\right)-
$$

\vspace{2mm}
$$
-{\bf 1}_{\{i_1=i_3\ne 0\}}{\bf 1}_{\{i_2=i_4\ne 0\}}
\sum\limits_{j_4,j_1=0}^p C_{j_4j_1j_4j_1}-
$$

\vspace{2mm}
$$
-{\bf 1}_{\{i_1=i_4\ne 0\}}{\bf 1}_{\{i_2=i_3\ne 0\}}
\sum\limits_{j_4,j_2=0}^p C_{j_4j_2j_2j_4}.
$$

\vspace{6mm}

In \cite{16} (see the proof of Theorem 8, p.~A.135),
\cite{20} (see the proof of Theorem 5.7, p.~A.388),
\cite{20xx}-\cite{20xxz} (see Chapter 2),
\cite{32} (see the proof of Theorem 4) we have proved that

$$
\hbox{\vtop{\offinterlineskip\halign{
\hfil#\hfil\cr
{\rm lim}\cr
$\stackrel{}{{}_{p\to \infty}}$\cr
}} }
\sum\limits_{j_4,j_1=0}^{p} C_{j_4j_4j_1j_1}=
\frac{1}{4}
\int\limits_t^T\int\limits_t^{t_4}dt_2dt_4,
$$

\vspace{2mm}
$$
\hbox{\vtop{\offinterlineskip\halign{
\hfil#\hfil\cr
{\rm lim}\cr
$\stackrel{}{{}_{p\to \infty}}$\cr
}} }
\sum\limits_{j_4,j_1=0}^{p} C_{j_4j_1j_4j_1}=0,
$$

\vspace{2mm}
$$
\hbox{\vtop{\offinterlineskip\halign{
\hfil#\hfil\cr
{\rm lim}\cr
$\stackrel{}{{}_{p\to \infty}}$\cr
}} }
\sum\limits_{j_4,j_2=0}^{p} C_{j_4j_2j_2j_4}=0.
$$

\vspace{6mm}

Then

\vspace{-3mm}

$$
\hbox{\vtop{\offinterlineskip\halign{
\hfil#\hfil\cr
{\rm lim}\cr
$\stackrel{}{{}_{p\to \infty}}$\cr
}} }
R_{T,t}^{(7)pppp}=0.
$$

\vspace{6mm}

Theorem 6 is proved.

\vspace{5mm}

\section{Theorems 1--6 from Point
of View of the Wong--Zakai Approximation}

\vspace{5mm}

The iterated Ito stochastic integrals and solutions
of Ito SDEs are complex and important func\-ti\-onals
from the independent components ${\bf f}_{s}^{(i)},$
$i=1,\ldots,m$ of the multidimensional
Wiener process ${\bf f}_{s},$ $s\in[0, T].$
Let ${\bf f}_{s}^{(i)p},$ $p\in\mathbb{N}$ 
be some approximation of
${\bf f}_{s}^{(i)},$
$i=1,\ldots,m$.
Suppose that 
${\bf f}_{s}^{(i)p}$
converges to
${\bf f}_{s}^{(i)},$
$i=1,\ldots,m$ if $p\to\infty$ in some sense and has
differentiable sample trajectories.

A natural question arises: if we replace 
${\bf f}_{s}^{(i)}$
by ${\bf f}_{s}^{(i)p},$
$i=1,\ldots,m$ in the functionals
mentioned above, will the resulting
functionals converge to the original
functionals from the components 
${\bf f}_{s}^{(i)},$
$i=1,\ldots,m$ of the multidimentional
Wiener process ${\bf f}_{s}$?
The answere to this question is negative 
in the general case. However, 
in the pioneering works of Wong E. and Zakai M. \cite{W-Z-1},
\cite{W-Z-2},
it was shown that under the special conditions and 
for some types of approximations 
of the Wiener process the answere is affirmative
with one peculiarity: the convergence takes place 
to the iterated Stratonovich stochastic integrals
and solutions of Stratonovich SDEs and not to iterated 
Ito stochastic integrals and solutions
of Ito SDEs.
The piecewise 
linear approximation 
as well as the regularization by convolution 
\cite{W-Z-1}-\cite{Watanabe} relate the 
mentioned types of approximations
of the Wiener process. The above approximation 
of stochastic integrals and solutions of SDEs 
is often called the Wong--Zakai approximation.

Let ${\bf w}_{\tau},$ $\tau\in[0, T]$ is a random vector with 
an $m+1$ components: ${\bf w}_{\tau}^{(i)}={\bf f}_{\tau}^{(i)}$ 
for $i=1,\ldots,m$ and 
${\bf w}_{\tau}^{(0)}=\tau,$\ 
${\bf f}_{\tau}^{(i)}$ $(i=1,\ldots,m)$
are independent standard Wiener processes.

It is well known that the following representation 
takes place \cite{Lipt}, \cite{7e}

\begin{equation}
\label{um1x}
{\bf w}_{\tau}^{(i)}-{\bf w}_{t}^{(i)}=
\sum_{j=0}^{\infty}\int\limits_t^{\tau}
\phi_j(s)ds\ \zeta_j^{(i)},\ \ \ \zeta_j^{(i)}=
\int\limits_t^T \phi_j(s)d{\bf w}_s^{(i)},
\end{equation}

\vspace{4mm}
\noindent
where $\tau\in[t, T],$ $t\ge 0,$
$\{\phi_j(x)\}_{j=0}^{\infty}$ is an arbitrary complete 
orthonormal system of functions in the space $L_2([t, T]),$ and
$\zeta_j^{(i)}$ are independent standard Gaussian 
random variables for various $i$ or $j.$
Moreover, the series (\ref{um1x}) converges for any $\tau\in [t, T]$
in the mean-square sense.

Let ${\bf w}_{\tau}^{(i)p}-{\bf w}_{t}^{(i)p}$ be 
the mean-square approximation of the process
${\bf w}_{\tau}^{(i)}-{\bf w}_{t}^{(i)},$
which has the following form

\vspace{-3mm}
\begin{equation}
\label{um1xx}
{\bf w}_{\tau}^{(i)p}-{\bf w}_{t}^{(i)p}=
\sum_{j=0}^{p}\int\limits_t^{\tau}
\phi_j(s)ds\ \zeta_j^{(i)}.
\end{equation}

\vspace{3mm}

From (\ref{um1xx}) we obtain

\vspace{-4mm}
\begin{equation}
\label{um1xxx}
d{\bf w}_{\tau}^{(i)p}=
\sum_{j=0}^{p}
\phi_j(\tau)\zeta_j^{(i)} d\tau.
\end{equation}

\vspace{4mm}

Consider the following iterated Riemann--Stieltjes
integral

\begin{equation}
\label{um1xxxx}
\int\limits_t^T
\psi_k(t_k)\ldots \int\limits_t^{t_2}\psi_1(t_1)
d{\bf w}_{t_1}^{(i_1)p_1}\ldots d{\bf w}_{t_k}^{(i_k)p_k},
\end{equation}

\vspace{4mm}
\noindent
where $i_1,\ldots,i_k=0,1,\ldots,m,$\ \ $p_1,\ldots,p_k\in\mathbb{N},$ 

\begin{equation}
\label{um1xxx1}
d{\bf w}_{\tau}^{(i)p}=
\left\{\begin{matrix}
d{\bf f}_{\tau}^{(i)p}\ &\hbox{\rm for}\ \ \ i=1,\ldots,m\cr\cr\cr
d\tau^p\ &\hbox{\rm for}\ \ \ i=0
\end{matrix}
,\right.
\end{equation}

\vspace{4mm}
\noindent
and $d{\bf f}_{\tau}^{(i)p},$ $d\tau^p$ are defined by the relation (\ref{um1xxx}).

Let us substitute (\ref{um1xxx}) into (\ref{um1xxxx})

\begin{equation}
\label{um1xxxx1}
\int\limits_t^T
\psi_k(t_k)\ldots \int\limits_t^{t_2}\psi_1(t_1)
d{\bf w}_{t_1}^{(i_1)p_1}\ldots d{\bf w}_{t_k}^{(i_k)p_k}=
\sum\limits_{j_1=0}^{p_1} \ldots \sum\limits_{j_k=0}^{p_k}
C_{j_k \ldots j_1}\prod\limits_{l=1}^k \zeta_{j_l}^{(i_l)},
\end{equation}

\vspace{4mm}
\noindent
where 
$$
\zeta_j^{(i)}=\int\limits_t^T \phi_j(s)d{\bf w}_s^{(i)}
$$ 

\vspace{2mm}
\noindent
are independent standard Gaussian random variables for various 
$i$ or $j$ (in the case when $i\ne 0$),
${\bf w}_{s}^{(i)}={\bf f}_{s}^{(i)}$ for
$i=1,\ldots,m$ and 
${\bf w}_{s}^{(0)}=s,$

$$
C_{j_k \ldots j_1}=\int\limits_t^T\psi_k(t_k)\phi_{j_k}(t_k)\ldots
\int\limits_t^{t_2}
\psi_1(t_1)\phi_{j_1}(t_1)
dt_1\ldots dt_k
$$

\vspace{4mm}
\noindent
is the Fourier coefficient.

To best of our knowledge \cite{W-Z-1}-\cite{Watanabe}
the approximations of the Wiener process
in the Wong--Zakai approximation must satisfy fairly strong
restrictions
\cite{Watanabe}
(see Definition 7.1, pp.~480--481).
Moreover, approximations of the Wiener process that are
similar to (\ref{um1xx})
were not considered in \cite{W-Z-1}, \cite{W-Z-2}
(also see \cite{Watanabe}, Theorems 7.1, 7.2).
Therefore, the proof of analogs of Theorems 7.1 and 7.2 \cite{Watanabe}
for approximations of the Wiener 
process based on its series expansion (\ref{um1x})
should be carried out separately.
Thus, the mean-square convergence of the right-hand side
of (\ref{um1xxxx1}) to the iterated Stratonovich stochastic integral 
(\ref{str})
does not follow from the results of the papers
\cite{W-Z-1}, \cite{W-Z-2} (also see \cite{Watanabe},
Theorems 7.1, 7.2).

From the other hand, Theorems 1--6 and Theorems 7--10 (see below) from this 
paper can be considered as the proof of the
Wong--Zakai approximation for the iterated Stratonovich stochastic integrals
(\ref{str}) of multiplicities 1 to 6 based on the
Riemann--Stieltjes integrals (\ref{um1xxxx}) and
approximation (\ref{um1xx}) of the Wiener process.
At that, the mentioned Riemann--Stieltjes integrals converge
(according to Theorems 1--6 and Theorems 7--10 (see below))
to the appropriate Stratonovich 
stochastic integrals (\ref{str}). Recall that
$\{\phi_j(x)\}_{j=0}^{\infty}$ (see (\ref{um1x}), (\ref{um1xx}), and
Theorems 3--10)
is a complete 
orthonormal system of Legendre polynomials or 
trigonometric functions 
in the space $L_2([t, T])$.

To illustrate the above reasoning, 
consider two examples for the case $k=2,$
$\psi_1(s),$ $\psi_2(s)\equiv 1;$ $i_1, i_2=1,\ldots,m.$

The first example relates to the piecewise linear approximation
of the multidimensional Wiener process (these approximations 
were considered in \cite{W-Z-1}-\cite{Watanabe}).

Let ${\bf b}_{\Delta}^{(i)}(t),$ $t\in[0, T]$ be the piecewise
linear approximation of the $i$th component ${\bf f}_t^{(i)}$
of the multidimensional standard Wiener process ${\bf f}_t,$
$t\in [0, T]$ with independent components
${\bf f}_t^{(i)},$ $i=1,\ldots,m,$ i.e.

$$
{\bf b}_{\Delta}^{(i)}(t)={\bf f}_{k\Delta}^{(i)}+
\frac{t-k\Delta}{\Delta}\Delta{\bf f}_{k\Delta}^{(i)},
$$

\vspace{3mm}
\noindent
where 

\vspace{-2mm}
$$
\Delta{\bf f}_{k\Delta}^{(i)}={\bf f}_{(k+1)\Delta}^{(i)}-
{\bf f}_{k\Delta}^{(i)},\ \ \
t\in[k\Delta, (k+1)\Delta),\ \ \ k=0, 1,\ldots, N-1.
$$

\vspace{4mm}

Note that w.~p.~1

\vspace{-1mm}
\begin{equation}
\label{pridum}
\frac{d{\bf b}_{\Delta}^{(i)}}{dt}(t)=
\frac{\Delta{\bf f}_{k\Delta}^{(i)}}{\Delta},\ \ \
t\in[k\Delta, (k+1)\Delta),\ \ \ k=0, 1,\ldots, N-1.
\end{equation}

\vspace{4mm}

Consider the following iterated Riemann--Stieltjes
integral

\vspace{1mm}
$$
\int\limits_0^T
\int\limits_0^{s}
d{\bf b}_{\Delta}^{(i_1)}(\tau)d{\bf b}_{\Delta}^{(i_2)}(s),\ \ \ 
i_1,i_2=1,\ldots,m.
$$

\vspace{4mm}

Using (\ref{pridum}) and additive property of the Riemann--Stieltjes integral, 
we can write w.~p.~1

\vspace{2mm}
$$
\int\limits_0^T
\int\limits_0^{s}
d{\bf b}_{\Delta}^{(i_1)}(\tau)d{\bf b}_{\Delta}^{(i_2)}(s)=
\int\limits_0^T
\int\limits_0^{s}
\frac{d{\bf b}_{\Delta}^{(i_1)}}{d\tau}(\tau)d\tau
\frac{d {\bf b}_{\Delta}^{(i_2)}}{d s}(s)
ds =
$$

\vspace{3mm}
$$
=
\sum\limits_{l=0}^{N-1}\int\limits_{l\Delta}^{(l+1)\Delta}
\left(
\sum\limits_{q=0}^{l-1}\int\limits_{q\Delta}^{(q+1)\Delta}
\frac{\Delta{\bf f}_{q\Delta}^{(i_1)}}{\Delta}d\tau+
\int\limits_{l\Delta}^{s}
\frac{\Delta{\bf f}_{l\Delta}^{(i_1)}}{\Delta}d\tau\right)
\frac{\Delta{\bf f}_{l\Delta}^{(i_2)}}{\Delta}ds=
$$

\vspace{3mm}
$$
=\sum\limits_{l=0}^{N-1}\sum\limits_{q=0}^{l-1}
\Delta{\bf f}_{q\Delta}^{(i_1)}
\Delta{\bf f}_{l\Delta}^{(i_2)}+
\frac{1}{\Delta^2}\sum\limits_{l=0}^{N-1}
\Delta{\bf f}_{l\Delta}^{(i_1)}
\Delta{\bf f}_{l\Delta}^{(i_2)}
\int\limits_{l\Delta}^{(l+1)\Delta}
\int\limits_{l\Delta}^{s}d\tau ds=
$$

\vspace{3mm}
\begin{equation}
\label{oh-ty}
=\sum\limits_{l=0}^{N-1}\sum\limits_{q=0}^{l-1}
\Delta{\bf f}_{q\Delta}^{(i_1)}
\Delta{\bf f}_{l\Delta}^{(i_2)}+
\frac{1}{2}\sum\limits_{l=0}^{N-1}
\Delta{\bf f}_{l\Delta}^{(i_1)}
\Delta{\bf f}_{l\Delta}^{(i_2)}.
\end{equation}

\vspace{6mm}

Using (\ref{oh-ty}), it 
is not difficult to show 
that

\vspace{1mm}
$$
\hbox{\vtop{\offinterlineskip\halign{
\hfil#\hfil\cr
{\rm l.i.m.}\cr
$\stackrel{}{{}_{N\to \infty}}$\cr
}} }
\int\limits_0^T
\int\limits_0^{s}
d{\bf b}_{\Delta}^{(i_1)}(\tau)d{\bf b}_{\Delta}^{(i_2)}(s)=
\int\limits_0^T
\int\limits_0^{s}
d{\bf f}_{\tau}^{(i_1)}d{\bf f}_{s}^{(i_2)}+
\frac{1}{2}{\bf 1}_{\{i_1=i_2\}}\int\limits_0^T ds=
$$

\vspace{3mm}
\begin{equation}
\label{uh-111}
=
\int\limits_0^{*T}
\int\limits_0^{*s}
d{\bf f}_{\tau}^{(i_1)}d{\bf f}_{s}^{(i_2)},
\end{equation}

\vspace{5mm}
\noindent
where $\Delta\to 0$ if $N\to\infty$ ($N\Delta=T$).

Obviously, (\ref{uh-111}) agrees with Theorem 7.1 (see \cite{Watanabe},
p.~486).

The next example relates to the approximation
of the Wiener process based on its series expansion
(\ref{um1x}) for $t=0$, where
$\{\phi_j(x)\}_{j=0}^{\infty}$ 
is a complete 
orthonormal system of Legendre polynomials or 
trigonometric functions 
in the space $L_2([0, T])$.

Consider the following iterated Riemann--Stieltjes
integral

\vspace{-1mm}
\begin{equation}
\label{abcd1}
\int\limits_0^T
\int\limits_0^{s}
d{\bf f}_{\tau}^{(i_1)p}d{\bf f}_{s}^{(i_2)p},\ \ \ 
i_1,i_2=1,\ldots,m,
\end{equation}

\vspace{3mm}
\noindent
where $d{\bf f}_{\tau}^{(i)p}$ is defined by the
relation
(\ref{um1xxx}).

Let us substitute (\ref{um1xxx}) into (\ref{abcd1}) 

\vspace{-1mm}
\begin{equation}
\label{set18}
\int\limits_0^T
\int\limits_0^{s}
d{\bf f}_{\tau}^{(i_1)p}d{\bf f}_{s}^{(i_2)p}=
\sum\limits_{j_1,j_2=0}^p
C_{j_2 j_1} \zeta_{j_1}^{(i_1)}\zeta_{j_2}^{(i_2)},
\end{equation}

\vspace{3mm}
\noindent
where 
$$
C_{j_2 j_1}=
\int\limits_0^T \phi_{j_2}(s)\int\limits_0^s
\phi_{j_1}(\tau)d\tau ds
$$

\vspace{3mm}
\noindent
is the Fourier coefficient; another notations 
are the same as in (\ref{um1xxxx1}).

As we noted above, approximations of the Wiener process that are
similar to (\ref{um1xx})
were not considered in \cite{W-Z-1}, \cite{W-Z-2}
(also see Theorems 7.1, 7.2 in \cite{Watanabe}).
Furthermore, the extension of the results of Theorems 7.1 and 7.2
\cite{Watanabe} to the case under consideration is
not obvious.

On the other hand, we can apply the theory built in Chapters 1 and 2
of the monographs \cite{20xx}-\cite{20xxz}. More precisely, 
using 
Theorems 3, 4 from this paper, 
we obtain from (\ref{set18}) the desired result

\vspace{-1mm}
$$
\hbox{\vtop{\offinterlineskip\halign{
\hfil#\hfil\cr
{\rm l.i.m.}\cr
$\stackrel{}{{}_{p\to \infty}}$\cr
}} }
\int\limits_0^T
\int\limits_0^{s}
d{\bf f}_{\tau}^{(i_1)p}d{\bf f}_{s}^{(i_2)p}=
\hbox{\vtop{\offinterlineskip\halign{
\hfil#\hfil\cr
{\rm l.i.m.}\cr
$\stackrel{}{{}_{p\to \infty}}$\cr
}} }
\sum\limits_{j_1,j_2=0}^p
C_{j_2 j_1} \zeta_{j_1}^{(i_1)}\zeta_{j_2}^{(i_2)}=
$$

\vspace{2mm}
\begin{equation}
\label{umen-bl}
=
\int\limits_0^{*T}
\int\limits_0^{*s}
d{\bf f}_{\tau}^{(i_1)}d{\bf f}_{s}^{(i_2)}.
\end{equation}

\vspace{5mm}

From the other hand, by Theorems 1, 2
(see (\ref{a2})) for the case
$k=2$ we obtain from (\ref{set18}) the following relation

\vspace{-2mm}
$$
\hbox{\vtop{\offinterlineskip\halign{
\hfil#\hfil\cr
{\rm l.i.m.}\cr
$\stackrel{}{{}_{p\to \infty}}$\cr
}} }
\int\limits_0^T
\int\limits_0^{s}
d{\bf f}_{\tau}^{(i_1)p}d{\bf f}_{s}^{(i_2)p}=
\hbox{\vtop{\offinterlineskip\halign{
\hfil#\hfil\cr
{\rm l.i.m.}\cr
$\stackrel{}{{}_{p\to \infty}}$\cr
}} }
\sum\limits_{j_1,j_2=0}^p
C_{j_2 j_1} \zeta_{j_1}^{(i_1)}\zeta_{j_2}^{(i_2)}=
$$

\vspace{2mm}
$$
=
\hbox{\vtop{\offinterlineskip\halign{
\hfil#\hfil\cr
{\rm l.i.m.}\cr
$\stackrel{}{{}_{p\to \infty}}$\cr
}} }
\sum\limits_{j_1,j_2=0}^p
C_{j_2 j_1} \biggl(\zeta_{j_1}^{(i_1)}\zeta_{j_2}^{(i_2)}-
{\bf 1}_{\{i_1=i_2\}}{\bf 1}_{\{j_1=j_2\}}\biggr)+
{\bf 1}_{\{i_1=i_2\}}\sum\limits_{j_1=0}^{\infty}
C_{j_1 j_1}=
$$

\vspace{2mm}
\begin{equation}
\label{umen-blx}
=
\int\limits_0^T
\int\limits_0^{s}
d{\bf f}_{\tau}^{(i_1)}d{\bf f}_{s}^{(i_2)}+
{\bf 1}_{\{i_1=i_2\}}\sum\limits_{j_1=0}^{\infty}
C_{j_1 j_1}.
\end{equation}

\vspace{5mm}

Since
$$
\sum\limits_{j_1=0}^{\infty}
C_{j_1 j_1}=\frac{1}{2}\sum\limits_{j_1=0}^{\infty}
\left(\int\limits_0^T \phi_j(\tau)d\tau\right)^2
=\frac{1}{2}
\left(\int\limits_0^T \phi_0(\tau)d\tau\right)^2=\frac{1}{2}
\int\limits_0^T ds,
$$

\vspace{5mm}
\noindent
then from standard relation between 
Ito and Stratonovich stochastic integrals and
(\ref{umen-blx}) we obtain (\ref{umen-bl}).

\vspace{5mm}

\section{Recent Results on Expansions of Iterated Stratonovich
Stochastic Integrals of Multiplicities 3 to 6}

\vspace{5mm}

Recently, a new approach to the expansion and mean-square 
approximation of iterated Stratonovich stochastic integrals has been obtained
\cite{20xx} (Sect.~2.10--2.16), 
\cite{25} (Sect.~5--11), \cite{arxiv-11} (Sect.~7--13),
\cite{32} (Sect.~13--19), \cite{new-art-1-xxy} (Sect.~4--9), \cite{new-art-1xxys}.
Let us formulate four theorems that were obtained using this approach.

\vspace{2mm}

{\bf Theorem 7}\ \cite{20xx}, \cite{25}, \cite{arxiv-11}, \cite{32}, \cite{new-art-1-xxy}.\
{\it Suppose 
that $\{\phi_j(x)\}_{j=0}^{\infty}$ is a complete orthonormal system of 
Legendre polynomials or trigonometric functions in the space $L_2([t, T]).$
Furthermore, let $\psi_1(\tau), \psi_2(\tau),$ $\psi_3(\tau)$ are continuously dif\-ferentiable 
nonrandom functions on $[t, T].$ 
Then, for the 
iterated Stra\-to\-no\-vich stochastic integral of third multiplicity

$$
J^{*}[\psi^{(3)}]_{T,t}={\int\limits_t^{*}}^T\psi_3(t_3)
{\int\limits_t^{*}}^{t_3}\psi_2(t_2)
{\int\limits_t^{*}}^{t_2}\psi_1(t_1)
d{\bf w}_{t_1}^{(i_1)}
d{\bf w}_{t_2}^{(i_2)}d{\bf w}_{t_3}^{(i_3)}\ \ \ (i_1,i_2,i_3=0,1,\ldots,m)
$$

\vspace{4mm}
\noindent
the following 
relations

\vspace{-1mm}
\begin{equation}
\label{fin1}
J^{*}[\psi^{(3)}]_{T,t}
=\hbox{\vtop{\offinterlineskip\halign{
\hfil#\hfil\cr
{\rm l.i.m.}\cr
$\stackrel{}{{}_{p\to \infty}}$\cr
}} }
\sum\limits_{j_1, j_2, j_3=0}^{p}
C_{j_3 j_2 j_1}\zeta_{j_1}^{(i_1)}\zeta_{j_2}^{(i_2)}\zeta_{j_3}^{(i_3)},
\end{equation}

\vspace{3mm}
\begin{equation}
\label{fin2}
{\sf M}\left\{\left(
J^{*}[\psi^{(3)}]_{T,t}-
\sum\limits_{j_1, j_2, j_3=0}^{p}
C_{j_3 j_2 j_1}\zeta_{j_1}^{(i_1)}\zeta_{j_2}^{(i_2)}\zeta_{j_3}^{(i_3)}\right)^2\right\}
\le \frac{C}{p}
\end{equation}

\vspace{5mm}
\noindent
are fulfilled, where $i_1, i_2, i_3=0,1,\ldots,m$ in {\rm (\ref{fin1})} and 
$i_1, i_2, i_3=1,\ldots,m$ in {\rm (\ref{fin2})},
constant $C$ is independent of $p,$

$$
C_{j_3 j_2 j_1}=\int\limits_t^T\psi_3(t_3)\phi_{j_3}(t_3)
\int\limits_t^{t_3}\psi_2(t_2)\phi_{j_2}(t_2)
\int\limits_t^{t_2}\psi_1(t_1)\phi_{j_1}(t_1)dt_1dt_2dt_3
$$

\vspace{4mm}
\noindent
and
$$
\zeta_{j}^{(i)}=
\int\limits_t^T \phi_{j}(\tau) d{\bf f}_{\tau}^{(i)}
$$ 

\vspace{2mm}
\noindent
are independent standard Gaussian random variables for various 
$i$ or $j$ {\rm (}in the case when $i\ne 0${\rm );} 
another notations are the same as in Theorems~{\rm 1, 2}.}

\vspace{2mm}

{\bf Theorem 8}\ \cite{20xx}, \cite{25}, \cite{arxiv-11}, \cite{32}, \cite{new-art-1-xxy}.\ 
{\it Let
$\{\phi_j(x)\}_{j=0}^{\infty}$ be a complete orthonormal system of 
Legendre polynomials or trigonometric functions in the space $L_2([t, T]).$
Furthermore, let $\psi_1(\tau), \ldots, \psi_4(\tau)$ be continuously dif\-ferentiable 
nonrandom functions on $[t, T].$ 
Then, for the 
iterated Stra\-to\-no\-vich stochastic integral of fourth multiplicity

\begin{equation}
\label{fin0}
J^{*}[\psi^{(4)}]_{T,t}={\int\limits_t^{*}}^T\psi_4(t_4)
{\int\limits_t^{*}}^{t_4}\psi_3(t_3)
{\int\limits_t^{*}}^{t_3}\psi_2(t_2)
{\int\limits_t^{*}}^{t_2}\psi_1(t_1)
d{\bf w}_{t_1}^{(i_1)}
d{\bf w}_{t_2}^{(i_2)}d{\bf w}_{t_3}^{(i_3)}d{\bf w}_{t_4}^{(i_4)}
\end{equation}

\vspace{4mm}
\noindent
the following 
relations

\begin{equation}
\label{fin3}
J^{*}[\psi^{(4)}]_{T,t}
=\hbox{\vtop{\offinterlineskip\halign{
\hfil#\hfil\cr
{\rm l.i.m.}\cr
$\stackrel{}{{}_{p\to \infty}}$\cr
}} }
\sum\limits_{j_1, j_2, j_3,j_4=0}^{p}
C_{j_4j_3 j_2 j_1}\zeta_{j_1}^{(i_1)}\zeta_{j_2}^{(i_2)}\zeta_{j_3}^{(i_3)}\zeta_{j_4}^{(i_4)},
\end{equation}

\vspace{3mm}

\begin{equation}
\label{fin4}
{\sf M}\left\{\left(
J^{*}[\psi^{(4)}]_{T,t}-
\sum\limits_{j_1, j_2, j_3, j_4=0}^{p}
C_{j_4 j_3 j_2 j_1}\zeta_{j_1}^{(i_1)}\zeta_{j_2}^{(i_2)}\zeta_{j_3}^{(i_3)}
\zeta_{j_4}^{(i_4)}
\right)^2\right\}
\le \frac{C}{p^{1-\varepsilon}}
\end{equation}

\vspace{5mm}
\noindent
are fulfilled, where $i_1, \ldots , i_4=0,1,\ldots,m$ in {\rm (\ref{fin0}),} {\rm (\ref{fin3})} 
and $i_1, \ldots, i_4=1,\ldots,m$ in {\rm (\ref{fin4}),}
constant $C$ does not depend on $p,$
$\varepsilon$ is an arbitrary
small positive real number 
for the case of complete orthonormal system of 
Legendre polynomials in the space $L_2([t, T])$
and $\varepsilon=0$ for the case of
complete orthonormal system of 
trigonometric functions in the space $L_2([t, T]),$

\vspace{1mm}
$$
C_{j_4 j_3 j_2 j_1}=
$$

$$
=
\int\limits_t^T\psi_4(t_4)\phi_{j_4}(t_4)
\int\limits_t^{t_4}\psi_3(t_3)\phi_{j_3}(t_3)
\int\limits_t^{t_3}\psi_2(t_2)\phi_{j_2}(t_2)
\int\limits_t^{t_2}\psi_1(t_1)\phi_{j_1}(t_1)dt_1dt_2dt_3dt_4;
$$

\vspace{5mm}
\noindent
another notations are the same as in Theorem~{\rm 7}.}

\vspace{2mm}

{\bf Theorem 9}\ \cite{20xx}, \cite{25}, \cite{arxiv-11}, \cite{32}, \cite{new-art-1-xxy}.\
{\it Assume 
that $\{\phi_j(x)\}_{j=0}^{\infty}$ is a complete orthonormal system of 
Legendre polynomials or trigonometric functions in the space $L_2([t, T])$
and $\psi_1(\tau), \ldots, \psi_5(\tau)$ are continuously dif\-ferentiable 
nonrandom functions on $[t, T].$ 
Then, for the 
iterated Stra\-to\-no\-vich stochastic integral of fifth multiplicity

\begin{equation}
\label{fin7}
J^{*}[\psi^{(5)}]_{T,t}={\int\limits_t^{*}}^T\psi_5(t_5)
\ldots
{\int\limits_t^{*}}^{t_2}\psi_1(t_1)
d{\bf w}_{t_1}^{(i_1)}
\ldots d{\bf w}_{t_5}^{(i_5)}
\end{equation}

\vspace{4mm}
\noindent
the following 
relations

\begin{equation}
\label{fin8}
J^{*}[\psi^{(5)}]_{T,t}
=\hbox{\vtop{\offinterlineskip\halign{
\hfil#\hfil\cr
{\rm l.i.m.}\cr
$\stackrel{}{{}_{p\to \infty}}$\cr
}} }
\sum\limits_{j_1,\ldots,j_5=0}^{p}
C_{j_5 \ldots j_1}\zeta_{j_1}^{(i_1)}\ldots \zeta_{j_5}^{(i_5)},
\end{equation}

\vspace{3mm}

\begin{equation}
\label{fin9}
{\sf M}\left\{\left(
J^{*}[\psi^{(5)}]_{T,t}-
\sum\limits_{j_1, \ldots, j_5=0}^{p}
C_{j_5 \ldots j_1}\zeta_{j_1}^{(i_1)}\ldots
\zeta_{j_5}^{(i_5)}
\right)^2\right\}
\le \frac{C}{p^{1-\varepsilon}}
\end{equation}

\vspace{5mm}
\noindent
are fulfilled, where $i_1, \ldots , i_5=0,1,\ldots,m$ in {\rm (\ref{fin7}),} {\rm (\ref{fin8})} 
and $i_1, \ldots, i_5=1,\ldots,m$ in {\rm (\ref{fin9}),}
constant $C$ is independent of $p,$
$\varepsilon$ is an arbitrary
small positive real number 
for the case of complete orthonormal system of 
Legendre polynomials in the space $L_2([t, T])$
and $\varepsilon=0$ for the case of
complete orthonormal system of 
trigonometric functions in the space $L_2([t, T]),$

\vspace{1mm}
$$
C_{j_5 \ldots j_1}=
\int\limits_t^T\psi_5(t_5)\phi_{j_5}(t_5)\ldots
\int\limits_t^{t_2}\psi_1(t_1)\phi_{j_1}(t_1)dt_1\ldots dt_5;
$$

\vspace{5mm}
\noindent
another notations are the same as in Theorems~{\rm 7, 8}.}

\vspace{2mm}

{\bf Theorem 10}\ \cite{20xx}, \cite{25}, \cite{arxiv-11}, \cite{32}, \cite{new-art-1xxys}.\
{\it Suppose that 
$\{\phi_j(x)\}_{j=0}^{\infty}$ is a complete orthonormal system of 
Legendre polynomials or trigonometric functions in the space $L_2([t, T]).$
Then, for the 
iterated Stratonovich stochastic integral of sixth multiplicity

\begin{equation}
\label{after10001qu1}
J_{T,t}^{*(i_1\ldots i_6)}={\int\limits_t^{*}}^T
\ldots
{\int\limits_t^{*}}^{t_2}
d{\bf w}_{t_1}^{(i_1)}
\ldots d{\bf w}_{t_6}^{(i_6)}
\end{equation}

\vspace{3mm}
\noindent
the following 
expansion 

\vspace{-1mm}
$$
J_{T,t}^{*(i_1\ldots i_6)}
=\hbox{\vtop{\offinterlineskip\halign{
\hfil#\hfil\cr
{\rm l.i.m.}\cr
$\stackrel{}{{}_{p\to \infty}}$\cr
}} }
\sum\limits_{j_1, \ldots, j_6=0}^{p}
C_{j_6 \ldots j_1}\zeta_{j_1}^{(i_1)}\ldots
\zeta_{j_6}^{(i_6)}
$$

\vspace{4mm}
\noindent
that converges in the mean-square sense is valid, where
$i_1, \ldots, i_6=0, 1,\ldots,m,$

$$
C_{j_6 \ldots j_1}=
\int\limits_t^T\phi_{j_6}(t_6)\ldots
\int\limits_t^{t_2}\phi_{j_1}(t_1)dt_1\ldots dt_6;
$$

\vspace{3mm}
\noindent
another notations are the same as in Theorems~{\rm 7--9}.}

\vspace{2mm}

The results of Theorems~7--10 section were developed in 
\cite{20xx} (Chapter~2), \cite{25}, \cite{arxiv-11}, \cite{32}, \cite{2024xx}-\cite{2025xxxaaa}.
In particular, analogues of Theorem~10 for iterated Stratonovich stochastic
integrals of multiplicities 7 and 8 were obtained in \cite{20xx} (Sect.~2.36, 2.37).
In addition, the variants of Thorems 7--9 
were obtained
for the case when $\{\phi_j(x)\}_{j=0}^{\infty}$ is an arbitrary complete orthonormal system
of functions in $L_2([t, T])$ \cite{20xx} (Sect.~2.1.4, 2.23, 2.24, 2.31--2.34),
\cite{25}, \cite{arxiv-11}, \cite{32}, \cite{2024xx}-\cite{2025xxxaaa}.

\end{document}